\documentclass{amsart}
\usepackage{amsmath} 
\usepackage{amssymb}

\newtheorem{theorem}{Theorem}[section] 
\newtheorem{claim}{Claim}[theorem]
\newtheorem{lemma}[theorem]{Lemma} 
\newtheorem{proposition}[theorem]{Proposition} 
\newtheorem{observation}[theorem]{Observation} 
\newtheorem{corollary}[theorem]{Corollary} 

\theoremstyle{definition}
\newtheorem{definition}[theorem]{Definition}
\newtheorem{example}[theorem]{Example}
\newtheorem{problem}[theorem]{Problem}

\theoremstyle{remark}
\newtheorem{remark}[theorem]{Remark}

\newtheorem{context}[theorem]{Context}

\numberwithin{equation}{section}

\newcommand{\conc}{{}^\frown\!}
\newcommand{\lh}{{\rm lh}\/}
\newcommand{\rest}{{\restriction}}
\newcommand{\vtl}{{\vartriangleleft}}
\newcommand{\dom}{{\rm dom}} 
\newcommand{\rk}{{\rm rk}} 
\newcommand{\rng}{{\rm rng}}
 
\newcommand{\suc}{{\rm succ}}
\newcommand{\mrot}{{\rm root}}
\newcommand{\set}{{\rm set}}
 
\newcommand{\gb}{{\mathfrak b}}

\newcommand{\bbC}{{\mathbb C}}
\newcommand{\gc}{{\mathfrak c}}

\newcommand{\gd}{{\mathfrak d\/}} 
\newcommand{\bH}{{\bf H}}
\newcommand{\cH}{{\mathcal H}}

\newcommand{\cF}{{\mathcal F}}

\newcommand{\cI}{{\mathcal I}}

\newcommand{\cL}{{\mathcal L}}

\newcommand{\bbP}{{\mathbb P}}
\newcommand{\cP}{{\mathcal P}}

\newcommand{\bp}{{\mathbf p}}
\newcommand{\bbQ}{{\mathbb Q}}

\newcommand{\cT}{{\mathcal T}}
\newcommand{\cU}{{\mathcal U}}
\newcommand{\cX}{{\mathcal X}}

\newcommand{\gY}{{\mathfrak Y}}

\newcommand{\cf}{{\rm cf}\/} 
 
\newcommand{\otp}{{\rm otp}\/} 
\newcommand{\st}{{\bf st}} 

\newcommand{\vare}{\varepsilon}
 
\newcommand{\Dl}{{{\mathcal D}\ell}}

\newcommand{\forces}{\Vdash} 
\newcommand{\bV}{{\bf V}} 
\newcommand{\lesdot}{\mathrel{\mathord{<}\!\!\raise 
0.8 pt\hbox{$\scriptstyle\circ$}}} 
\newcommand{\tegamexi}{{\Game^{{\rm pr}}_{\name{\cU}_\xi,\bp,\bar{\mu}}}} 
\newcommand{\tezero}{{\Game^{\pr}_{\cU,\bp,\bar{\mu}}}}
\newcommand{\teone}{{\Game^{\pr}_{E,\bp,\bar{\mu}}}}
\newcommand{\bgame}{{\Game^{\rm rbB}_{\bp}}}
\newcommand{\supgame}{{\Game^{B+}_{\bp}}}
\newcommand{\sergame}{{\Game^{\rm servant}_{S,D}}}
\newcommand{\masgame}{{\Game^{\rm master}_{S,D}}}
\newcommand{\dsgame}{{\Game^{\bf 2ser}_{S,D,\bar{\mu}}}}
\newcommand{\dmgame}{{\Game^{\bf 2mas}_{S,D,\bar{\mu}}}}

\newcommand{\fEE}{{\mathbb Q}^{\bar{E}}_E}
\newcommand{\pEE}{{\mathbb P}^{\bar{E}}_E}
\newcommand{\qFF}{{\mathbb Q}^{\bH}_{\bar{F},F}}
\newcommand{\pFF}{{\mathbb P}^{\bH}_{\bar{F},F}}
\newcommand{\rcbgame}{{\Game^{\rm rcB}_D}}

\newcommand{\tefo}{{\bbQ^2_D}}
\newcommand{\pos}{{\rm pos}}
\newcommand{\pr}{{\rm pr}}

\newcount\skewfactor
\def\mathunderaccent#1#2 {\let\theaccent#1\skewfactor#2
\mathpalette\putaccentunder}
\def\putaccentunder#1#2{\oalign{$#1#2$\crcr\hidewidth
\vbox to.2ex{\hbox{$#1\skew\skewfactor\theaccent{}$}\vss}\hidewidth}}
\def\name{\mathunderaccent\tilde-3 }

\begin{document}

\title{Lords of the iteration}

\author{Andrzej Ros{\l}anowski}
\address{Department of Mathematics\\
 University of Nebraska at Omaha\\
 Omaha, NE 68182-0243, USA}
\email{roslanow@member.ams.org}
\urladdr{http://www.unomaha.edu/logic}
\thanks{The first author would like to thank the Hebrew University of
  Jerusalem and the Lady Davis Fellowship Trust for awarding him with
  {\em Sch\"onbrunn Visiting Professorship\/} under which this research was
  carried out.}

\author{Saharon Shelah}
\address{Einstein Institute of Mathematics\\
Edmond J. Safra Campus, Givat Ram\\
The Hebrew University of Jerusalem\\
Jerusalem, 91904, Israel\\
 and  Department of Mathematics\\
 Rutgers University\\
 New Brunswick, NJ 08854, USA}
\email{shelah@math.huji.ac.il}
\urladdr{http://shelah.logic.at}
\thanks{Both authors acknowledge support from the United States-Israel
Binational Science Foundation (Grant no. 2002323). This is publication
888 of the second author.}

\subjclass{Primary 03E40; Secondary:03E35}
\date{July 2010}
\keywords{iterated forcing, $\lambda$--support, iteration theorems,}

\begin{abstract}
We introduce several properties of forcing notions which imply that their
$\lambda$--support iterations are $\lambda$--proper. Our methods and
techniques refine those studied in \cite{RoSh:655}, \cite{RoSh:777},
\cite{RoSh:860} and \cite{RoSh:890}, covering some new forcing notions
(though the exact relation of the new properties to the old ones remains
undecided).
\end{abstract}

\maketitle

\section{Introduction}
Since the beginning of 1980s it has been known that the theory of proper
forcing does not admit naive generalization to the context of larger
cardinals and iterations with larger supports. The evidence of that was
given already in Shelah \cite{Sh:b} (see \cite[Appendix 3.6(2)]{Sh:f}). It
seems that the first steps towards developing the theory of forcing iterated
with uncountable supports were done in Shelah \cite{Sh:587}, \cite{Sh:667},
but the properties introduced there were aimed at situations when we do not
want to add new subsets of $\lambda$ (corresponding to the case of {\em no
  new reals in CS iterations of proper forcing notions\/}). Later Ros{\l}anowski
and Shelah \cite{RoSh:655} introduced an iterable property called {\em
  properness over semi-diamonds\/} and then Eisworth \cite{Ei03} proposed
an iterable relative of it. These properties work nicely for
$\lambda$--support iterations (where $\lambda=\lambda^{<\lambda}$ is
essentially arbitrary) and forcings adding new subsets of $\lambda$, but the
price to pay is that many natural examples are not covered. If we restrict
ourselves to inaccessible $\lambda$, then the properties given by
Ros{\l}anowski and Shelah \cite{RoSh:777,RoSh:860,RoSh:890} may occur
useful. Those papers give both iteration theorems and new examples of
forcing notions for which the theorems apply. 

In the present paper we further advance the theory and we give results
applicable to both the case of inaccessible $\lambda$ as well as those
working for successor cardinals. The tools developed here may be treated as
yet another step towards {\em comparing and contrasting the structure of
${}^\lambda\lambda$ with that of ${}^\omega\omega$}.  That line of research
already has received some attention in the literature (see e.g., Cummings
and Shelah \cite{CuSh:541}, Shelah and Spasojevi{\' c} \cite{ShSj:643} or
Zapletal \cite{Za97}). Also with better iteration theorems one may hope for
further generalizations of Ros{\l}anowski and Shelah \cite{RoSh:470} to the
context of uncountable cardinals. (Initial steps in the latter direction
were presented in Ros{\l}anowski and Shelah \cite{RoSh:777}.) However, while
we do give some examples of forcing notions to which our properties apply,
we concentrate on the development of the theory of forcing leaving the {\em 
  real\/} applications for further investigations. The need for the
development of such general theory was indirectly stated by Hytinnen and
Rautila in \cite{HyRa01}, where they commented:
\begin{quote}
{\em Our proof is longer than the one in \cite{MkSh:398} partly because we are
not able to utilize the general theory of proper forcing, especially the
iteration lemma, but we have to prove everything ``from scratch''. } 
\end{quote} 
We believe that the present paper brings us substantially closer to the
right general iteration theorems for iterations with uncountable supports. 
\medskip

In the first section we introduce {\em $\Dl$--parameters\/} (which will play
an important role in our definitions) and a slight generalization of the
B--bounding property from \cite{RoSh:860}.  We also define a canonical
example for testing usefulness of our iteration theorems: the forcing $\fEE$
in which conditions are complete $\lambda$--trees in which along each
$\lambda$-branch the set of splittings forms a set from a filter $E$ (and
the splitting at $\nu$ is into a set from a filter $E_\nu$ on
$\lambda$). The main result of the first section (Theorem \ref{like860})
says that we may iterate with $\lambda$--supports forcing notions $\fEE$,
provided $\lambda$ is inaccessible and $E$ is always the same and has some
additional properties. 

If we want to iterate forcing notions like $\fEE$ but with different $E$ on
each coordinate (when the result of the first section is not applicable), we
may decide to use very orthogonal filters. Section 2 presents an iteration
theorem \ref{lordsthm} which is tailored for such situation. Also here we
need the assumption that $\lambda$ is inaccessible. 

The following section introduces {\em B--noble forcing notions\/} and the
iteration theorem \ref{nobleiteration} for them. The main gain here is that
it allows us to iterate (with $\lambda$--supports) forcing notions like
$\fEE$ even if $\lambda$ is not inaccessible. The fourth section gives more
examples of forcing notions and shows a possible application. In Corollary
\ref{improve860} we substantially improve a result from \cite{RoSh:860}
showing that dominating numbers associated with different filters may be
distinct even if $\lambda$ is a successor. 

The fifth section shows that some of closely related forcing notions may
have different properties. Section 6 presents yet another property that is
useful in $\lambda$--support iterations (for inaccessible $\lambda$): {\em
  reasonably merry forcing notions}. This property has the flavour of
putting together {\em being B--bounding\/} (of \cite{RoSh:860}) with {\em
  being fuzzy proper\/} (of \cite{RoSh:777}). We also give an example of a
forcing notion which is reasonably merry but which was not covered by
earlier properties.  We conclude the paper with a section listing open
problems.
\medskip

This research is a natural continuation of papers mentioned earlier
(\cite{Sh:587}, \cite{Sh:667}, \cite{RoSh:655}, \cite{RoSh:777},
\cite{RoSh:860} and \cite{RoSh:890}). All our iteration proofs are based on
{\em trees of conditions\/} and the arguments are similar to those from the
earlier works. While we tried to make this presentation self-contained, the
reader familiar with the previous papers will definitely find the proofs
presented here easier to follow (as several technical aspects do re-occur).

\subsection{Notation}
Our notation is rather standard and compatible with that of classical
textbooks (like Jech \cite{J}). In forcing we keep the older convention that
{\em a stronger condition is the larger one}. 

\begin{enumerate}
\item Ordinal numbers will be denoted be the lower case initial letters of
the Greek alphabet ($\alpha,\beta,\gamma,\delta\ldots$) and also by $i,j$
(with possible sub- and superscripts). 

Cardinal numbers will be called $\kappa,\lambda,\mu$; {\bf $\lambda$ will
  be always assumed to be a regular uncountable cardinal such that
  $\lambda^{<\lambda}=\lambda$\/} (we may forget to mention this).

Also, $\chi$ will denote a {\em sufficiently large\/} regular cardinal;
$\cH(\chi)$ is the family of all sets hereditarily of size less than
$\chi$. Moreover, we fix a well ordering $<^*_\chi$ of $\cH(\chi)$.

\item We will consider several games of two players. One player will be
called {\em Generic\/} or {\em Complete\/} or just {\em COM\/}, and we
will refer to this player as ``she''. Her opponent will be called {\em
Antigeneric\/} or {\em Incomplete} or just {\em INC\/} and will be
referred to as ``he''.

\item For a forcing notion $\bbP$, almost all $\bbP$--names for objects in
  the extension via $\bbP$ will be denoted with a tilde below (e.g.,
  $\name{\tau}$, $\name{X}$). There will be some exceptions to this rule,
  however.  $\Gamma_\bbP$ will stand for the canonical $\bbP$--name for the
  generic filter in $\bbP$. Also some (names for) normal filters generated
  in the extension from objects in the ground model will be denoted by $D$, 
  $D^{\bbP}$ or $D[\bbP]$.  

  The weakest element of $\bbP$ will be denoted by $\emptyset_\bbP$ (and we
  will always assume that there is one, and that there is no other condition
  equivalent to it). All forcing notions under considerations are assumed to
  be atomless. 

  By ``$\lambda$--support iterations'' we mean iterations in which domains
  of conditions are of size $\leq\lambda$. However, on some occasions we
  will pretend that conditions in a $\lambda$--support iteration
  $\bar{\bbQ}=\langle\bbP_\zeta,\name{\bbQ}_\zeta:\zeta< \zeta^* \rangle$
  are total functions on $\zeta^*$ and for $p\in\lim(\bar{\bbQ})$ and
  $\alpha\in\zeta^*\setminus\dom(p)$ we will let
  $p(\alpha)=\name{\emptyset}_{\name{\bbQ}_\alpha}$.

\item By ``a sequence'' we mean ``a function defined on a set of ordinals''
  (so the domain of a sequence does not have to be an ordinal). For two
  sequences $\eta,\nu$ we write $\nu\vtl\eta$ whenever $\nu$ is a proper
  initial segment of $\eta$, and $\nu \trianglelefteq\eta$ when either
  $\nu\vtl\eta$ or $\nu=\eta$.  The length of a sequence $\eta$ is the order
  type of its domain and it is denoted by $\lh(\eta)$. 

\item  A tree is a $\vtl$--downward closed set of sequences. A complete
  $\lambda$--tree is a tree $T\subseteq {}^{<\lambda}\lambda$ such that
  every $\vtl$-chain of size less than $\lambda$ has an $\vtl$-bound in $T$
  and for each $\eta\in T$ there is $\nu\in T$ such that $\eta\vtl\nu$.

Let $T$ be a $\lambda$--tree. For $\eta\in T$ we let 
\[\suc_T(\eta)=\{\alpha<\lambda:\eta\conc\langle\alpha\rangle\in T\}\quad
\mbox{ and }\quad (T)_\eta=\{\nu\in T:\nu\vtl\eta\mbox{ or }\eta
\trianglelefteq \nu\}.\]
We also let $\mrot(T)$ be the shortest $\eta\in T$ such that
$|\suc_T(\eta)|>1$ and $\lim_\lambda(T)=\{\eta\in{}^\lambda\lambda: (\forall
\alpha<\lambda)(\eta\rest\alpha\in T)\}$.
\end{enumerate}

\subsection{Background on trees of conditions}

\begin{definition}
\label{comp}
Let $\bbP$ be a forcing notion.
\begin{enumerate}
\item For a condition $r\in\bbP$, let $\Game_0^\lambda(\bbP,r)$ be the
following game of two players, {\em Complete} and  {\em Incomplete}:   
\begin{quotation}
\noindent the game lasts at most $\lambda$ moves and during a play the
players construct a sequence $\langle (p_i,q_i): i<\lambda\rangle$ of pairs
of conditions from $\bbP$ in such a way that $(\forall j<i<\lambda)(r\leq
p_j\leq q_j\leq p_i)$ and at the stage $i<\lambda$ of the game, first 
Incomplete chooses $p_i$ and then Complete chooses $q_i$.  
\end{quotation}
Complete wins if and only if for every $i<\lambda$ there are legal moves for
both players. 
\item We say that the forcing notion $\bbP$ is {\em strategically
$({<}\lambda)$--complete\/} if Complete has a winning strategy in the game 
$\Game_0^\lambda(\bbP,r)$ for each condition $r\in\bbP$. 
\item Let $N\prec (\cH(\chi),\in,<^*_\chi)$ be a model such that
${}^{<\lambda} N\subseteq N$, $|N|=\lambda$ and $\bbP\in N$. We say that a
condition $p\in\bbP$ is {\em $(N,\bbP)$--generic in the standard
sense\/} (or just: {\em $(N,\bbP)$--generic\/}) if for every
$\bbP$--name $\name{\tau}\in N$ for an ordinal we have $p\forces$``
$\name{\tau}\in N$ ''. 
\item $\bbP$ is {\em $\lambda$--proper in the standard sense\/} (or just:
{\em $\lambda$--proper\/}) if there is $x\in \cH(\chi)$  such that for
every model $N\prec (\cH(\chi),\in,<^*_\chi)$ satisfying  
\[{}^{<\lambda} N\subseteq N,\quad |N|=\lambda\quad\mbox{ and }\quad\bbP,x
  \in N, \]
and every condition $q\in N\cap\bbP$ there is an $(N,\bbP)$--generic
condition $p\in\bbP$ stronger than $q$.
\end{enumerate}
\end{definition}

\begin{remark}
Let us recall that if $\bbP$ is either strategically
$({<}\lambda^+)$--complete or $\lambda^+$--cc, then $\bbP$ is
$\lambda$--proper. Also, if $\bbP$ is $\lambda$--proper then
\begin{itemize}
\item $\lambda^+$ is not collapsed in forcing by $\bbP$, moreover
\item for every set of ordinals $A\in\bV^\bbP$ of size $\lambda$ there is a
  set $A^+\in\bV$ of size $\lambda$ such that $A\subseteq A^+$.  
\end{itemize}
\end{remark}

\begin{definition}
[Compare {\cite[Def. A.1.7]{RoSh:777}}, see also {\cite[Def. 2.2]{RoSh:860}}]
\label{dA.5}
\qquad
\begin{enumerate}
\item Let $\gamma$ be an ordinal, $\emptyset\neq w\subseteq \gamma$.
{\em A $(w,1)^\gamma$--tree\/} is a pair $\cT=(T,\rk)$ such that  
\begin{itemize}
\item $\rk:T\longrightarrow w\cup\{\gamma\}$, 
\item if $t\in T$ and $\rk(t)=\vare$, then $t$ is a sequence $\langle
(t)_\zeta: \zeta\in w\cap\vare\rangle$, 
\item $(T,\vtl)$ is a tree with root $\langle\rangle$ and 
\item if $t\in T$, then there is $t'\in T$ such that $t\trianglelefteq t'$
  and $\rk(t')=\gamma$.
\end{itemize}
\item If, additionally, $\cT=(T,\rk)$ is such that every chain in $T$ has a
  $\vtl$--upper bound in $T$, we will call it  {\em a standard
    $(w,1)^\gamma$--tree\/}  

We will keep the convention that $\cT^x_y$ is $(T^x_y,\rk^x_y)$.
\item Let $\bar{\bbQ}=\langle\bbP_i,\name{\bbQ}_i:i<\gamma\rangle$ be a
$\lambda$--support iteration. {\em A tree of conditions in $\bar{\bbQ}$} is
a system $\bar{p}=\langle p_t:t\in T\rangle$ such that  
\begin{itemize}
\item $(T,\rk)$ is a $(w,1)^\gamma$--tree for some $w\subseteq
\gamma$, 
\item $p_t\in\bbP_{\rk(t)}$ for $t\in T$, and
\item if $s,t\in T$, $s\vtl t$, then $p_s=p_t\rest\rk(s)$. 
\end{itemize}
If, additionally, $(T,\rk)$ is a standard tree, then $\bar{p}$ is called
{\em a standard tree of conditions}.  
\item Let $\bar{p}^0,\bar{p}^1$ be trees of conditions in $\bar{\bbQ}$,
  $\bar{p}^i=\langle p^i_t:t\in T\rangle$. We write $\bar{p}^0\leq
  \bar{p}^1$ whenever for each $t\in T$ we have $p^0_t\leq p^1_t$.   
\end{enumerate}
\end{definition}

Note that our standard trees and trees of conditions are a special case of
that \cite[Def. A.1.7]{RoSh:777} when $\alpha=1$. 

\begin{proposition}
\label{oldstuff}
Assume that $\bar{\bbQ}=\langle\bbP_i,\name{\bbQ}_i:i<\gamma\rangle$ is a
$\lambda$--support iteration such that for all $i<\gamma$ we have 
\[\forces_{\bbP_i}\mbox{`` $\name{\bbQ}_i$ is strategically
$({<}\lambda)$--complete ''.}\]
Suppose that $\bar{p}=\langle p_t:t\in T\rangle$ is a tree of conditions in
$\bar{\bbQ}$, $|T|<\lambda$, and $\cI\subseteq\bbP_\gamma$ is open
dense. Then there is a tree of conditions $\bar{q}=\langle q_t:t\in
T\rangle$ such that $\bar{p}\leq \bar{q}$ and $(\forall t\in T)(\rk(t)=
\gamma\ \Rightarrow\ q_t\in\cI)$. 
\end{proposition}

\begin{proof}
  This is essentially \cite[Proposition A.1.9]{RoSh:777} and the proof there
  applies here without changes.
\end{proof}

\section{$\Dl$--parameters}
In this section we introduce $\Dl$--parameters and we use them to get a
possible slight improvement of \cite[Theorem 3.1]{RoSh:860} (in Theorem
\ref{like860}). We also define our canonical testing forcing $\fEE$ to which
this result can be applied.

\begin{definition}
\label{Dlpar}
\begin{enumerate}
\item {\em A pre-$\Dl$--parameter on $\lambda$} is a triple
  $\bp=(\bar{P},S,D)=(\bar{P}^\bp,S^\bp,D^\bp)$ such that 
\begin{itemize}
\item $D$ is a proper uniform normal filter on $\lambda$, $S\in D$, 
\item $\bar{P}=\langle P_\delta:\delta\in S\rangle$ and $P_\delta\in
  [{}^\delta\delta]^{<\lambda}$ for each $\delta\in S$.
\end{itemize}
\item For a function $f\in {}^\lambda\lambda$ and a pre-$\Dl$--parameter
  $\bp=(\bar{P},S,D)$ we let 
\[\set^\bp(f)=\{\delta\in S:f\rest\delta\in P_\delta\}.\]
\item We say that a pre-$\Dl$--parameter $\bp=(\bar{P},S,D)$ is {\em a
  $\Dl$--parameter on $\lambda$} if $\set^\bp(f)\in D$ for every $f\in
  {}^\lambda\lambda$.  
\end{enumerate}
\end{definition}

\begin{example}
\begin{enumerate}
\item If $\lambda$ is strongly inaccessible, $D$ is the filter generated by
  club subsets of $\lambda$ and $P_\delta={}^\delta\delta$, $\bar{P}=\langle
  P_\delta:\delta<\lambda\rangle$, then $(\bar{P},\lambda,D)$ is a
  $\Dl$--parameter on $\lambda$. 
\item $\diamondsuit^+_\lambda$ is a statement asserting existence of a
  $\Dl$--parameter with the filter generated by clubs of $\lambda$. 
\item $\diamondsuit_\lambda$ implies the existence of a $\Dl$--parameter
  $(\bar{P},S,D)$ such that $|P_\delta|=|\delta|$.
\item For more instances of the existence of $\Dl$--parameters we refer the
  reader to Shelah \cite[\S 3]{Sh:460}. 
\end{enumerate}
\end{example}

\begin{definition}
Let $\bp$ be a pre-$\Dl$--parameter on $\lambda$ and $\bbQ$ be a forcing
notion not collapsing $\lambda$. In $\bV^\bbQ$ we define 
\begin{itemize}
\item $D^\bp[\bbQ]=D^{\bp[\bbQ]}$ is the normal filter generated by $D^\bp
  \cup\{\set^\bp(f):f\in {}^\lambda\lambda\}$,
\item $\bp[\bbQ]=(\bar{P}^\bp,S^\bp,D^{\bp[\bbQ]})$.
\end{itemize}
\end{definition}

\begin{remark}
  If $\bbQ$ is a strategically $({<}\lambda)$--complete forcing notion and
  $D$ is a (proper) normal filter on $\lambda$, then in $\bV^\bbQ$ the
  normal filter on $\lambda$ generated by $D\cap\bV$ is also a proper 
  filter. Abusing notation, we will denote this filter by $D$ (or
  $D^\bbQ$). The filter $D^\bp[\bbQ]$ can be larger, but it is still a
  proper filter, provided $\bp$ is a $\Dl$--parameter.
\end{remark}

\begin{lemma}
  \label{Dlprop}
Assume that $\bp=(\bar{P},S,D)$ is a $\Dl$--parameter on $\lambda$ and
$\bbQ$ is a strategically $({<}\lambda)$--complete forcing notion. Then 
$\forces_\bbQ\emptyset\notin D^\bp[\bbQ]$. Consequently, $\forces_\bbQ$``
$\bp[\bbQ]$ is a $\Dl$--parameter on $\lambda$''. 
\end{lemma}

\begin{proof}
  Assume that $p\in\bbQ$ and $\name{A}_\delta$ is a $\bbQ$--name for 
  an element of $D\cap\bV$ and $\name{f}_\delta$ is a $\bbQ$--name for  
  an element of ${}^\lambda\lambda$ (for $\delta<\lambda$). Using the
  strategic completeness of $\bbQ$ build a sequence $\langle
  p_\delta,A_\delta,f_\delta:\delta<\lambda\rangle$ such that for each
  $\delta<\lambda$:
\begin{enumerate}
\item[(i)]  $p_\delta\in\bbQ$, $p\leq p_0\leq p_\alpha\leq p_\delta$ for
  $\alpha<\delta$, 
\item[(ii)]  $A_\delta\in D\cap\bV$, $f_\delta\in {}^\lambda\lambda$ and 
\item[(iii)] $p_\delta\forces_\bbQ$`` $\name{A}_\delta=A_\delta$ and
  $\name{f}_\alpha\rest\delta=f_\alpha\rest \delta$ for all $\alpha\leq
  \delta$ ''. 
\end{enumerate}
Since $\bp$ is a $\Dl$--parameter, we know that $B=
\mathop{\triangle}\limits_{\delta< \lambda} A_\delta\cap
\mathop{\triangle}\limits_{\delta< \lambda}\set^\bp(f_\delta)\in D$. Let
$\delta\in B$. Then  
\[p_\delta\forces_\bbQ\mbox{`` } \delta\in \mathop{\triangle}\limits_{
  \alpha< \lambda}\name{A}_\alpha\mbox{ and }\name{f}_\alpha\rest \delta=
f_\alpha\rest\delta\in P_\delta\mbox{ for all }\alpha<\delta\mbox{ '',}\]
so  $p_\delta\forces_\bbQ\mbox{`` } \delta\in \mathop{\triangle}\limits_{
  \alpha< \lambda}\name{A}_\alpha\cap\mathop{\triangle}\limits_{\alpha<
  \lambda}\set^\bp(\name{f}_\alpha)$ ''.
\end{proof}

\begin{lemma}
\label{Dlmodels}
Assume that $\lambda^{<\lambda}=\lambda$, $\bp=(\bar{P},S,D)$ is a
$\Dl$--parameter on $\lambda$, $\bbQ$ is a strategically
$({<}\lambda)$--complete forcing notion and $N\prec (\cH(\chi),\in,
<^*_\chi)$ is such that $\bp\in N$, $|N|=\lambda$ and
${}^{<\lambda}N\subseteq N$. Let $\langle N_\delta:\delta<\lambda\rangle$ be
an increasing continuous sequence of elementary submodels of $N$ such that
$\bp\in N_0$, $\delta\subseteq N_\delta$, $P_\delta\subseteq N_{\delta+1}$,
$\langle N_\vare:\vare\leq\delta\rangle\in N_{\delta+1}$ and
$|N_\delta|<\lambda$ (for $\delta<\lambda$). Then
\[\forces_\bbQ\mbox{`` }\big(\forall A\subseteq N\big)\big(\big\{\delta<
\lambda: A\cap N_\delta\in N_{\delta+1}\big\}\in D[\bbQ]\big)\mbox{ ''.}\] 
\end{lemma}

\begin{proof}
We may find an increasing continuous sequence $\langle\alpha_\delta: \delta<
\lambda\rangle\subseteq\lambda$ and a bijection $f:N\longrightarrow\lambda$
such that $f[N_\delta]=\alpha_\delta$ and $f\rest N_\delta\in N_{\delta+1}$
(for $\delta<\lambda$). For $A\subseteq N$ let $\varphi_A:\lambda
\longrightarrow 2$ be such that $\varphi_A(\alpha)=1$ if and only if
$f^{-1}(\alpha)\in A$. Plainly, if $\delta=\alpha_\delta$ and
$\varphi_A\rest \delta\in P_\delta$, then $A\cap N_\delta\in N_{\delta+1}$.  
\end{proof}

Let $\bbC^\lambda_0$ be a forcing notion consisting of all pairs
$(\alpha,f)$ such that $\alpha<\lambda$ and $f\in\prod\limits_{\beta
  <\alpha} (\beta+1)$ ordered by the extension (so $(\alpha,f)\leq
(\alpha',f')$ if and only if $f\subseteq f'$). Thus it is a
$({<}\lambda)$--complete forcing notion which is an incarnation of the 
$\lambda$--Cohen forcing notion.

\begin{proposition}
  \label{badprop}
Assume $\lambda$ is strongly inaccessible. If $\bp=(\bar{P},S,D)$ is a
$\Dl$--parameter on $\lambda$ such that $(\forall\delta\in S)(|P_\delta|
\leq |\delta|)$, then $\forces_{\bbC^\lambda_0}$`` $D^{\bbC^\lambda_0}\neq
D[\bbC^\lambda_0]$ ''. 
\end{proposition}

\begin{proof}
Let $\name{f}$ be the canonical $\bbC^\lambda_0$--name for the generic
function in $\prod\limits_{\alpha<\lambda}(\alpha+1)$, so
$(\alpha,f)\forces_{\bbC^\lambda_0} f\subseteq\name{f}$. Plainly,
$\forces_{\bbC^\lambda_0}\set^\bp(\name{f})\in D[\bbC^\lambda_0]$ and we are
going to argue that $\forces_{\bbC^\lambda_0}\lambda\setminus\set^\bp(
\name{f})\in \big(D^{\bbC^\lambda_0}\big)^+$. To this end, suppose that
$p\in\bbC^\lambda_0$ and $\name{A}_\delta$ is a $\bbC^\lambda_0$--name for 
an element of $D\cap\bV$ (for $\delta<\lambda$). By induction on
$\xi<\lambda$ choose $\langle\alpha_\xi,B_\xi,\bar{p}^\xi:\xi<\lambda
\rangle$ so that   
\begin{enumerate}
\item[$(\alpha)$] $\langle\alpha_\xi:\xi<\lambda\rangle$ is an increasing
  continuous sequence of ordinals below $\lambda$, 
\item[$(\beta)$]  $B_\xi\in D$, $\bar{p}^\xi=\langle p^\xi_\sigma:\sigma\in
  \prod\limits_{\zeta<\xi}(\alpha_\zeta+1)\rangle\subseteq\bbC^\lambda_0$,
  $p^0_{\langle\rangle}=p=(\alpha_0,f^0_{\langle\rangle})$, 
\item[$(\gamma)$]  if $\sigma\in \prod\limits_{\zeta<\xi}(\alpha_\zeta+1)$,
  then $p^\xi_\sigma=(\alpha_\xi,f^\xi_\sigma)$ and $f^\xi_\sigma(
  \alpha_\zeta)=\sigma(\zeta)$ for $\zeta<\xi$, 
\item[$(\delta)$]  if $\xi<\xi'$, $\sigma'\in\prod\limits_{\zeta<\xi'}
  (\alpha_\zeta+1)$ and $\sigma=\sigma'\rest\xi$, then $p^\xi_\sigma\leq
  p^{\xi'}_{\sigma'}$, 
\item[$(\vare)$]  if $\xi<\lambda$ is limit, $\sigma\in
  \prod\limits_{\zeta<\xi}(\alpha_\zeta+1)$, then ($\alpha_\xi=\sup(
  \alpha_\zeta:\zeta<\xi)$ and) $f^\xi_\sigma=\bigcup\limits_{\zeta<\xi}
  f^\zeta_{\sigma\rest\zeta}$, and 
\item[$(\zeta)$] $p^{\xi+1}_\sigma\forces B_\xi\subseteq\name{A}_\xi$ for
  every $\sigma\in \prod\limits_{\zeta\leq\xi}(\alpha_\zeta+1)$.
\end{enumerate}
(Remember that $\lambda$ is inaccessible, so $\big|\prod\limits_{\zeta<\xi}
(\alpha_\zeta+1)\big|<\lambda$ for each $\xi$.) Next, consider the set $B=
\mathop{\triangle}\limits_{\xi<\lambda}B_\xi\in D$. Let $\delta\in B\cap S$
be a limit ordinal. Since $|P_\delta|\leq|\delta|<\prod\limits_{\xi<\delta}
|\alpha_\xi|$, we may pick $\sigma\in\prod\limits_{\xi<\delta}(\alpha_\xi
+1)$ such that $f^\delta_\sigma\notin P_\delta$. Then
$p^\delta_\sigma\forces \delta\in \mathop{\triangle}\limits_{\xi<\lambda}
\name{A}_\xi \setminus\set^\bp(\name{f})$. 
\end{proof}

\begin{definition}
  \label{Bplus}
Let $\bp=(\bar{P},S,D)$ be a $\Dl$--parameter on $\lambda$, $\bbQ$ be a
strategically $({<}\lambda)$--complete forcing notion.
\begin{enumerate}
\item For a condition $p\in\bbQ$ we define a game $\bgame(p,\bbQ)$ between  
two players, Generic and Antigeneric, as follows. A play of $\bgame(p,\bbQ)$
lasts $\lambda$ steps and during a play a sequence   
\[\Big\langle I_\alpha,\langle p^\alpha_t,q^\alpha_t:t\in I_\alpha\rangle:
\alpha<\lambda\Big\rangle\]   
is constructed. Suppose that the players have arrived to a stage
$\alpha<\lambda$ of the game. Now,
\begin{enumerate}
\item[$(\aleph)_\alpha$]  first Generic chooses a set $I_\alpha$ of
  cardinality $<\lambda$ and a system $\langle p^\alpha_t:t\in I_\alpha
  \rangle$ of conditions from $\bbQ$,\footnote{Note that no relation between
    $p^\alpha_t$ and $p^\beta_s$ for $\beta<\alpha$ is required to hold.}   
\item[$(\beth)_\alpha$]  then Antigeneric answers by picking a system
  $\langle q^\alpha_t:t\in I_\alpha\rangle$ of conditions from $\bbQ$ such
  that $(\forall t\in I_\alpha)(p^\alpha_t\leq q^\alpha_t)$.
\end{enumerate}
At the end, Generic wins the play $\big\langle I_\alpha, \langle p^\alpha_t,
q^\alpha_t:t\in I_\alpha\rangle:\alpha<\lambda\big\rangle$ of
$\bgame(p,\bbQ)$ if and only if 
\begin{enumerate}
\item[$(\circledast)^\bp_{\rm rbB}$] there is a condition $p^*\in\bbQ$
  stronger than $p$ and such that 
\[p^*\forces_{\bbQ}\mbox{`` }\big\{\alpha<\lambda:\big(\exists t\in
I_\alpha\big)\big(q^\alpha_t \in\Gamma_\bbQ\big)\big\}\in D[\bbQ]\mbox{
  ''.}\] 
\end{enumerate}
\item A forcing notion $\bbQ$ is {\em reasonably B--bounding over $\bp$\/}
  if for any $p\in\bbQ$, Generic has a winning strategy in the game
  $\bgame(p,\bbQ)$.      
\end{enumerate}
\end{definition}

\begin{remark}
  The notion introduced in \ref{Bplus} is almost the same as the one of
  \cite[Definition 3.1(2),(5)]{RoSh:860}. The difference is that in
  $(\circledast)^\bp_{\rm rbB}$ we use the filter $D[\bbQ]$ and not
  $D^\bbQ=D$, so potentially we have a weaker property here. We do not know,
  however, if there exists a forcing notion which is reasonably B--bounding
  over $\bp$ and not reasonably B--bounding over $D$. (See Problem
  \ref{prob0}.) 

  In a similar fashion we may also modify the property of being nicely
  double {\rm b}--bounding (see \cite[Definition 2.9(2),(4)]{RoSh:890}) and
  get the parallel iteration theorem.
\end{remark}

\begin{theorem}
  \label{like860}
Assume that
\begin{enumerate}
\item $\lambda$ is a strongly inaccessible cardinal and $\bp$ is a
  $\Dl$--parameter on $\lambda$,
\item $\bar{\bbQ}=\langle\bbP_\alpha,\name{\bbQ}_\alpha:\alpha<\gamma
  \rangle$ is a $\lambda$--support iteration,
\item for every $\alpha<\lambda$, $\forces_{\bbP_\alpha}$``
  $\name{\bbQ}_\alpha$ is reasonably B--bounding over $\bp[\bbP_\alpha]$ ''.  
\end{enumerate}
Then 
\begin{enumerate}
\item[(a)] $\bbP_\gamma=\lim(\bar{\bbQ})$ is $\lambda$--proper,
\item[(b)] if $\name{\tau}$ is a $\bbP_\gamma$--name for a function from
  $\lambda$ to $\bV$, $p\in\bbP_\gamma$, then there are $q\geq p$ and 
  $\langle A_\xi:\xi<\lambda\rangle$ such that $(\forall\xi<\lambda)(
  |A_\xi|<\lambda)$ and 
\[q\forces\mbox{`` }\{\xi<\lambda:\name{\tau}(\xi) \in A_\xi\}\in
D^{\bp}[\bbP_\gamma]\mbox{ ''.}\]  
\end{enumerate}
\end{theorem}

\begin{proof}
  The proof is essentially the same as that of \cite[Theorem 3.1]{RoSh:860}
  with a small modification at the end (in Claim 3.1 there); compare with
  the proof of Theorem \ref{lordsthm} here and specifically with \ref{cl3}. 
\end{proof}

\begin{definition}
  \label{exBp}
Let $\bar{E}=\langle E_\nu:\nu\in {}^{<\lambda}\lambda\rangle$ be a system
of $({<}\lambda)$--complete non-principal filters on $\lambda$ and let $E$
be a normal filter on $\lambda$. We define a forcing notion $\fEE$ as
follows. 

\noindent {\bf A condition $p$ in $\fEE$} is a complete $\lambda$--tree
$p\subseteq {}^{<\lambda}\lambda$ such that 
\begin{itemize}
\item for every $\nu\in p$, either $|\suc_p(\nu)|=1$ or $\suc_p(\nu)\in
  E_\nu$, and
\item for every $\eta\in\lim_\lambda(p)$ the set $\{\alpha<\lambda:
  \suc_p(\eta\rest\alpha)\in E_{\eta\rest\alpha}\}$ belongs to $E$.
\end{itemize}

\noindent {\bf The order $\leq=\leq_{\fEE}$} is the reverse inclusion:
$p\leq q$ if and only if ($p,q\in\fEE$ and ) $q\subseteq p$.
\end{definition}

\begin{proposition}
  \label{Qisgood}
Assume that $\bar{E},E$ are as in \ref{exBp}. Let $\bp=(\bar{P},S,D)$ be a
$\Dl$--parameter on $\lambda$ such that $\lambda\setminus S\in E$.
\begin{enumerate}
\item $\fEE$ is a $({<}\lambda)$--complete forcing notion of size
  $2^\lambda$. 
\item $\fEE$ is reasonably B--bounding over $\bp$.
\item If $\lambda$ is strongly inaccessible and $(\forall\delta\in S)(
  |P_\delta| \leq |\delta|)$, then $\forces_{\fEE} D^{\fEE}\neq D[\fEE]$. 
\end{enumerate}
\end{proposition}

\begin{proof}
(1)\quad Should be clear.
\smallskip

\noindent (2)\quad Let $p\in\fEE$. We are going to describe a strategy $\st$
for Generic in $\bgame(p,\fEE)$. In the course of the play, Generic
constructs aside a sequence $\langle T_\xi:\xi<\lambda\rangle$ so that if
$\big\langle I_\xi,\langle p^\xi_t,q^\xi_t:t\in I_\xi\rangle:\xi< \lambda 
\big\rangle$ is the sequence formed by the innings of the two players, then
the following conditions are satisfied.
\begin{enumerate}
\item[(a)] $T_\xi\in\fEE$ and if $\xi<\zeta<\lambda$ then $p=T_0\supseteq
  T_\xi\supseteq T_\zeta$ and $T_\zeta\cap {}^\xi\lambda= T_\xi\cap {}^\xi
  \lambda$. 
\item[(b)] If $\zeta<\lambda$ is limit, then
  $T_\zeta=\bigcap\limits_{\xi<\zeta} T_\xi$.  
\item[(c)] If $\xi\in S$ then
  \begin{itemize}
\item $I_\xi=P_\xi\cap T_\xi$ and $p^\xi_t=(T_\xi)_t$ for $t\in I_\xi$,  
\item $T_{\xi+1}=\bigcup\{q^\xi_t: t\in I_\xi\}\cup\bigcup\big\{(T_\xi)_\nu: 
  \nu\in {}^\xi\lambda\cap T_\xi\setminus I_\xi\big\}$.
  \end{itemize}
\item[(d)] If $\xi\notin S$, then $I_\xi=\emptyset$ and $T_{\xi+1}=T_\xi$.  
\end{enumerate}
Conditions (a)--(d) fully describe the strategy $\st$. Let us argue that it
is a winning strategy and to this end suppose that $\big\langle I_\xi,
\langle p^\xi_t,q^\xi_t:t\in I_\xi\rangle:\xi<\lambda\big\rangle$ is a play
of $\bgame(p,\fEE)$ in which Generic uses $\st$ and constructs aside the
sequence $\langle T_\xi:\xi<\lambda\rangle$ so that (a)--(d) are satisfied.
Put $p^*=\bigcap\limits_{\xi<\lambda} T_\xi\subseteq {}^{<\lambda}
\lambda$. It follows from (a)+(b) that $p^*$ is a complete $\lambda$--tree
and for each $\nu\in p^*$ either $|\suc_{p^*}(\nu)|=1$ or $\suc_{p^*}(\nu)
\in E_\nu$. Suppose now that $\eta\in\lim_\lambda(p^*)$ and for $\xi<
\lambda$ let $B_\xi\stackrel{\rm def}{=}\{\alpha<\lambda: \suc_{T_\xi}(\eta
\rest\alpha)\in E_{\eta\rest\alpha}\}$. Since $\eta\in\lim_\lambda(T_\xi)$
for each $\xi<\lambda$, we know that $B_\xi\in E$. Let 
\[B=\mathop{\triangle}\limits_{\xi<\lambda}B_\xi\cap \{\xi<\lambda:\xi\mbox{ 
  is limit and }\xi\notin S\}.\]
It follows from our assumptions that $B\in E$. For each $\alpha\in B$ we
know that $\suc_{T_\xi}(\eta\rest\alpha)\in E_{\eta\rest\alpha}$ for
$\xi<\alpha$ and $T_\alpha=\bigcap\limits_{\xi<\alpha}T_\xi$, so
$\suc_{T_\alpha}(\eta\rest\alpha)\in E_{\eta\rest\alpha}$. Moreover,
$T_\beta\cap {}^{\alpha+1}\lambda=T_\alpha\cap {}^{\alpha+1}\lambda$ for all
$\beta>\alpha$ (remember (a)+(d)) and consequently $\suc_{p^*}(\eta \rest
\alpha)=\suc_{T_\alpha}(\eta\rest\alpha)\in E_{\eta\rest\alpha}$. Thus we
have shown that $p^*\in\fEE$.  

Let $\name{W}$ be a $\fEE$--name given by $\forces_{\fEE}\name{W}=
\bigcup\{\mrot(p):p\in\Gamma_{\fEE}\}$. It should be clear that
$\forces_{\fEE}\name{W}\in {}^\lambda\lambda$ and thus $\forces_{\fEE}
\set^\bp(\name{W})\in D[\fEE]$. Plainly, if $\xi\in S$ and $t\in
{}^\xi\lambda \cap p^*$, then $(p^*)_t\geq q^\xi_t$ and hence 
\[p^*\forces_{\fEE}\mbox{`` if }\xi\in\set^\bp(\name{W}),
\mbox{ then }\name{W}\rest\xi\in I_\xi\mbox{ and }q^\xi_{\name{W}\rest\xi} 
\in \Gamma_{\fEE}\mbox{ '',}\]
so Generic won the play.
\smallskip

\noindent (3)\quad We are going to show that $\forces_{\fEE}\set^\bp
(\name{W})\notin D^{\fEE}$. To this end suppose that $p\in\fEE$ and
$\name{A}_\xi$ (for $\xi<\lambda$) are $\fEE$--names for elements of
$D$. Let $\st$ be the winning strategy of Generic in $\bgame(p,\fEE)$
described in part (2) above. Consider a play $\big\langle I_\xi,\langle
p^\xi_t,q^\xi_t:t\in I_\xi\rangle:\xi<\lambda\big\rangle$ of
$\bgame(p,\fEE)$ in which 
\begin{enumerate}
\item[$(*)_1$] Generic follows $\st$ and constructs aside a sequence
  $\langle T_\xi:\xi<\lambda\rangle$,
\item[$(*)_2$] Antigeneric plays so that at a stage $\xi\in S$ he picks a
  set $B_\xi\in D$ and conditions $q^\xi_t\geq p^\xi_t$ (for $t\in I_\xi$)
  such that  
\[\big(\forall t\in I_\xi\big)\big(q^\xi_t\forces B_\xi\subseteq
\bigcap\limits_{\zeta\leq\xi}\name{A}_\zeta\big).\]   
\end{enumerate}
Let $p^*=\bigcap\limits_{\xi<\lambda}T_\xi$ be the condition determined at
the end of part (2) and let $B=\mathop{\triangle}\limits_{\xi<\lambda}
B_\xi$. Choose an increasing continuous sequence $\langle\gamma_\xi:
\xi<\lambda\rangle\subseteq\lambda$ and a complete $\lambda$--tree
$T\subseteq p^*$ such that for every $\xi<\lambda$ we have
\begin{enumerate}
\item[$(*)_3$] if $\nu\in T\cap {}^{\gamma_\xi}\lambda$, then $|\{\rho\in
  T\cap {}^{\gamma_{\xi+1}}\lambda:\nu\vtl\rho\}|=|\gamma_\xi|$ and  
\item[$(*)_4$] if $\rho\in T\cap {}^{\gamma_{\xi+1}}\lambda$, then
  $\rho\rest\alpha\in P_\alpha$ for some $\alpha\in (\gamma_\xi,
  \gamma_{\xi+1}) \cap S$.
\end{enumerate}
(The choice can be done by induction on $\xi$; remember that $\bp$ is a 
$\Dl$--parameter and $\lambda$ is assumed to be inaccessible.) Pick a limit
ordinal $\xi\in B\cap S$ such that $\xi=\gamma_\xi$. Since $|T\cap
{}^\xi\lambda|>|\gamma_\xi|$, we may choose $\nu\in T\cap {}^\xi\lambda
\setminus P_\xi$. Put $q=(p^*)_\nu$. Then $q\geq p^*\geq p$ and
$q\forces_{\fEE} \xi\in \mathop{\triangle}\limits_{\xi<\lambda} \name{A}_\xi
\setminus \set^\bp(\name{W})$ (remember $(*)_2+(*)_4$).
\end{proof}

\section{Iterations with lords}
Theorem \ref{like860} can be used for $\lambda$-support iteration of forcing
notions $\fEE$ when on each coordinate we have the same filter $E$. But if
we want to use different filters on various coordinates we have serious
problems. However, if we move to the other extreme: having very orthogonal
filters we may use a different approach to argue that the limit of the
iteration is $\lambda$--proper. 

\begin{definition}
\label{purity}
\begin{enumerate}
\item {\em A forcing notion with $\lambda$--complete
  $(\kappa,\mu)$--purity\/} is a triple $(\bbQ,\leq,\leq_\pr)$ such that
  $\leq,\leq_\pr$ are transitive reflexive (binary) relations on $\bbQ$ such
  that 
\begin{enumerate}
\item[(a)]  $\leq_\pr\ \subseteq\ \leq$,
\item[(b)]  both $(\bbQ,\leq)$ and $(\bbQ,\leq_\pr)$ are strategically 
  $({<}\lambda)$--complete, 
\item[(c)]  for every $p\in\bbQ$ and a $(\bbQ,\leq)$--name $\name{\tau}$ for
  an ordinal below $\kappa$, there are a set $A$ of size less than $\mu$
  and a condition $q\in\bbQ$ such that $p\leq_\pr q$ and $q$ forces (in 
  $(\bbQ,\leq)$) that ``$\name{\tau}\in A$''. 
\end{enumerate}
\item If $(\bbQ,\leq,\leq_\pr)$ is a forcing notion with $\lambda$--complete
  $(\kappa,\mu)$--purity for every $\kappa$, then we say that it has
  {\em $\lambda$--complete $(*,\mu)$--purity}.  
\item If $(\bbQ,\leq,\leq_\pr)$ is a forcing notion with
  $\lambda$--complete $(\kappa,\mu)$--purity, then all our forcing terms
  (like ``forces'', ``name'' etc) refer to $(\bbQ,\leq)$. The relation
  $\leq_\pr$ has an auxiliary character only and if we want to refer to it
  we add ``purely'' (so ``$q$ is stronger than $p$'' means $p\leq q$, and
  ``$q$ is purely stronger than $p$'' means that $p\leq_\pr q$).  
\end{enumerate}
\end{definition}

\begin{definition}
\label{purgame}
Let $\bbQ=(\bbQ,\leq,\leq_\pr)$ be a forcing notion with $\lambda$--complete 
$(*,\lambda^+)$--purity, $\bp=(\bar{P},S,D)$ be a $\Dl$--parameter on
$\lambda$, $\cU$ be a normal filter on $\lambda$ and $\bar{\mu}= \langle
\mu_\alpha:\alpha<\lambda\rangle$ be a sequence of cardinals below
$\lambda$.  
\begin{enumerate}
\item For a condition $p\in\bbQ$ we define a game $\tezero(p,\bbQ)$ between 
two players, Generic and Antigeneric, as follows. A play of $\tezero(p,
\bbQ)$ lasts $\lambda$ steps and during a play a sequence  
\[\Big\langle\ell_\alpha,\langle p^\alpha_t,q^\alpha_t:t\in \mu_\alpha
\rangle: \alpha<\lambda\Big\rangle\]  
is constructed. So suppose that the players have arrived to a stage
$\alpha<\lambda$ of the game. Now,
\begin{enumerate}
\item[$(\aleph)^\pr_\alpha$]  first Antigeneric pics $\ell_\alpha\in
  \{0,1\}$.  
\item[$(\beth)^\pr_\alpha$]  After this, Generic chooses a system $\langle
  p^\alpha_t:t\in \mu_\alpha\rangle$ of paiwise incompatible conditions from
  $\bbQ$, and 
\item[$(\gimel)^\pr_\alpha$]
  Antigeneric answers with a system of conditions $q^\alpha_t\in\bbQ$ (for
  $t\in \mu_\alpha$) such that for each $t\in \mu_\alpha$:
\begin{itemize}
\item $p^\alpha_t\leq q^\alpha_t$, and   
\item if $\ell_\alpha=1$, then $p^\alpha_t\leq_\pr q^\alpha_t$.  
\end{itemize}
\end{enumerate}
At the end, Generic wins the play 
\[\Big\langle\ell_\alpha, \langle p^\alpha_t,q^\alpha_t:t\in \mu_\alpha
\rangle: \alpha<\lambda\Big\rangle\]  
if and only if either $\{\alpha<\lambda:\ell_\alpha=1\}\notin \cU$, or 
\begin{enumerate}
\item[$(\circledast)^\bp_\pr$] there is a condition $p^*\in\bbQ$ stronger
  than $p$ and such that
\[p^*\forces_{\bbQ}\mbox{`` }\big\{\alpha<\lambda:\big(\exists t\in
\mu_\alpha\big)\big(q^\alpha_t \in\Gamma_\bbQ\big)\}\in D[\bbQ]\mbox{
  ''.}\] 
\end{enumerate}
\item We say that the forcing notion $\bbQ$ (with $\lambda$--complete
  $(*,\lambda^+)$--purity) is {\em purely B$^*$--bounding over $\cU,\bp, 
  \bar{\mu}$\/} if for any $p\in\bbQ$, Generic has a winning strategy in the
  game $\tezero(p,\bbQ)$.     
\end{enumerate}
\end{definition}

\begin{remark}
Note that in the definition of the game $\tezero(p,\bbQ)$ the size of the
index set used at stage $\alpha$ is declared to be $\mu_\alpha$ (while in
the related game $\bgame(p,\bbQ)$ we required just
$|I_\alpha|<\lambda$). The reason for this is that otherwise in the proof of
the iteration theorem for the current case we could have
problems with deciding the size of the set $I_\alpha$; compare clause
$(*)_4$ of the proof of Theorem \ref{lordsthm}.
\end{remark}

\begin{observation}
  \label{expur}
Assume $\bar{E},E$ are as in \ref{exBp}. For $p,q\in \fEE$ let $p\leq_\pr q$
mean that $p\leq q$ and $\mrot(p)=\mrot(q)$. Then 
\begin{enumerate}
\item $(\fEE,\leq,\leq_\pr)$ is a forcing notion with $\lambda$--complete
  $(*,\lambda^+)$--purity,
\item if, additionally, each $E_\nu$ (for $\nu\in {}^{<\lambda}\lambda$) is
  an ultrafilter on $\lambda$, then $(\fEE,\leq,\leq_\pr)$ has
  $(\kappa,2)$--purity for every $\kappa<\lambda$.  
\end{enumerate}
\end{observation}

\begin{proposition}
  \label{exrc}
Assume that $\bar{E},E$ are as in \ref{exBp}, $\bp=(\bar{P},S,D)$ is a
$\Dl$--parameter on $\lambda$ and $\bar{\mu}=\langle\mu_\alpha:\alpha<
\lambda\rangle$ is a sequence of non-zero cardinals below $\lambda$ such
that $(\forall\alpha\in S)(|P_\alpha|\leq\mu_\alpha)$. Then $(\fEE,\leq,
\leq_\pr)$ is purely B$^*$--bounding over $E,\bp,\bar{\mu}$.
\end{proposition}

\begin{proof}
Let $p\in \fEE$ and let $\st$ be the strategy described in the proof of
\ref{Qisgood}(2) with a small modification that we start the construction
with $\xi_0=\lh(\mrot(p))+1$ (so $T_{\xi_0}=p$ and the first $\xi_0$ steps
of the play are not relevant). Then we also replace clauses (c)+(d) there by 
\begin{enumerate}
\item[(cd)] If $\xi\geq\xi_0$ then
  \begin{itemize}
\item $I_\xi\subseteq {}^\xi\lambda\cap T_\xi$ is of size $\mu_\xi$, and
  $p^\xi_t=(T_\xi)_t$ for $t\in I_\xi$, and 
\item if $\xi\in S$ then $P_\xi\cap T_\xi\subseteq I_\xi$, and 
\item $T_{\xi+1}=\bigcup\{q^\xi_t: t\in I_\xi\}\cup\bigcup\big\{(T_\xi)_\nu: 
  \nu\in {}^\xi\lambda\cap T_\xi\setminus I_\xi\big\}$.
  \end{itemize}
\end{enumerate}
(So, in particular, Antigeneric's choice of $\ell_\xi$ has no influence on
the answers by Generic.) We are going to show that $\st$ is a winning
strategy for Generic in $\teone(p,\fEE)$. To this end suppose that
$\big\langle \ell_\xi,\langle p^\xi_t,q^\xi_t:t\in I_\xi\rangle: \xi<\lambda
\big\rangle$ is a play of $\teone(p,\fEE)$ in which Generic follows $\st$
(we identify $I_\xi$ with $|I_\xi|=\mu_\xi$) and $\langle T_\xi:\xi< \lambda
\rangle$ is the sequence of side objects constructed in the course of the
play. Assume $A=\{\xi<\lambda: \ell_\xi=1\}\in E$ (otherwise Generic wins by
default). Like in \ref{Qisgood}(2), put $p^*=\bigcap\limits_{\xi<\lambda}
T_\xi$. To argue that $p^*\in\fEE$ we note that if
$\eta\in\lim_\lambda(p^*)$ and   
\[\delta\in \mathop{\triangle}\limits_{\xi<\lambda}\big\{\alpha< \lambda: 
\suc_{T_\xi}(\eta\rest\alpha)\in E_{\eta\rest\alpha}\big\}\cap A\cap \{\xi< 
\lambda: \xi>\xi_0\mbox{ is limit }\},\]
then 
\[\suc_{p^*}(\eta\rest\delta)=\left\{
\begin{array}{ll}
\suc_{T_\delta}(\eta\rest\delta)&\mbox{ if }\eta\rest\delta\notin I_\delta\\ 
\suc_{q^\delta_{\eta\rest\delta}}(\eta\rest\delta)&\mbox{ if }\eta\rest
\delta\in I_\delta\\ 
\end{array}\right.
\in E_{\eta\rest\delta}.\]
Exactly as in \ref{Qisgood}(2) we justify that $p^*$ witnesses
$(\circledast)^\bp_\pr$. 
\end{proof}

\begin{lemma}
  \label{lemlor}
Assume that
\begin{enumerate}
\item $\lambda$ is strongly inaccessible,
\item $\bar{\bbQ}=\langle\bbP_\alpha,\name{\bbQ}_\alpha:\alpha<\gamma
  \rangle$ is a $\lambda$--support iteration, $w\subseteq\gamma$,
  $|w|<\lambda$, $\alpha_0\in w$,
\item for every $\alpha<\gamma$, 
\[\forces_{\bbP_\alpha}\mbox{`` }\name{\bbQ}_\alpha=(\name{\bbQ}_\alpha,
\leq,\leq_\pr) \mbox{ is a forcing notion with $\lambda$--complete
  $(*,\lambda^+)$--purity '',}\] 
\item $\bbP_{\alpha_0}$ is $\lambda$--proper,
\item $\cT=(T,\rk)$ is a standard $(w,1)^\gamma$--tree, $|T|<\lambda$, 
\item $\bar{p}=\langle p_t:t\in T\rangle$ is a standard tree of conditions
  in $\bbP_\gamma$, and 
\item $\name{\tau}$ is a $\bbP_\gamma$--name for an ordinal.
\end{enumerate}
Then there are a set $A$ of size $\lambda$ and a standard tree of conditions
$\bar{q}=\langle q_t:t\in T\rangle\subseteq\bbP_\gamma$ such that 
\begin{enumerate}
\item[(a)] $\big(\forall t\in T\big)\big(\rk(t)=\gamma\ \Rightarrow\
  q_t\forces\name{\tau}\in A\big)$, and
\item[(b)] $\bar{p}\leq\bar{q}$ and if $t\in T$, $\rk(t)>\alpha_0$ then
  $q_{t\rest\alpha_0}\forces_{\bbP_{\alpha_0}} p_t(\alpha_0)\leq_\pr
  q_t(\alpha_0)$.  
\end{enumerate}
\end{lemma}

\begin{proof}
Let us start with the following observation.  

\begin{claim}
\label{cl2}
If $p\in\bbP_\gamma$ then there are a set $A_0$ of size $\lambda$ and a
condition $q\geq p$ such that $q\forces_{\bbP_\gamma}\name{\tau}\in A_0$ and
$q\rest\alpha_0\forces_{\bbP_{\alpha_0}}p(\alpha_0)\leq_\pr q(\alpha_0)$.  
\end{claim}

\begin{proof}[Proof of the Claim]
Let us look at $\bbP_\gamma$ as the result of 3 stage composition
$\bbP_{\alpha_0}*\name{\bbQ}_{\alpha_0}*\name{\bbP}_{(\alpha_0+1),\gamma}$,
where $\name{\bbP}_{(\alpha_0+1),\gamma}$ is a $\bbP_{\alpha_0+1}$--name for
the following forcing notion. The set of conditions in
$\name{\bbP}_{(\alpha_0+1),\gamma}$ is $\{r\rest (\alpha_0,\gamma):r\in
\bbP_\gamma\}$ (so it belongs to $\bV$); the order of
$\name{\bbP}_{(\alpha_0+1),\gamma}$ is such that if $G_{\alpha_0+1}\subseteq
\bbP_{\alpha_0+1}$ is generic over $\bV$, then 
\[\begin{array}{ll}
\bV[G_{\alpha_0+1}]\models&\mbox{`` }r
\leq_{\name{\bbP}_{(\alpha_0+1),\gamma}[G_{\alpha_0}]} s\mbox{ if and only
  if}\\ 
&\ \mbox{ there is } q\in G_{\alpha_0+1}\mbox{ such that }q \conc r
\leq_{\bbP_\gamma} q\conc s\mbox{ ''.} 
\end{array}\]
Now, pick a $\bbP_{\alpha_0+1}$--name $(\name{r},\name{\alpha})$ such that 
\[p\rest (\alpha_0+1)\forces_{\bbP_{\alpha_0+1}}\mbox{`` }p\rest
(\alpha_0,\gamma)\leq \name{r}\mbox{ and }\name{r}\forces\name{\alpha}
=\name{\tau}\mbox{ ''} \]
and then choose a $\bbP_{\alpha_0}$--name $\name{A}^*$ for a subset of
$\name{\bbP}_{(\alpha_0+1),\gamma}\times {\rm ON}$ and a
$\bbP_{\alpha_0}$--name $q(\alpha_0)$ for a condition in
$\name{\bbQ}_{\alpha_0}$ such that 
\[\begin{array}{ll}
p\rest\alpha_0\forces_{\bbP_{\alpha_0}}&\mbox{`` }p(\alpha_0)\leq_\pr
q(\alpha_0)\mbox{ and }|\name{A}^*|=\lambda\mbox{ and}\\
&\ \ q(\alpha_0)\forces_{\name{\bbQ}_{\alpha_0}}\big(\exists (s,\beta)\in
\name{A}^*\big)\big(\name{r}=s\ \&\ \name{\alpha}=\beta\big)\mbox{ ''.}  
\end{array}\]
Since $\bbP_{\alpha_0}$ is $\lambda$--proper, we may choose a set
$A^+\subseteq \name{\bbP}_{(\alpha_0+1),\gamma}\times{\rm ON}$ of size
$\lambda$ and a condition $q\rest \alpha_0\geq p\rest\alpha_0$ such that
$q\rest\alpha_0\forces\name{A}^*\subseteq A^+$. Then 
\[q\rest(\alpha_0+1)\forces_{\bbP_{\alpha_0+1}}\big(\exists (s,\beta)\in A^+
\big)\big(\name{r}=s\ \&\ \name{\alpha}=\beta\big).\]
Put $A=\{\beta:(\exists s)((s,\beta)\in A^+)\}$. Now we may easily define
$q\rest (\alpha_0,\gamma)$ so that $\dom\big(q\rest (\alpha_0,\gamma)\big)
=\bigcup\{\dom(s):(\exists\beta)((s,\beta)\in A^+)\}$ and\\ 
$q\rest (\alpha_0+1)\forces_{\bbP_{\alpha_0+1}}\mbox{`` }\name{r}
\leq_{\name{\bbP}_{(\alpha_0+1),\gamma}} q\rest (\alpha_0,\gamma)\mbox{ and
} \name{\alpha}\in A$ ''.
\end{proof}

Fix an enumeration $\langle t_\zeta:\zeta\leq\zeta^*\rangle$ of $\{t\in T:
\rk(t)=\zeta\}$ (so $\zeta^*<\lambda$). For each $\alpha\in\gamma\setminus
\{\alpha_0\}$ fix a $\bbP_\alpha$--name $\name{\st}^0_\alpha$ for a winning
strategy of Complete in the game $\Game^\lambda_0\big((\name{\bbQ}_\alpha,
\leq),\name{\emptyset}_{\name{\bbQ}_\alpha}\big)$ such that as long as
Incomplete plays $\name{\emptyset}_{\name{\bbQ}_\alpha}$, Complete answers
with $\name{\emptyset}_{\name{\bbQ}_\alpha}$ as well. Let
$\name{\st}^\zeta_\pr$ be the $<^*_\chi$--first $\bbP_{\alpha_0}$--name for
a winning strategy of Complete in $\Game^\lambda_0\big((
\name{\bbQ}_{\alpha_0},\leq_\pr), p_{t_\zeta}(\alpha_0)\big)$ (for
$\zeta\leq\zeta^*$). Note that if $\zeta,\zeta'\leq\zeta^*$ and
$t_\zeta\rest (\alpha_0+1)= t_{\zeta'}\rest (\alpha_0+1)$, then
$\name{\st}^\zeta_\pr=\name{\st}^{\zeta'}_\pr$. 

By induction on $\zeta\leq\zeta^*$ we choose a sequence $\langle
\bar{p}^\zeta,\bar{q}^\zeta,A^\zeta:\zeta\leq\zeta^*\rangle$ so that the
following demands are satisfied.
\begin{enumerate}
\item[(i)] $\bar{p}^\zeta=\langle p^\zeta_t:t\in T\rangle, \bar{q}^\zeta=
  \langle q^\zeta_t:t\in T\rangle$ are standard trees of conditions,
  $A^\zeta$ is a set of ordinals of size $\lambda$.
\item[(ii)] If $\vare<\zeta\leq\zeta^*$, then $\bar{p}\leq \bar{p}^\vare
  \leq \bar{q}^\vare\leq\bar{p}^\zeta$ and $A^\vare\subseteq A^\zeta$.
\item[(iii)] $p^\zeta_{t_\zeta}\forces_{\bbP_\gamma}\name{\tau}\in A^\zeta$.
\item[(iv)] If $\alpha\in\gamma\setminus\{\alpha_0\}$, $\zeta,\xi\leq
  \zeta^*$, then 
\[\begin{array}{ll}
  q^\xi_{t_\zeta\rest\alpha}\forces_{\bbP_\alpha}&\mbox{`` }\langle
  p^\vare_{t_\zeta\rest\alpha}(\alpha), q^\vare_{t_\zeta\rest\alpha} 
  (\alpha): \vare\leq\xi\rangle\mbox{ is a result of a play of
  }\Game^\lambda_0 \big((\name{\bbQ}_\alpha,\leq),
  \name{\emptyset}_{\name{\bbQ}_\alpha}\big)\\ 
&\ \mbox{ in which Complete uses }\name{\st}^0_\alpha\mbox{ ''.} 
\end{array}\]
\item[(v)]  If $\zeta,\xi\leq\zeta^*$, then 
\[\begin{array}{ll}
  q^\xi_{t_\zeta\rest\alpha_0}\forces_{\bbP_{\alpha_0}}&\mbox{`` }\langle 
  p^\vare_{t_\zeta\rest\alpha_0}(\alpha_0), q^\vare_{t_\zeta\rest\alpha_0}  
  (\alpha_0): \vare\leq\xi\rangle\mbox{ is a result of a play}\\
&\ \mbox{ of }\Game^\lambda_0 \big((\name{\bbQ}_{\alpha_0},\leq_\pr),
  p_{t_\zeta}(\alpha_0)\big)\mbox{ in which Complete uses
  }\name{\st}^\zeta_\pr\mbox{ ''.}  
\end{array}\]
\end{enumerate}
Suppose that we have determined $\bar{p}^\vare,\bar{q}^\vare,A^\vare$ for
$\vare<\zeta\leq\zeta^*$. First we choose $\bar{p}'=\langle p'_t:t\in
T\rangle\subseteq\bbP_\gamma$. If $\xi=0$ then we set
$\bar{p}'=\bar{p}$. Otherwise we choose $\bar{p}'$ so that for $t\in T$ we
have: 
\begin{enumerate}
\item[(vi)] $\dom(p'_t)=\bigcup\limits_{\vare<\zeta}\dom(q^\vare_t)$, and  
\item[(vii)] if $\alpha\in\dom(p'_t)\setminus\{\alpha_0\}$, then
$p'(\alpha)$ is the $<^*_\chi$--first $\bbP_\alpha$--name for a condition
in $\name{\bbQ}_\alpha$ such that $p'_t\rest\alpha\forces_{\bbP_\alpha}
(\forall\vare<\zeta)(q^\vare_t(\alpha)\leq p'_t(\alpha))$, and
\item[(viii)] $p'(\alpha_0)$ is the $<^*_\chi$--first
  $\bbP_{\alpha_0}$--name for a condition in $\name{\bbQ}_{\alpha_0}$ such
  that 
\[p'_t\rest\alpha_0 \forces_{\bbP_{\alpha_0}}(\forall\vare<\zeta)(
  q^\vare_t(\alpha_0) \leq_\pr p'_t(\alpha_0)).\]
\end{enumerate}
The choice is possible by (iv)+(v), and since we pick ``the
$<^*_\chi$--first names'' we easily see that $\bar{p}'$ is a standard tree
of conditions. Now we use \ref{cl2} to find a set $A^\zeta$ of size
$\lambda$ and a condition $p^\zeta_{t_\zeta}\in\bbP_\gamma$ such that  
\[\bigcup\limits_{\vare<\zeta}A^\vare\subseteq A^\zeta,\quad
p'_{t_\zeta}\leq p^\zeta_{t_\zeta},\quad p^\zeta_{t_\zeta}\rest\alpha_0
\forces_{\bbP_{\alpha_0}} p'_{t_\zeta}(\alpha_0)\leq_\pr p^\zeta_{t_\zeta}
(\alpha_0)\ \mbox{ and }\ p^\zeta_{t_\zeta}\forces_{\bbP_\gamma}
\name{\tau}\in A^\zeta.\]
Next, for each $t\in T$ we let $p^\zeta_t\in\bbP_{\rk(t)}$ be such that 
\begin{enumerate}
\item[(ix)] if $s=t\cap t_\zeta$, then\\
$p^\zeta_t\rest\rk(s)=p^\zeta_{t_\zeta}\rest\rk(s)$ and $p^\zeta_t\rest
[\rk(s),\rk(t))=p'_t\rest [\rk(s),\rk(t))$.  
\end{enumerate}
Clearly, $\bar{p}^\zeta=\langle p^\zeta_t:t\in T\rangle$ is a standard tree
of conditions satisfying the relevant parts of the demands in (ii)--(v). Now
we choose a tree of conditions $\bar{q}^\zeta=\langle q^\zeta_t:t\in
T\rangle$ so that the requirements of (iv)+(v) hold (for this we proceed
like in (vi)--(viii) above). 

After the construction is carried out we note that $\bar{q}^{\zeta^*}$ and
$A^{\zeta^*}$ are as required in the assertion of the lemma.  
\end{proof}

\begin{theorem}
\label{lordsthm}
Assume that 
\begin{enumerate}
\item $\lambda$ is strongly inaccessible, $\bar{\mu}=\langle\mu_\alpha:
  \alpha<\lambda\rangle$ is a sequence of cardinals below $\lambda$,
  $\bp=(\bar{P},S,D)$ is a $\Dl$--parameter on $\lambda$, and  
\item $\bar{\bbQ}=\langle\bbP_\alpha,\name{\bbQ}_\alpha:\alpha<\gamma
  \rangle$ is a $\lambda$--support iteration,
\item $\name{\cU}_\alpha$ is a $\bbP_\alpha$--name for a normal filter on
 $\lambda$ (for $\alpha<\gamma$), 
\item $A_{\alpha,\beta}\subseteq \lambda$ is such that
 $\forces_{\bbP_\alpha} A_{\alpha,\beta}\in\name{\cU}_\alpha$ and
 $\forces_{\bbP_\beta} \lambda\setminus A_{\alpha,\beta}\in
 \name{\cU}_\beta$ (for $\alpha<\beta<\gamma$),  and 
\item for every $\alpha<\gamma$,    
\[\forces_{\bbP_\alpha}\mbox{`` $\name{\bbQ}_\alpha$ is purely
  B$^*$--bounding over $\name{\cU}_\alpha,\bp[\bbP_\alpha],\bar{\mu}$ ''.}\] 
\end{enumerate}
Then $\bbP_{\gamma}=\lim(\bar{\bbQ})$ is $\lambda$--proper. 
\end{theorem}

\begin{proof}
  The arguments follow closely the lines of the arguments for
  \cite[Thm. 3.1, 3.2]{RoSh:860} and \cite[Thm. 2.12]{RoSh:890}. The proof
  is by induction on $\gamma$, so assume that we know also that each
  $\bbP_\alpha$ is $\lambda$--proper for $\alpha<\lambda$.

Let $N\prec (\cH(\chi),\in,<^*_\chi)$ be such that ${}^{<\lambda}N\subseteq
N$, $|N|=\lambda$ and $\bar{\bbQ},\langle A_{\alpha,\beta}:\alpha<\beta<
\gamma\rangle,\bp,\ldots\in N$. Let $p\in N\cap\bbP_\gamma$ and $\langle
\name{\tau}_\delta:\delta<\lambda\rangle$ list all $\bbP_\gamma$--names for 
ordinals from $N$. For each $\xi\in N\cap\gamma$ fix a $\bbP_\xi$--name
$\name{\st}^0_\xi\in N$ for a winning strategy of Complete in
$\Game_0^\lambda(\name{\bbQ}_\xi,\name{\emptyset}_{\name{\bbQ}_\xi})$ such
that it instructs Complete to play $\name{\emptyset}_{\name{\bbQ}_\xi}$ as
long as her opponent plays $\name{\emptyset}_{\name{\bbQ}_\xi}$. 

By induction on $\delta<\lambda$ we will choose  
\begin{enumerate}
\item[$(\otimes)_\delta^a$] \qquad $\cT_\delta,\bar{p}^\delta_*,
\bar{q}^\delta_*,r^-_\delta, r_\delta,w_\delta,Z_\delta,\alpha_\delta$, and  
\item[$(\otimes)_\delta^b$] \qquad $\ell_{\delta,\xi},
  \name{\bar{p}}_{\delta,\xi}, \name{\bar{q}}_{\delta,\xi}$ and
  $\name{\st}_\xi$ for $\xi\in N\cap \gamma$,    
\end{enumerate}
so that the following demands are satisfied. 
\begin{enumerate}
\item[$(*)_0$] All objects listed in $(\otimes)_\delta^a+(\otimes)_\delta^b$
  belong to $N$. After stage $\delta<\lambda$ of the construction, the
  objects in $(\otimes)_\delta^a$ are known as well as those in
  $(\otimes)_\delta^b$ for $\xi\in w_\delta$. 
\item[$(*)_1$] $r^-_\delta,r_\delta\in \bbP_\gamma$, $r^-_0(0)=r_0(0)=p(0)$, 
$w_\delta\subseteq\gamma$, $|w_\delta|=|\delta|+1$, $\bigcup\limits_{\alpha<  
\lambda}\dom(r_\alpha)=\bigcup\limits_{\alpha<\lambda} w_\alpha= N\cap
\gamma$, $w_0=\{0\}$, $w_\delta\subseteq w_{\delta+1}$ and if $\delta$ is
limit then $w_\delta=\bigcup\limits_{\alpha<\delta} w_\alpha$.
\item[$(*)_2$] For each $\alpha<\delta<\lambda$ we have $(\forall\xi\in
w_{\alpha+1})(r_\alpha(\xi)=r_\delta(\xi))$ and $p\leq r_\alpha^-\leq
r_\alpha\leq r^-_\delta\leq r_\delta$. 
\item[$(*)_3$] If $\xi\in(\gamma\setminus w_\delta)\cap N$, then 
\[\begin{array}{ll}
r_\delta\rest\xi\forces&\mbox{`` the sequence }\langle r^-_\alpha(\xi),
r_\alpha(\xi):\alpha\leq\delta\rangle\mbox{ is a legal partial play of }\\
&\quad\Game_0^\lambda\big(\name{\bbQ}_\xi,
\name{\emptyset}_{\name{\bbQ}_\xi}\big)\mbox{ in which Complete follows
}\name{\st}^0_\xi\mbox{ ''}  
\end{array}\]
and if $\xi\in w_{\delta+1}\setminus w_\delta$, then $\name{\st}_\xi\in N$
is a $\bbP_\xi$--name  for a winning strategy of Generic in
$\tegamexi(r_\delta(\xi),\name{\bbQ}_\xi)$. (And $\st_0\in N$ is a winning
strategy of Generic in ${\Game^{{\rm pr}}_{\name{\cU}_0,\bp,\bar{\mu}}}
(p(0),\bbQ_0)$.)    
\item[$(*)_4$] $\cT_\delta=(T_\delta,\rk_\delta)$ is a standard $(w_\delta,
1)^\gamma$--tree, $T_\delta=\bigcup\limits_{\alpha\leq\gamma}
\prod\limits_{\xi\in w_\delta\cap\alpha}\mu_\delta$.
\item[$(*)_5$] $\bar{p}^\delta_*=\langle p^\delta_{*,t}:t\in T_\delta
\rangle$ and $\bar{q}^\delta_*=\langle q^\delta_{*,t}:t\in T_\delta\rangle$
are standard trees of conditions, $\bar{p}^\delta_*\leq\bar{q}^\delta_*$.  
\item[$(*)_6$] For $t\in T_\delta$ we have that $\dom(p^\delta_{*,t})=
  \big(\dom(p)\cup\bigcup\limits_{\alpha<\delta}\dom(r_\alpha)\cup
  w_\delta\big) \cap \rk_\delta(t)$ and for each $\xi\in
  \dom(p^\delta_{*,t})\setminus w_\delta$: 
\[\begin{array}{ll}
p^\delta_{*,t}\rest\xi\forces_{\bbP_\xi}&\mbox{`` if the set } \{
r_\alpha(\xi):\alpha<\delta\}\cup\{p(\xi)\}\mbox{ has an upper bound in
}\name{\bbQ}_\xi,\\  
&\mbox{\quad then $p^\delta_{*,t}(\xi)$ is such an upper bound ''.}
  \end{array}\]
\item[$(*)_7$] For $\xi\in N\cap\gamma$, $\ell_{\delta,\xi}\in \{0,1\}$ and
  $\name{\bar{p}}_{\delta,\xi},\name{\bar{q}}_{\delta,\xi}$ are
  $\bbP_\xi$--names for sequences of conditions in $\name{\bbQ}_\xi$ of
  length $\mu_\delta$.
\item[$(*)_8$] If either $\xi=0=\beta$ or $\xi\in w_{\beta+1}\setminus
  w_\beta$, $\beta<\lambda$, then 
\[\begin{array}{r}
\forces_{\bbP_\xi}\mbox{`` }\langle\ell_{\alpha,\xi},
\name{\bar{p}}_{\alpha,\xi},\name{\bar{q}}_{\alpha,\xi}:\alpha<\lambda
\rangle\mbox{ is a play of }\tegamexi(r_\beta(\xi),\name{\bbQ}_\xi)\\  
\mbox{ in which Generic uses $\name{\st}_\xi$ ''.}
  \end{array}\]
\item[$(*)_9$] $\alpha_\delta\in w_\delta$ ($\alpha_\delta$ will be called
  {\em the lord of stage $\delta$}) and 
\[\big(\forall\beta\in w_\delta\setminus\{\alpha_\delta\}\big)\big(\delta
\notin \bigcap\big\{\lambda\setminus A_{\xi,\beta}:\xi\in w_\delta\cap
\beta\big\}\cap \bigcap\big\{A_{\beta,\xi}:\xi\in w_\delta\setminus
(\beta+1) \big\}\big).\] 
\item[$(*)_{10}$] $\ell_{\delta,\xi}=0$ for $\xi\in N\cap\gamma\setminus
  \{\alpha_\delta\}$ and $\ell_{\delta,\alpha_\delta}=1$.
\item[$(*)_{11}$] If $t\in T_\delta$, $\rk_\delta(t)=\xi<\gamma$, then for
  each $\vare<\mu_\delta$
\[q^\delta_{*,t}\forces_{\bbP_\xi}\mbox{`` } \name{\bar{p}}_{\delta,\xi}
(\vare)=p^\delta_{*,t\conc\langle\vare\rangle}(\xi)\mbox{ and }
\name{\bar{q}}_{\delta,\xi}(\vare)=q^\delta_{*,t\conc\langle\vare\rangle} 
(\xi)\mbox{ ''.}\] 
\item[$(*)_{12}$] If $t_0,t_1\in T_\delta$, $\rk_\delta(t_0)=
\rk_\delta(t_1)$ and $\xi\in w_\delta\cap\rk_\delta(t_0)$,
$t_0\rest\xi=t_1\rest\xi$ but $\big(t_0\big)_\xi\neq \big(t_1\big)_\xi$, 
then 
\[p^\delta_{*,t_0\rest\xi}\forces_{\bbP_\xi}\mbox{`` the conditions  
$p^\delta_{*,t_0}(\xi),p^\delta_{*,t_1}(\xi)$ are incompatible ''.}\] 
\item[$(*)_{13}$] $Z_\delta$ is a set of ordinals, $|Z_\delta|=\lambda$ and
  for each $t\in T_\delta$ with $\rk_\delta(t)=\gamma$ we have
  $q^\delta_{*,t}\forces_{\bbP_\gamma}\big(\forall\alpha\leq \delta\big)
  \big(\name{\tau}_\alpha\in Z_\delta\big)$.
\item[$(*)_{14}$] $\dom(r_\delta^-)=\dom(r_\delta)= \bigcup\limits_{t\in
    T_\delta} \dom(q^\delta_{*,t})\cup\dom(p)$ and if $t\in T_\delta$,
$\xi\in\dom(r_\delta)\cap \rk_\delta(t)\setminus w_\delta$, and
$q^\delta_{*,t}\rest\xi\leq q\in\bbP_\xi$, $r_\delta\rest\xi\leq q$, then   
\[\begin{array}{ll}
q\forces_{\bbP_\xi}&\mbox{`` if the set }\{r_\alpha(\xi):\alpha<\delta\}
\cup\{q^\delta_{*,t}(\xi), p(\xi)\}\mbox{ has an upper bound in }  
\name{\bbQ}_\xi,\\ 
&\mbox{\quad then $r_\delta^-(\xi)$ is such an upper bound ''.}
  \end{array}\]
\end{enumerate}

First we fix an increasing continuous sequence $\langle w_\alpha:\alpha<
\lambda\rangle$ of subsets of $N\cap\gamma$ such that the relevant demands
in $(*)_1$ are satisfied. Now, suppose that we have arrived to a stage
$\delta<\lambda$ of the construction and all objects listed in
$(\otimes)_\alpha^a$ and relevant cases of $(\otimes)_\alpha^b$ (see
$(*)_0$) have been determined for $\alpha<\delta$.

To ensure $(*)_0$, all choices below are made in $N$ (e.g., each time we
choose an object with some properties, we pick the $<^*_\chi$--first such
object). 

If $\delta$ is a successor ordinal and $\xi\in w_\delta\setminus
w_{\delta-1}$, then we let $\name{\st}_\xi\in N$ be a $\bbP_\xi$--name for a
winning strategy of Generic in $\tegamexi(r_{\delta-1}(\xi), 
\name{\bbQ}_\xi)$. We also put $\ell_{\alpha,\xi}=0$ for all $\alpha<\delta$
and we pick $\name{\bar{p}}_{\alpha,\xi}, \name{\bar{q}}_{\alpha,\xi}$ (for
$\alpha<\delta$) so that the suitable parts of $(*)_7+(*)_8$ at $\xi$ are
satisfied. 

Clause $(*)_4$ fully describes $\cT_\delta$. Now we choose the lord of stage
$\delta$. If for some $\beta\in w_\delta$ we have 
\[\delta\in \bigcap\big\{\lambda\setminus A_{\xi,\beta}:\xi\in w_\delta\cap 
\beta\big\}\cap \bigcap\big\{A_{\beta,\xi}:\xi\in w_\delta\setminus
(\beta+1) \big\},\] 
then $\alpha_\delta$ is equal to this $\beta$ (note that there is at most
one $\beta\in w_\delta$ with the required property). Otherwise we let
$\alpha_\delta=0$. Then we put $\ell_{\delta,\alpha_\delta}=1$ and
$\ell_{\delta,\xi}=0$ for all $\xi\in w_\delta\setminus\{\alpha_\delta\}$.

Next, for each $\xi\in w_\delta$ we choose a $\bbP_\xi$--name
$\name{\bar{p}}_{\delta,\xi}$ such that 
\[\begin{array}{ll}
\forces_{\bbP_\xi}&\mbox{`` }\name{\bar{p}}_{\delta,\xi}=\langle
\name{\bar{p}}_{\delta,\xi}(\vare):\vare<\mu_\delta\rangle\mbox{ is given to
Generic by }\name{\st}_\xi\\ 
&\quad\mbox{as an answer to }\langle\ell_{\alpha,\xi},
\name{\bar{p}}_{\alpha,\xi},\name{\bar{q}}_{\alpha,\xi}:\alpha<\delta\rangle
\conc\langle\ell_{\delta,\xi}\rangle\mbox{ ''}.
  \end{array}\]
(Note that $\forces_{\bbP_\xi}$`` conditions $\name{\bar{p}}_{\delta,\xi}(
\vare_0), \name{\bar{p}}_{\delta,\xi}(\vare_1)$ are incompatible'' whenever
$\vare_0<\vare_1<\mu_\delta$ and $\xi\in w_\delta$.) After this we may
choose a tree of conditions $\bar{p}^\delta_*=\langle p^\delta_{*,t}: t\in
T_\delta\rangle$ such that for each $t\in T_\delta$: 
\begin{itemize}
\item $\dom(p^\delta_{*,t})=\big(\dom(p)\cup \bigcup\limits_{\alpha<\delta} 
  \dom(r_\alpha) \cup w_\delta\big) \cap \rk_\delta(t)$ and 
\item for $\xi\in\dom(p^\delta_{*,t})\setminus w_\delta$,
  $p^\delta_{*,t}(\xi)$ is the $<^*_\chi$--first $\bbP_\xi$--name for a
  condition in $\name{\bbQ}_\xi$ such that  
\[\begin{array}{ll}
p^\delta_{*,t}\rest\xi\forces_{\bbP_\xi}&\mbox{`` if the set } \{
r_\alpha(\xi):\alpha<\delta\}\cup\{p(\xi)\}\mbox{ has an upper bound in
}\name{\bbQ}_\xi,\\  
&\mbox{\quad then $p^\delta_{*,t}(\xi)$ is such an upper bound '',}
\end{array}\]
\item $p^\delta_{*,t}(\xi)=\name{\bar{p}}_{\alpha,\xi}\big((t)_\xi\big)$ for
  $\xi\in \dom(p^\delta_{*,t})\cap w_\delta$. 
\end{itemize}
Using Lemma \ref{lemlor} we may pick a tree of conditions $\bar{q}^\delta_*
=\langle q^\delta_{*,t}:t\in T_\delta\rangle$ and a set $Z_\delta$ of
ordinals such that 
\begin{itemize}
\item $\bar{p}^\delta_*\leq\bar{q}^\delta_*$, $|Z_\delta|=\lambda$,
\item if $t\in T_\delta$, $\rk_\delta(t)=\gamma$ then $q^\delta_{*,t}
  \forces_{\bbP_\gamma} (\forall\alpha\leq\delta)(\name{\tau}_\alpha\in
  Z_\delta)$,
\item if $t\in T_\delta$, $\rk_\delta(t)>\alpha_\delta$ then
  $q^\delta_{*,t\rest \alpha_\delta}\forces_{\bbP_{\alpha_\delta}}
  p^\delta_{*,t}( \alpha_\delta)\leq_{\pr} q^\delta_{*,t}(\alpha_\delta)$. 
\end{itemize}
Note that if $\xi\in w_\delta$ and $\vare_0<\vare_1<\mu_\delta$, $t\in
T_\delta$, $\rk_\delta(t)=\xi$, then 
\[q^\delta_{*,t}\forces_{\bbP_\xi}\mbox{`` the conditions }q^\delta_{*,
  t\conc\langle\vare_0\rangle}(\xi), q^\delta_{*,t\conc\langle
  \vare\rangle} (\xi)\mbox{ are incompatible ''}.\]
Hence we have no problems with finding $\bbP_\xi$--names
$\name{\bar{q}}_{\delta,\xi}$ (for $\xi\in w_\delta$) such that 
\begin{itemize}
\item $\forces_{\bbP_\xi}$`` $\name{\bar{q}}_{\delta,\xi}=\langle
  \name{\bar{q}}_{\delta,\xi}(\vare):\vare<\mu_\delta\rangle$ is a sequence
  of conditions in $\name{\bbQ}_\xi$ '',
\item $\forces_{\bbP_\xi}$`` $(\forall\vare<\mu_\delta)
  (\name{\bar{p}}_{\delta,\xi}(\vare)\leq\name{\bar{q}}_{\delta,\xi}(\vare))$ 
  '' and $\forces_{\bbP_{\alpha_\delta}}$`` $(\forall\vare<\mu_\delta)
  (\name{\bar{p}}_{\delta,\alpha_\delta}(\vare)\leq_\pr
  \name{\bar{q}}_{\delta,\alpha_\delta}(\vare))$ '',
\item if $t\in T_\delta$, $\rk_\delta(t)>\xi$, then $q^\delta_{*,t\rest\xi}
  \forces_{\bbP_\xi}q^\delta_{*,t}(\xi)=\name{\bar{q}}_{\delta,\xi}\big(
  (t)_\xi\big)$.  
\end{itemize}
Now we define $r^-_\delta,r_\delta\in\bbP_\gamma$ so that
\[\dom(r^-_\delta)=\dom(r_\delta)=\bigcup\limits_{t\in T_\delta}\dom(
q^\delta_{*,t})\cup \dom(p)\]
 and 
\begin{itemize}
\item if $\xi\in w_{\alpha+1}$, $\alpha<\delta$, then $r^-_\delta(\xi) =
  r_\delta(\xi)=r_\alpha(\xi)$, 
\item if $\xi\in\dom(r^-_\delta)\setminus w_\delta$, then $r^-_\delta(\xi)$
  is the $<^*_\chi$--first $\bbP_\xi$--name for an element of
  $\name{\bbQ}_\xi$ such that   
\[\begin{array}{ll}
r^-_\delta\rest\xi\forces_{\bbP_\xi}&\mbox{`` }r^-_\delta(\xi)\mbox{ is an
upper bound of }\{r_\alpha(\xi):\alpha<\delta\}\cup\{p(\xi)\}\mbox{ and }\\  
&\mbox{\quad if }t\in T_\delta,\ \ \rk_\delta(t)>\xi,\mbox{ and }
q^\delta_{*,t}\rest\xi\in\Gamma_{\bbP_\xi}\mbox{ and the set}\\
&\quad\{r_\alpha(\xi):\alpha<\delta\}\cup\{q^\delta_{*,t}(\xi),
  p(\xi)\}\mbox{ has an upper bound in }\name{\bbQ}_\xi,\\ 
&\mbox{\quad then $r^-_\delta(\xi)$ is such an upper bound '',}
  \end{array}\]
and $r_\delta(\xi)$ is the $<^*_\chi$--first $\bbP_\xi$--name for an element 
  of $\name{\bbQ}_\xi$ such that 
\[\begin{array}{ll}
r_\delta\rest\xi\forces_{\bbP_\xi}&\mbox{`` }r_\delta(\xi)\mbox{ is given to
  Complete by }\name{\st}^0_\xi\mbox{ as the answer to }\ \\ 
&\quad\langle r^-_\alpha(\xi),r_\alpha(\xi):\alpha<\delta\rangle\conc
  \langle r^-_\delta(\xi)\rangle \mbox{  ''} 
  \end{array}\]
\end{itemize}
Note that by a straightforward induction on $\xi\in\dom(r_\delta)$ one
easily applies $(*)_3$ from previous stages to show that
$r_\delta^-,r_\delta$ are well defined and $r_\delta\geq r^-_\delta\geq
r_\alpha,p$ for $\alpha<\delta$. If $\delta=0$ we also stipulate
$r^-_0(0)=r_0(0)=p(0)$.    

This completes the description of stage $\delta$ of our construction. One
easily verifies that the demands $(*)_0$--$(*)_{14}$ are satisfied.

After completing all $\lambda$ stages of the construction, for each $\xi\in
N\cap\gamma$ we look at the sequence $\langle\ell_{\delta,\xi},
\name{\bar{p}}_{\delta,\xi},\name{\bar{q}}_{\delta,\xi}:\delta<\lambda
\rangle$. For $\delta<\lambda$ such that $\xi\in w_\delta$ let 
\[B^\xi_\delta=\bigcap\big\{\lambda\setminus A_{\zeta,\xi}:\zeta \in
w_\delta\cap\xi\big\}\cap \bigcap\big\{A_{\xi,\zeta}:\zeta\in w_\delta
\setminus (\xi+1)\big\},\] 
and for $\delta<\lambda$ such that $\xi\notin w_\delta$ put
$B^\xi_\delta=\lambda$. It follows from our assumptions that
$\forces_{\bbP_\xi} (\forall\delta<\lambda)(B^\xi_\delta\in\name{\cU}_\xi)$
and thus also  $\forces_{\bbP_\xi}\mathop{\triangle}\limits_{\alpha<
  \lambda} B^\xi_\alpha\in\name{\cU}_\xi$. Note that if $\delta$ is a limit 
ordinal, $\xi\in w_\delta$ and $\delta\in \mathop{\triangle}\limits_{\alpha
  <\lambda} B^\xi_\alpha$, then also $\delta\in B^\xi_\delta$ and hence
$\xi=\alpha_\delta$ (remember $(*)_9$) and $\ell_{\delta,\xi}=1$ (by
$(*)_{10}$). Consequently for each $\xi\in N\cap\gamma$ 
\[\forces_{\bbP_\xi}\mbox{`` }\{\delta<\lambda:\ell_{\delta,\xi}=1\} \in
\name{\cU}_\xi\mbox{ ''.}\]
Therefore, for every $\xi\in N\cap\gamma$ we may pick a $\bbP_\xi$--name
$q(\xi)$ for a condition in $\name{\bbQ}_\xi$ such that 
\begin{itemize}
\item if $\xi\in w_{\beta+1}\setminus w_\beta$, $\beta<\lambda$ (or
  $\xi=0=\beta$), then   
\[\forces_{\bbP_\xi}\mbox{`` }q(\xi)\geq r_{\beta}(\xi)\mbox{ and }
q(\xi) \forces_{\name{\bbQ}_\xi}\big\{\delta\!<\!\lambda\!:\big(\exists 
\vare\!<\!\mu_\delta\big)\big(\name{\bar{q}}_{\delta,\xi}(\vare)\in
\Gamma_{\name{\bbQ}_\xi}\big)\big\}\in D[\bbP_{\xi+1}]\mbox{ ''.}\] 
\end{itemize}
This determines a condition $q\in\bbP_\gamma$ (with $\dom(q)=N\cap \gamma$)
and easily $(\forall\beta<\lambda)(p\leq r_\beta\leq q)$ (remember
$(*)_2$). For each $\xi\in N\cap\gamma$ fix $\bbP_{\xi+1}$--names
$\name{C}^\xi_i,\name{f}^\xi_i$ (for $i<\lambda$) such that 
\[\begin{array}{ll}
q\rest(\xi+1)\forces_{\bbP_{\xi+1}}&\mbox{`` }\big(\forall i<\lambda\big)
\big(\name{C}^\xi_i\in D\cap \bV\ \&\ \name{f}^\xi_i\in {}^\lambda\lambda
\big)\mbox{ and}\\
&\ \ \ \big(\forall\delta\in \mathop{\triangle}\limits_{i<\lambda}
\name{C}^\xi_i\cap \mathop{\triangle}\limits_{i<\lambda}\set^\bp(
\name{f}^\xi_i)\big)\big(\exists\vare<\mu_\delta\big)\big(
\name{\bar{q}}_{\delta,\xi}(\vare)\in\Gamma_{\name{\bbQ}_\xi}\big)\mbox{
  ''.} \end{array}\]

\begin{claim}
\label{cl3}
For each limit ordinal $\delta<\lambda$, the condition $q$ forces (in
$\bbP_\gamma$) that 
\[\mbox{`` }\big(\forall\xi\in w_\delta\big)\big(
\delta\in\mathop{\triangle}\limits_{i<\lambda}\name{C}^\xi_i\cap
\mathop{\triangle}\limits_{i<\lambda}\set^\bp(\name{f}^\xi_i)\big)\
\Rightarrow\ \big(\exists t\in T_\delta\big)\big(\rk_\delta(t)=\gamma\ \&\ 
q^\delta_{*,t}\in \Gamma_{\bbP_\gamma}\big)\mbox{ ''.}\]
\end{claim}

\begin{proof}[Proof of the Claim]
The proof is essentially the same as that for \cite[Claim 3.1]{RoSh:860},
however for the sake of completeness we will present it fully. Suppose that
$r\geq q$ and a limit ordinal $\delta<\lambda$ are such that  
\begin{enumerate}
\item[$(\boxplus)_a$] \quad $r\forces_{\bbP_\gamma}$`` $\big(\forall\xi\in 
  w_\delta\big)\big(\delta\in\mathop{\triangle}\limits_{i<\lambda}
  \name{C}^\xi_i\cap \mathop{\triangle}\limits_{i<\lambda}\set^\bp(
  \name{f}^\xi_i)\big)$ ''.  
\end{enumerate}
For each $\zeta<\gamma$ fix a $\bbP_\zeta$--name $\name{\st}^*_\zeta$ for a
winning strategy of Complete in $\Game^\lambda_0(\name{\bbQ}_\zeta,
\name{\emptyset}_{\name{\bbQ}_\zeta})$ such that as long as Incomplete plays
$\name{\emptyset}_{\name{\bbQ}_\zeta}$, Complete answers with
$\name{\emptyset}_{\name{\bbQ}_\zeta}$ as well. 
 
We are going to show that there is $t\in T_\delta$ such that $\rk_\delta(t)
=\gamma$ and the conditions $q^\delta_{*,t}$ and $r$ are compatible. Let
$\langle\vare_\alpha:\alpha\leq\alpha^*\rangle=w_\delta\cup\{\gamma\}$ be
the increasing enumeration. By induction on $\alpha\leq\alpha^*$ we will
choose conditions $r^*_\alpha,r^{**}_\alpha\in\bbP_{\vare_\alpha}$ and
$t=\langle (t)_{\vare_\alpha}:\alpha<\alpha^*\rangle\in T_\delta$ such that
letting $t^\alpha_\circ=\langle (t)_{\vare_\beta}:\beta<\alpha\rangle\in
T_\delta$ we have
\begin{enumerate}
\item[$(\boxplus)_b$] $q^\delta_{*,t^\alpha_\circ}\leq r^*_\alpha$ and 
$r\rest\vare_\alpha\leq r^*_\alpha$,  
\item[$(\boxplus)_c$] for every $\beta<\alpha$ and $\zeta<\vare_\alpha$,
\[\begin{array}{ll}
r^{**}_\beta\forces_{\bbP_\zeta}&\mbox{`` }\langle r^*_{\beta'}(\zeta),
r^{**}_{\beta'}(\zeta):\beta'<\beta\rangle\mbox{ is a legal partial play of 
}\Game^\lambda_0(\name{\bbQ}_\zeta,\name{\emptyset}_{\name{\bbQ}_\zeta})\\
&\ \ \mbox{ in which Complete uses her winning strategy $\name{\st}^*_\zeta$
  ''.}\end{array}\]   
\end{enumerate}
Suppose that $\alpha\leq\alpha^*$ is a limit ordinal and we have already
defined $t^\alpha_\circ=\langle (t)_{\vare_\beta}: \beta<\alpha\rangle$ and 
$\langle r^*_\beta, r^{**}_\beta: \beta<\alpha\rangle$. Let $\xi=\sup(
\vare_\beta:\beta<\alpha)$. It follows  from $(\boxplus)_c$ that we may find
a condition $r^*_\alpha\in\bbP_{\vare_\alpha}$ such that $r^*_\alpha\rest
\xi\in\bbP_\xi$ is stronger than all $r^{**}_\beta$ (for $\beta<\alpha$)
and also $r^*_\alpha\rest [\xi,\vare_\alpha)=r\rest [\xi,\vare_\alpha)$. 
Clearly  $r\rest\vare_\alpha\leq r^*_\alpha$ and also $q^\delta_{*,
t^\alpha_\circ}\rest\xi\leq r^*_\alpha\rest\xi$ (remember $(\boxplus)_b$ for
$\beta<\alpha$). Now by induction on $\zeta\in [\xi,\vare_\alpha)$ we argue
that $q^\delta_{*,t^\alpha_\circ}\rest\zeta\leq r^*_\alpha\rest\zeta$. So
suppose that $\xi\leq \zeta<\vare_\alpha$ and we know $q^\delta_{*,
t^\alpha_\circ}\rest\zeta\leq r^*_\alpha\rest\zeta$. It follows from
$(*)_3+(*)_5+(*)_6$ that $r^*_\alpha\rest\zeta\forces (\forall 
i<\delta)(r_i(\zeta)\leq p^\delta_{*,t^\alpha_\circ}(\zeta)\leq
q^\delta_{*,t^\alpha_\circ}(\zeta))$ and therefore we may use $(*)_{14}$ to
conclude that   
\[r^*_\alpha\rest\zeta\forces_{\bbP_\zeta} q^\delta_{*,t^\alpha_\circ}
(\zeta)\leq r_\delta(\zeta)\leq q(\zeta)\leq r(\zeta)=r^*_\alpha(\zeta).\]  
Finally we let $r^{**}_\alpha\in\bbP_{\vare_\alpha}$ be a condition such
that for each $\zeta<\vare_\alpha$
\[\begin{array}{ll}
r^{**}_\alpha\rest\zeta\forces_{\bbP_\zeta}&\mbox{`` $r^{**}_\alpha(\zeta)$
  is given to Generic by $\name{\st}^*_\zeta$ as the answer to}\\
&\ \ \ \langle r^*_\beta(\zeta),r^{**}_\beta(\zeta):\beta<\zeta\rangle\conc
\langle r^*_\alpha(\zeta)\rangle\mbox{ ''.}
  \end{array}\]
Now suppose that $\alpha=\beta+1\leq\alpha^*$ and we have already defined
$r^*_\beta, r^{**}_\beta\in\bbP_{\vare_\beta}$ and $t^\beta_\circ\in
T_\delta$. It follows from the choice of $q$ and $(\boxplus)_a+(\boxplus)_b
+ (*)_{11}$ that   
\[r^{**}_\beta\forces_{\bbP_{\vare_\beta}}\mbox{`` }r(\vare_\beta)
\forces_{\name{\bbQ}_{\vare_\beta}}\big(\exists\vare<\mu_\delta\big )\big(
q^\delta_{*,t^\beta_\circ\conc\langle\vare\rangle}(\vare_\beta)\in
\Gamma_{\name{\bbQ}_{\vare_\beta}}\big)\mbox{ ''.}\]  
Therefore we may choose $\vare=(t)_{\vare_\beta}<\mu_\delta$ (thus defining
$t^\alpha_\circ$) and a condition $r^*_\alpha\in\bbP_{\vare_\alpha}$ such
that 
\begin{itemize}
\item $r^{**}_\beta\leq r^*_\alpha\rest\vare_\beta$ and 
\[r^*_\alpha\rest\vare_\beta\forces_{\bbP_{\vare_\beta}}\mbox{`` }
r^*_\alpha(\vare_\beta)\geq r(\vare_\beta)\ \&\ r^*_\alpha(\vare_\beta)\geq
q^\delta_{*,t^\alpha_\circ}(\vare_\beta)\mbox{ '',}\]
\item $r^*_\alpha\rest (\vare_\beta,\vare_\alpha)=r\rest (\vare_\beta,
  \vare_\alpha)$. 
\end{itemize}
Exactly like in the limit case we argue that $r^*_\alpha$ has the desired
properties and then in the same manner as there we define $r^{**}_\alpha$.

We finish the proof of the claim noting that $t=t^{\alpha^*}_\circ\in
T_\delta$ and the condition $r^*_{\alpha^*}$ are such that
$r^*_{\alpha^*}\geq r$ and $r^*_{\alpha^*}\geq q^\delta_{*,t}$.   
\end{proof}

Let us use \ref{cl3} to argue that $q$ is $(N,\bbP_\gamma)$--generic. To
this end suppose $\name{\tau}\in N$ is a $\bbP_\gamma$--name for an ordinal,
say $\name{\tau}=\name{\tau}_\alpha$, $\alpha<\lambda$, and let $q'\geq
q$. Since $\bbP_\gamma$ is strategically $({<}\lambda)$--complete we may
build an increasing sequence $\langle q'_i:i<\lambda\rangle$ of conditions
above $q'$ and a sequence $\langle C^\xi_i,f^\xi_i:\xi\in N\cap\gamma,\
i<\lambda\rangle$ such that $C^\xi_i\in D\cap\bV$, $f^\xi_i\in {}^\lambda
\lambda$ and for each $\xi\in w_i$ 
\[q'_i\forces\big(\forall j\leq i\big)\big(\name{C}^\xi_j=C^\xi_j\ \&\
\name{f}^\xi_j\rest i=f^\xi_j\rest i\big).\]
The set $\big\{\delta<\lambda:\big(\forall\xi\in w_\delta\big)\big(\forall
j<\delta\big)\big(\delta\in C^\xi_j\cap\set^\bp(f^\xi_j)\big)\big\}$ is in
$D$, so we may choose a limit ordinal $\delta>\alpha$ such that for each
$\xi\in w_\delta$ we have 
\[\delta\in \mathop{\triangle}\limits_{i<\lambda} C^\xi_i\cap
\mathop{\triangle}\limits_{i<\lambda} \set^\bp(f^\xi_i).\]
Then $q'_\delta\forces\big(\forall\xi\in w_\delta\big)\big(\delta\in
\mathop{\triangle}\limits_{i<\lambda}\name{C}^\xi_i\cap
\mathop{\triangle}\limits_{i<\lambda}\set^\bp(\name{f}^\xi_i)\big)$ and
therefore, by \ref{cl3}, 
\[q'_\delta\forces\big(\exists t\in T_\delta\big)\big(\rk_\delta(t)=\gamma\
\& \ q^\delta_{*,t}\in\Gamma_{\bbP_\gamma}\big).\]
Using $(*)_{13}$ we conclude $q'_\delta\forces\name{\tau}\in Z_\delta$ and
hence $q'_\delta\forces\name{\tau}\in N$.
\end{proof}

\begin{remark}
  Naturally, we want to apply Theorem \ref{lordsthm} to $\gamma=
  \lambda^{++}$ in a model where $2^\lambda=\lambda^+$, so one may ask if
  the assumptions $(3)+(4)$ of \ref{lordsthm} can be satisfied in such a
  universe. But they are not so unusual: suppose that we start with
  $\bV\models \lambda^{<\lambda}=\lambda\ \&\ 2^\lambda=\lambda^+$. Consider
  the following forcing notion $\bbC^\lambda_{\lambda^+}$.

\noindent {\bf A condition $p\in \bbC^\lambda_{\lambda^+}$} is a function
$p:\dom(p)\longrightarrow 2$ such that $\dom(p)\subseteq \lambda^+\times
\lambda$ and $|\dom(p)|<\lambda$. 

\noindent{\bf The order} is the inclusion.\\
Plainly, $\bbC^\lambda_{\lambda^+}$ is a $({<}\lambda)$--complete
$\lambda^+$--cc forcing notion of size $\lambda^+$. Suppose now that
$G\subseteq\bbC^\lambda_{\lambda^+}$ is generic over $\bV$ and let us work
in $\bV[G]$. Put $f=\bigcup G$ (so $f:\lambda^+\times\lambda\longrightarrow
2$) and for $\alpha<\lambda^+$ and $i<2$ let $A^i_\alpha=\{\xi<\lambda: 
f(\alpha,\xi)=i\}$. For a function $h:\lambda^+\longrightarrow 2$ let $U_h$ be
the normal filter generated by the family $\big\{A^{h(\alpha)}_\alpha:
\alpha< \lambda^+\big\}$. One easily verifies that each $U_h$ is a proper
(normal) filter and plainly if $h,h':\lambda^+\longrightarrow 2$ are
distinct, say $h(\alpha)=0$, $h'(\alpha)=1$, then $A^0_\alpha\in U_h$ and
$\lambda\setminus A^0_\alpha=A^1_\alpha\in U_{h'}$. 
\end{remark}

\section{Noble iterations}
The iteration theorems \ref{like860} and \ref{lordsthm} have one common
drawback: they assume that $\lambda$ is strongly inaccessible. In this
section we introduce a property slightly stronger than {\em being
  B--bounding over $\bp$} and we show the corresponding iteration
theorem. The main gain is that the only assumption on $\lambda$ is
$\lambda=\lambda^{<\lambda}$. 

\begin{definition}
\label{supergame}
Let $\bbQ=(\bbQ,\leq)$ be a forcing notion and $\bp=(\bar{P},\lambda,D)$ be
a $\Dl$--parameter on $\lambda$.  
\begin{enumerate}
\item For a condition $p\in\bbQ$ we define a game $\supgame(p,\bbQ)$ between  
two players, Generic and Antigeneric, as follows. A play of $\supgame(p, 
\bbQ)$ lasts $\lambda$ steps during which the players construct a sequence
$\big\langle f_\alpha,\cX_\alpha,\bar{p}_\alpha,\bar{q}_\alpha:\alpha<
\lambda\big\rangle$ such that 
\begin{enumerate}
\item[(a)] $f_\alpha:\alpha\longrightarrow\bbQ$ and $f_\beta\subseteq
  f_\alpha$ for $\beta<\alpha$, 
\item[(b)] $\cX_\alpha\subseteq P_\alpha$ and for every $\eta\in \cX_\alpha$
  the sequence $\langle f_\alpha\big(\eta(\xi)\big):\xi< \alpha \rangle
  \subseteq\bbQ$ has an upper bound in $\bbQ$ and if
  $\eta_0,\eta_1\in\cX_\alpha$ are distinct, then for some $\xi<\alpha$ 
  the conditions $f_\alpha\big(\eta_0(\xi)\big),f_\alpha\big(\eta_1(\xi)
  \big)$ are incompatible,  
\item[(c)] $\bar{p}_\alpha=\langle p^\eta_\alpha:\eta\in\cX_\alpha\rangle
  \subseteq\bbQ$ is a system of conditions in $\bbQ$ such that
  $\big(\forall\xi<\alpha\big)\big(f_\alpha\big(\eta(\xi)\big)
  \leq p^\eta_\alpha\big)$ for $\eta\in\cX_\alpha$, 
\item[(d)] $\bar{q}_\alpha=\langle q^\eta_\alpha:\eta\in\cX_\alpha\rangle
  \subseteq\bbQ$ is a system of conditions in $\bbQ$ such that
  $(\forall\eta\in \cX_\alpha)(p^\eta_\alpha\leq q^\eta_\alpha)$  
\end{enumerate}
The choices of the objects listed above are done so that at stage
$\alpha<\lambda$ of the play:  
\begin{enumerate}
\item[$(\aleph)_\alpha$]  first Generic picks a function $f_\alpha: \alpha
  \longrightarrow\bbQ$ with the property described in (a) above (so if
  $\alpha$ is limit, then $f_\alpha=\bigcup\limits_{\beta<\alpha}
  f_\beta$). She also chooses $\cX_\alpha$, $\bar{p}_\alpha$ satisfying the 
  demands of (b)+(c) (note that $\cX_\alpha$ could be empty). 
\item[$(\beth)_\alpha$]  Then Antigeneric decides a system
  $\bar{q}_\alpha$ as in (d).
\end{enumerate}
At the end, Generic wins the play $\big\langle f_\alpha,\cX_\alpha,
\bar{p}_\alpha,\bar{q}_\alpha: \alpha<\lambda\big\rangle$ if and only if   
\begin{enumerate}
\item[$(\circledast)^{B+}_{\bp}$] there is a condition $p^*\in\bbQ$ stronger 
  than $p$ and such that 
\[p^*\forces_{\bbQ}\mbox{`` }\big\{\alpha<\lambda:\big(\exists \eta\in
\cX_\alpha\big)\big(q^\eta_\alpha\in\Gamma_\bbQ\big)\big\}\in D[\bbQ]\mbox{ 
  ''.}\] 
\end{enumerate}
\item A forcing notion $\bbQ$ is {\em B--noble over $\bp$\/} if it is
  strategically $({<}\lambda)$--complete and Generic has a winning strategy
  in the game $\supgame(p,\bbQ)$ for every $p\in\bbQ$.       
\end{enumerate}
\end{definition}

Note that in the above definition we assumed that $\bar{P}=\langle P_\delta:
\delta<\lambda\rangle$. This was caused only to simplify the description of
the game -- if the domain of $\bar{P}$ is $S\in D$, then we may extend it to 
$\lambda$ in some trivial way without changing the resulting properties.

\begin{observation}
\label{nobleisbounding}
  If $\bp$ is a $\Dl$--parameter and a forcing notion $\bbQ$ is B--noble
  over $\bp$, then $\bbQ$ is reasonably B--bounding over $\bp$.
\end{observation}

\begin{theorem}
\label{nobleiteration}
Assume that 
\begin{enumerate}
\item $\lambda=\lambda^{<\lambda}$ and $\bp=(\bar{P},\lambda,D)$ is a
  $\Dl$--parameter on $\lambda$, and   
\item $\bar{\bbQ}=\langle\bbP_\xi,\name{\bbQ}_\xi:\xi<\gamma\rangle$ is a
  $\lambda$--support iteration such that for every $\xi<\gamma$,    
\[\forces_{\bbP_\xi}\mbox{`` $\name{\bbQ}_\xi$ is B--noble over $\bp[
  \bbP_\xi]$ ''.}\]  
\end{enumerate}
Then 
\begin{enumerate}
\item[(a)] $\bbP_{\gamma}=\lim(\bar{\bbQ})$ is $\lambda$--proper, 
\item[(b)] for each $\bbP_\gamma$--name $\name{\tau}$ for a function from
  $\lambda$ to $\bV$ and a condition $p\in \bbP_\gamma$ there are $q\in
  \bbP_\gamma$ and $\langle A_\alpha:\alpha<\lambda\rangle$ such that
  $|A_\alpha|<\lambda$ (for $\alpha<\lambda$) and $q\geq p$ and 
\[q\forces_{\bbP_\gamma}\mbox{`` }\{\alpha<\lambda:\name{\tau}(\alpha)\in
  A_\alpha\}\in D[\bbP_\gamma]\mbox{ ''.}\]   
\end{enumerate}
\end{theorem}

\begin{proof}
(a) \quad Assume that $N\prec (\cH(\chi),\in,<^*_\chi)$ is such that
  ${}^{<\lambda}N\subseteq N$, $|N|=\lambda$ and $\bar{\bbQ},\bp,\ldots\in
  N$. Let $p\in N\cap\bbP_\gamma$. Choose $\bar{N}=\langle N_\delta:\delta<
  \lambda\rangle$ and $\bar{\alpha}=\langle \alpha_\delta:\delta<\lambda
  \rangle$  such that 
\begin{itemize}
\item  $\bar{N}$ is an increasing continuous sequence of elementary
    submodels of $N$,
\item  $\bar{\alpha}$ is an increasing continuous sequence of ordinals below
  $\lambda$, 
\item $N=\bigcup\limits_{\delta<\lambda} N_\delta$, $\bar{\bbQ},\bp,p,
  \ldots\in N_0$, $\delta\subseteq N_\delta$, $P_\delta\subseteq
  N_{\delta+1}$, $\bar{N}\rest(\delta+1)\in N_{\delta+1}$,
  $|N_\delta|<\lambda$ and 
\item $\alpha_\delta+\otp(N_\delta\cap\gamma)+888<\alpha_{\delta+1}$,
  $\bar{\alpha} \rest (\delta+2)\in N_{\delta+1}$.
  \end{itemize}
Put $w_\delta=N_\delta\cap\gamma$ and for each $\xi<\gamma$ let
$\name{\st}^0_\xi$ be the $<^*_\chi$--first $\bbP_\xi$--name for a winning
strategy of Complete in $\Game_0^\lambda(\name{\bbQ}_\xi,
\name{\emptyset}_{\name{\bbQ}_\xi})$ such that it instructs Complete to play
$\name{\emptyset}_{\name{\bbQ}_\xi}$ as long as her opponent plays
$\name{\emptyset}_{\name{\bbQ}_\xi}$. We also assume that whenever possible,
$\name{\emptyset}_{\name{\bbQ}_\xi}$ is the $<^*_\chi$--first name for the
answer by $\name{\st}^0_\xi$ to a particular sequence of names. Note that
$\langle\name{\st}^0_\xi:\xi<\gamma\rangle\in N_0$.  
\medskip

By induction on $\delta<\lambda$ we will construct
\begin{enumerate}
\item[$(\otimes)_\delta$] $\cT_\delta,\bar{p}^\delta,\bar{q}^\delta,
  r^-_\delta, r_\delta$ and $\name{f}_{\delta,\xi},\name{\cX}_{\delta,\xi},
  \name{\bar{p}}_{\delta,\xi}, \name{\bar{q}}_{\delta,\xi},\name{\st}_\xi$
  and $\bar{p}^{\delta,\xi}$ for $\xi\in w_\delta$    
\end{enumerate}
so that the following conditions $(*)_0$--$(*)_{16}$ are satisfied. 
\begin{enumerate}
\item[$(*)_0$] Objects listed in $(\otimes)_\delta$ form the
  $<^*_\chi$--first tuple with the properties described in
  $(*)_1$--$(*)_{16}$ below. Consequently, the sequence 
\[\langle\mbox{ objects listed in }(\otimes)_\vare:\vare<\delta\rangle\]
is definable from $\bar{N}\rest\alpha_{\delta}$, $\bar{\alpha}\rest\delta$,
$\bar{\bbQ}$, $\bp$, $p$ (in the language $\cL(\in,<^*_\chi)$), so if
$\delta=\alpha_\delta$ is limit, then this sequence belongs to
$N_{\delta+1}$. Also, objects listed in $(\otimes)_\delta$ are known after
stage $\delta$ (and they all belong to $N$).
\item[$(*)_1$] $r^-_\delta,r_\delta\in \bbP_\gamma$, $w_\delta\subseteq
  \dom(r^-_\delta)=\dom(r_\delta)$ and $r^-_0(\xi)=r_0(\xi)=p(\xi)$ for
  $\xi\in w_0$. 
\item[$(*)_2$] For each $\vare<\delta<\lambda$ we have $(\forall\xi\in
  w_{\vare+1})(r_\vare(\xi)=r_\delta(\xi))$ and $p\leq r_\vare^-\leq r_\vare
  \leq r^-_\delta\leq r_\delta$.  
\item[$(*)_3$] If $\xi\in\dom(r_\delta)\setminus w_\delta$,
  then  
\[\begin{array}{ll}
r_\delta\rest\xi\forces_{\bbP_\xi}&\mbox{`` the sequence }\langle
r^-_\vare(\xi),r_\vare(\xi):\vare\leq\delta\rangle\mbox{ is a legal
  partial play of }\\ 
&\quad\Game_0^\lambda\big(\name{\bbQ}_\xi,
\name{\emptyset}_{\name{\bbQ}_\xi}\big)\mbox{ in which Complete follows
}\name{\st}^0_\xi\mbox{ ''}  
\end{array}\]
and if $\xi\in w_{\delta+1}\setminus w_\delta$, then $\name{\st}_\xi$ is the
$<^*_\chi$--first $\bbP_\xi$--name for a winning strategy of Generic in 
$\supgame(r_\delta(\xi),\name{\bbQ}_\xi)$. (And for $\xi\in w_0$,
$\name{\st}_\xi$ is the $<^*_\chi$--first $\bbP_\xi$--name for a winning 
strategy of Generic in $\supgame(p(\xi),\name{\bbQ}_\xi)$. Note that
$\name{\st}_\xi\in N_0$ for $\xi\in w_0$ and $\name{\st}_\xi\in
N_{\alpha_{\delta+1}}$ for $\xi\in w_{\delta+1}$.)     
\item[$(*)_4$] $\cT_\delta=(T_\delta,\rk_\delta)$ is a $(w_\delta,
  1)^\gamma$--tree, $T_\delta\subseteq\bigcup\limits_{\alpha\leq\gamma}   
\prod\limits_{\xi\in w_\delta\cap\alpha}\big(P_\delta\cup\{*\}\big)$.
(Note that we do not require here that $\cT_\delta$ is standard, so some
chains in $\cT_\delta$ may have no $\vtl$--bounds.)
\item[$(*)_5$] $\bar{p}^\delta=\langle p^\delta_t:t\in T_\delta\rangle$ and
  $\bar{q}^\delta=\langle q^\delta_t:t\in T_\delta\rangle$ are trees of
  conditions, $\bar{p}^\delta\leq\bar{q}^\delta$.    
\item[$(*)_6$] For $t\in T_\delta$ we have that $\dom(p^\delta_t) \supseteq 
  \big(\dom(p)\cup\bigcup\limits_{\alpha<\delta}\dom(r_\alpha)\cup
  w_\delta\big) \cap \rk_\delta(t)$ and for each $\xi\in\dom(p^\delta_t)
  \setminus w_\delta$:  
\[\begin{array}{ll}
p^\delta_t\rest\xi\forces_{\bbP_\xi}&\mbox{`` if the set } \{
r_\vare(\xi):\vare<\delta\}\cup\{p(\xi)\}\mbox{ has an upper bound in
}\name{\bbQ}_\xi,\\  
&\mbox{\quad then $p^\delta_t(\xi)$ is such an upper bound ''.}
  \end{array}\]
\item[$(*)_7$] If $\xi\in w_{\beta+1}\setminus w_\beta$, $\beta<\delta$, 
  then 
\[\begin{array}{r}
\forces_{\bbP_\xi}\mbox{`` }\langle \name{f}_{\vare,\xi},
\name{\cX}_{\vare,\xi}, \name{\bar{p}}_{\vare,\xi},
\name{\bar{q}}_{\vare,\xi}:\vare<\delta \rangle\mbox{ is a partial play of
}\\ 
\supgame(r_\beta(\xi),\name{\bbQ}_\xi)\mbox{ in which Generic uses
  $\name{\st}_\xi$ ''.} 
  \end{array}\]
\item[$(*)_8$] $\dom(r^-_\delta)=\dom(r_\delta)=\bigcup\limits_{t\in
    T_\delta}\dom(q^\delta_t)$ and if $t\in T_\delta$, $\xi\in\dom(r_\delta)
  \cap\rk_\delta(t)\setminus w_\delta$, and $q^\delta_t\rest\xi\leq
  q\in\bbP_\xi$, $r_\delta\rest\xi\leq q$, then    
\[\begin{array}{ll}
q\forces_{\bbP_\xi}&\mbox{`` if the set }\{r_\alpha(\xi):\alpha<\delta\}
\cup\{q^\delta_t(\xi),p(\xi)\}\mbox{ has an upper bound in }
\name{\bbQ}_\xi,\\   
&\mbox{\quad then $r^-_\delta(\xi)$ is such an upper bound ''.}
  \end{array}\]
\item[$(*)_{9}$] $\bar{p}^{\delta,\xi}=\langle p^{\delta,\xi}_t:t\in
  T_\delta\ \&\ \rk_\delta(t)\leq\xi\rangle\subseteq\bbP_\xi$ is a tree of
  conditions (for $\xi\in w_\delta\cup\{\gamma\}$), $\bar{p}^{\delta,\gamma}
  =\bar{p}^\delta$.
\item[$(*)_{10}$] If $\zeta,\xi\in w_\delta\cup\{\gamma\}$, $\zeta<\xi$ and
  $t\in T_\delta$, $\rk_\delta(t)=\zeta$, then $p^{\delta,\zeta}_t\leq
  p^{\delta,\xi}_t$. 
\end{enumerate}
The demands $(*)_{11}$--$(*)_{16}$ formulated below are required only if
$\delta=\alpha_\delta$ is a limit ordinal.

\begin{enumerate}
\item[$(*)_{11}$] If $t\in T_\delta$, $\rk_\delta(t)=\xi\in w_\delta$ and
  $X^\delta_t=\{(s)_\xi:t\vtl s\in T_\delta\}$, then
\begin{itemize}
\item   either $\emptyset\neq X^\delta_t\subseteq P_\delta$ and
  $p^{\delta,\xi}_t\forces_{\bbP_\xi}$`` $\name{\cX}_{\delta, 
    \xi}= X_t^\delta$ '',
\item  or $X^\delta_t=\{*\}$ and $p^{\delta,\xi}_t\forces_{\bbP_\xi}$``
  $\name{\cX}_{\delta,\xi}=\emptyset$ ''.
  \end{itemize}
\item[$(*)_{12}$] If $s\in T_\delta$, $\rk_\delta(s)=\zeta$, $\xi\in
  w_\delta\cap\zeta$ and $(s)_\xi\neq *$, then 
\[p^{\delta,\zeta}_s\rest\xi\forces_{\bbP_\xi}\mbox{`` }
\name{\bar{p}}_{\delta,\xi} ((s)_\xi)\leq p^{\delta,\zeta}_s(\xi)\mbox{
  ''.}\] 
\item[$(*)_{13}$] If $\xi\in w_\delta\cup\{\gamma\}$, $\otp(w_\delta\cap\xi) 
  = \zeta$ then 
\begin{itemize}
\item $\{t\in T_\delta:\rk_\delta(t)\leq\xi\}\subseteq N_{\delta+\zeta+1}$,
  $\{t\in T_\delta:\rk_\delta(t)\leq\xi\}, \langle \bar{p}^{\delta,\beta}:
  \beta\in (w_\delta\cup\{\gamma\})\cap (\xi+1)\rangle\in
  N_{\delta+\zeta+2}$, and  
\item if $\zeta$ is limit, then $\langle \bar{p}^{\delta,\beta}: 
  \beta\in w_\delta\cap \xi\rangle\in N_{\delta+\zeta+1}$, and if $\bar{t}= 
  \langle t_i:i<i^*\rangle\in N_{\delta+\zeta+1}$ is a $\vtl$--chain in
  $\{t\in T_\delta:\rk_\delta(t)<\xi\}$ with $\sup(\rk_\delta(t_i): i<i^*)
  =\sup(w_\delta\cap\xi)$, then $\bar{t}$ has a $\vtl$--bound in $T_\delta$.  
\end{itemize}
\item[$(*)_{14}$]  If $t\in T_\delta$, $\xi\in w_\delta\cap\rk_\delta(t)$
  and $(t)_\xi\neq *$, then  
\[q^\delta_t\rest\xi\forces_{\bbP_\xi}\mbox{`` }\name{\bar{q}}_{\delta,\xi}( 
(t)_\xi)=q^\delta_t(\xi)\mbox{ ''.}\] 
\item[$(*)_{15}$] If $t_0,t_1\in T_\delta$, $\rk_\delta(t_0)=
  \rk_\delta(t_1)$ and $\xi\in w_\delta\cap\rk_\delta(t_0)$, $t_0\rest 
  \xi=t_1\rest\xi$ but $\big(t_0\big)_\xi\neq \big(t_1\big)_\xi$, then 
\[p^\delta_{t_0\rest\xi}\forces_{\bbP_\xi}\mbox{`` the conditions  
$p^\delta_{t_0}(\xi),p^\delta_{t_1}(\xi)$ are incompatible ''.}\] 
\item[$(*)_{16}$] If $\name{\tau}\in N_\delta$ is a $\bbP_\gamma$--name for
  an ordinal and $t\in T_\delta$ satisfies $\rk_\delta(t)=\gamma$, then the
  condition $q^\delta_t$ forces a value to $\name{\tau}$. 
\end{enumerate}
The rule $(*)_0$ (and conditions $(*)_1$--$(*)_{16}$) actually fully
determines our objects, but we should argue that at each stage there exist
objects with properties listed in $(*)_1$--$(*)_{16}$.

Suppose we have arrived to a stage $\delta<\lambda$ of the construction and
all objects listed in $(\otimes)_\beta$ for $\beta<\delta$ have been
determined so that all relevant demands are satisfied, in particular the
sequence $\langle\mbox{ objects listed in }(\otimes)_\vare:\vare<\delta
\rangle$ is definable from $\bar{N}\rest\alpha_\delta$,
$\bar{\alpha}\rest\delta$, $\bar{\bbQ},\bp,p$.  

If $\delta$ is a successor ordinal and $\xi\in w_\delta\setminus
w_{\delta-1}$, then we let $\name{\st}_\xi$ be the $<^*_\chi$--first
$\bbP_\xi$--name for a winning strategy of Generic in $\supgame(r_{\delta
-1}(\xi),\name{\bbQ}_\xi)$. We also pick the $<^*_\chi$--first
sequence $\langle\name{f}_{\vare,\xi},\name{\cX}_{\vare,\xi},
\name{\bar{p}}_{\vare,\xi},\name{\bar{q}}_{\vare,\xi}:\vare<\delta\rangle$
so that $(*)_7$ is satisfied. Then assuming that $\delta$ is not limit or
$\delta\neq\alpha_\delta$ we may find objects listed in $(\otimes)_\delta$
so that the demands in $(*)_1$--$(*)_{10}$ are satisfied and $|\{t\in
T_\delta:\rk_\delta(t)=\gamma\}|=1$.  

So suppose now that $\delta=\alpha_\delta$ is a limit ordinal. For each
$\xi\in w_\delta$ we let $\name{f}_{\delta,\xi}$ be the $<^*_\chi$--first
$\bbP_\xi$--name such that  $\forces_{\bbP_\xi}$``$\name{f}_{\delta,\xi} = 
\bigcup\limits_{\alpha< \delta}\name{f}_{\alpha,\xi}$'', and
$\name{\cX}_{\delta,\xi}, \name{\bar{p}}_{\delta,\xi}$ be the
$<^*_\chi$--first $\bbP_\xi$--names such that   
\[\begin{array}{ll}
\forces_{\bbP_\xi}&\mbox{`` }\name{f}_{\delta,\xi},\name{\cX}_{\delta,\xi}, 
\name{\bar{p}}_{\delta,\xi} \mbox{ are given to Generic by }\name{\st}_\xi\\ 
&\ \ \mbox{ as the answer to }\langle\name{f}_{\vare,\xi},
\name{\cX}_{\vare,\xi},\name{\bar{p}}_{\vare,\xi},
\name{\bar{q}}_{\vare,\xi}: \vare<\delta\rangle\mbox{ ''.}\end{array}\] 
Note that 
\[\langle\mbox{ objects listed in }(\otimes)_\vare:\vare<\delta
\rangle\conc\langle \name{f}_{\delta,\xi},\name{\cX}_{\delta,\xi},
\name{\bar{p}}_{\delta,\xi}: \xi\in w_\delta\rangle\in N_{\delta+1}.\] 
Now by induction on $\xi\in w_\delta\cup\{\gamma\}$ we will choose $\{t\in
T_\delta: \rk_\delta(t)\leq\xi\}$ and $\bar{p}^{\delta,\xi}$ and auxiliary
objects $\bar{p}^{*,\xi}$ so that, in addition to demands
$(*)_{9}$--$(*)_{13}$ we also have:
\begin{enumerate}
\item[$(*)_{17}$] $\bar{p}^{*,\xi}=\langle p^{*,\xi}_t:t\in T_\delta\ \&\
  \rk_\delta(t)\leq \xi\rangle\subseteq\bbP_\xi$ is a tree of conditions,
  $\bar{p}^{*,\xi}\leq \bar{p}^{\delta,\xi}$ and $\dom(p^{*,\xi}_t)\supseteq 
  \big(\dom(p)\cup \bigcup\limits_{\vare<\delta}\dom(r_\vare)\cup w_\delta
  \big)\cap\rk_\delta(t)$ whenever $t\in T_\delta$, $\rk_\delta(t)\leq \xi$,
  and 
\item[$(*)_{18}$] if $\xi_0<\xi_1$ are from $w_\delta\cup\{\gamma\}$, $t\in
  T_\delta$, $\rk_\delta(t)=\xi_0$, then $p^{\delta,\xi_0}_t\leq
  p^{*,\xi_1}_t$, and 
\item[$(*)_{19}$] if $t\in T_\delta$, $\rk_\delta(t)=\xi$, then
  $\dom(p^{*,\xi}_t)=\dom(p^{\delta,\xi}_t)$ and for $\beta\in
  \dom(p^{\delta,\xi}_t)$ we have
\[\begin{array}{ll}
p^{\delta,\xi}_t\rest\beta\forces_{\bbP_\beta}&\mbox{``  the sequence
}\langle p^{*,\zeta}_{t\rest\zeta}(\beta),p^{\delta,\zeta}_{t\rest\zeta}(
\beta):\zeta\in (w_\delta\cup\{\gamma\})\cap (\xi+1)\rangle\mbox{ is }\\ 
&\ \mbox{ a legal partial play of }\Game^\lambda_0(\name{\bbQ}_\beta, 
\name{\emptyset}_{\name{\bbQ}_\beta})\mbox{ in which Complete uses }\\
&\ \mbox{ the winning strategy }\name{\st}^0_\beta\mbox{ ''.} 
\end{array}\]
\end{enumerate}
To take care of clause $(*)_{13}$, each time we pick an object, we choose
the $<^*_\chi$--first one with the respective property.
\smallskip

\noindent{\sc Case 1:}\quad $\otp(w_\delta\cap\xi)=\zeta+1$ is a successor
ordinal.\\ 
Let $\xi_0=\max(w_\delta\cap\xi)$ and suppose that we have defined
$T^*=\{t\in T_\delta: \rk_\delta(t)\leq \xi_0\}$ and $\bar{p}^{*,\xi_0}, 
\bar{p}^{\delta,\xi_0}$ satisfying the relevant demands of
$(*)_{9}$--$(*)_{19}$. Let $t\in T^*$ be such that $\rk_\delta(t)=\xi_0$. It
follows from $(*)_{11}$ that either $p^{\delta,\xi_0}_t\forces$``
$\name{\cX}_{\delta,\xi_0}=\emptyset$ '' or $p^{\delta,\xi_0}_t\forces$``
$\name{\cX}_{\delta,\xi_0}=X_t^\delta$ '' for some non-empty set $X_t^\delta
\subseteq P_\delta$. In the former case stipulate $X_t=\{*\}$. Note that
necessarily $X_t^\delta\subseteq N_{\delta+1}$ and $X_t^\delta\in  
N_{\delta+\zeta+2}$ (remember $(*)_{13}$). We declare that 
\[\{t\in T_\delta:\rk_\delta(t)\leq\xi\}=T^*\cup\{t\cup \{(\xi_0,a)\}: t\in
T^*\ \&\ \rk_\delta(t)=\xi_0\ \&\ a\in X_t^\delta\}.\] 
Plainly,  $|\{t\in T_\delta:\rk_\delta(t)\leq\xi\}|<\lambda$ and even
$\{t\in T_\delta:\rk_\delta(t)\leq\xi\}\subseteq N_{\delta+\zeta+1}$ and
$\{t\in T_\delta:\rk_\delta(t)\leq\xi\}\in N_{\delta+\zeta+2}$ (again, by
  $(*)_{13}$). Choose a tree of conditions $\bar{p}^+=\langle p^+_t: t\in
  T_\delta\ \&\ \rk_\delta(t)\leq\xi\rangle\subseteq \bbP_\xi$ so that 
  \begin{itemize}
\item $\dom(p^+_t)\supseteq\big(\dom(p)\cup \bigcup\limits_{\vare<\delta} 
    \dom(r_\vare)\cup w_\delta\big)\cap\rk_\delta(t)$ for $t\in T_\delta$,
    $\rk_\delta(t)\leq \xi$, 
\item if $t\in T_\delta$, $\rk_\delta(t)<\xi$ then
  $p^+_t=p^{\delta,\xi_0}_t$,
\item if $t\in T_\delta$, $\rk_\delta(t)=\xi$ and $(t)_{\xi_0}\neq *$, then
  $p^+_t(\xi_0)$ is a $\bbP_{\xi_0}$--name such that
\[p^{\delta,\xi_0}_{t\rest\xi}\forces_{\bbP_{\xi_0}}\mbox{`` }
\name{\bar{p}}_{\delta,\xi_0}\big((t)_{\xi_0}\big)\leq p^+_t(\xi_0)\mbox{
  '',}\] 
\item if $t\in T_\delta$, $\rk_\delta(t)=\xi$ and $\beta\in \dom(p^+_t)
  \setminus (\xi_0+1)$, then
\[\begin{array}{ll}
p^+_t\rest\beta\forces_{\bbP_\beta}&\mbox{`` if the set } \{
r_\vare(\beta):\vare<\delta\}\cup\{p(\beta)\}\mbox{ has an upper bound in 
}\name{\bbQ}_\beta,\\  
&\mbox{\quad then $p^+_t(\beta)$ is such an upper bound ''.}
  \end{array}\]
\end{itemize}
(Note: $\bar{p}^+\in N_{\delta+\zeta+2}$.) Next we may use Proposition
\ref{oldstuff} to pick a tree of conditions $\bar{p}^{*,\xi}=\langle
p^{*,\xi}_t: t\in T_\delta\ \&\ \rk_\delta(t)\leq\xi\rangle$ such that
$\bar{p}^+\leq\bar{p}^{*,\xi}$ and
\begin{itemize}
\item if $\xi<\gamma$, $t\in T_\delta$, $\rk_\delta(t)=\xi$, then either
  $p^{*,\xi}_t \forces\name{\cX}_{\delta,\xi}=\emptyset$ or for some
  non-empty set $X^\delta_t\subseteq P_\delta$ we have $p^{*,\xi}_t\forces  
  \name{\cX}_{\delta,\xi}= X_t^\delta$.
\end{itemize}
(Again, by our rule of picking ``the $<^*_\chi$--first'',
$\bar{p}^{*,\xi}\in N_{\delta+\zeta+2}$.) Then we choose a tree of
conditions $\bar{p}^{\delta,\xi}=\langle p^{\delta,\xi}_t: t\in T_\delta\
\&\ \rk_\delta(t)\leq\xi\rangle$ so that $\bar{p}^{*,\xi}\leq
\bar{p}^{\delta,\xi}$ and for every $t\in T_\delta$ with $\rk_\delta(t)=\xi$
we have $\dom(p^{\delta,\xi}_t)=\dom(p^{*,\xi}_t)$ and for $\beta\in
\dom(p^{\delta,\xi}_t)$, $p^{\delta,\xi}_t(\beta)$ is the $<^*_\chi$--first
$\bbP_\beta$--name for a condition in $\name{\bbQ}_\beta$ such that 
\[\begin{array}{ll}
p^{\delta,\xi}_t\rest\beta\forces_{\bbP_\beta}&\mbox{``
  $p^{\delta,\xi}_t(\beta)$ is given to Generic by $\name{\st}^0_\beta$ as
  the answer to}\\ 
&\ \ \ \langle p^{*,\vare}_{t\rest\vare}(\beta),
p^{\delta,\vare}_{t\rest\vare}(\beta): \vare\in w_\delta\cap\xi\rangle\conc
\langle  p^{*,\xi}_t(\beta)\rangle\mbox{ ''.} 
  \end{array}\]
Note that, by the rule of picking ``the $<^*_\chi$--first'',
$\bar{p}^{\delta,\xi}\in N_{\delta+\zeta+2}$. It should be also clear that
$\bar{p}^{*,\xi},\bar{p}^{\delta,\xi}$ satisfy all the relevant demands stated
in $(*)_{9}$--$(*)_{19}$. 
\smallskip

\noindent{\sc Case 2:}\quad $\otp(w_\delta\cap\xi)=\zeta$ is a limit
ordinal.\\  
Suppose we have defined $\{t\in T_\delta:\rk_\delta(t)\leq\vare\}$ and
$\bar{p}^{*,\vare},\bar{p}^{\delta,\vare}$ for $\vare\in w_\delta\cap
\xi$. By our rule of choosing ``the $<^*_\chi$--first objects'', we know
that the sequence $\big\langle\{t\in T_\delta:\rk_\delta(t)\leq\vare\},
\bar{p}^{*,\vare}, \bar{p}^{\delta,\vare}: \vare\in w_\delta\cap\xi\big
\rangle$ belongs to $N_{\delta+\zeta+1}$.  We also know that $\{t\in
T_\delta:\rk_\delta(t)<\xi\}\subseteq N_{\delta+\zeta}$. Let $T^+$ be the
set of all limit branches in $\big(\{t\in T_\delta:\rk_\delta(t)<\xi\},\vtl
\big)$, so elements of $T^+$ are sequences $s=\langle (s)_\vare: \vare\in
w_\delta\cap \xi\rangle$ such that $s\rest\vare=\langle
(s)_{\vare'}:\vare'\in w_\delta\cap\vare\rangle\in \{t\in T_\delta:
\rk_\delta(t)\leq\vare\}$ for $\vare\in w_\delta\cap\xi$. (Of course,
$T^+\in N_{\delta+\zeta+1}$.) We put 
\[\{t\in T_\delta:\rk_\delta(t)\leq\xi\}=\{t\in T_\delta:\rk_\delta(t)<
\xi\} \cup\big(T^+\cap N_{\delta+\zeta+1}\big)\in N_{\delta+\zeta +2}.\] 
Due to $(*)_{19}$ at stages $\vare\in w_\delta\cap\xi$, we may choose a tree
of conditions $\bar{p}^+=\langle p^+_t:t\in T_\delta\ \&\ \rk_\delta(t)\leq
\xi\rangle\subseteq\bbP_\xi$ such that 
\begin{itemize}
\item $\dom(p^+_t)\supseteq\big(\dom(p)\cup \bigcup\limits_{\vare<\delta} 
    \dom(r_\vare)\cup w_\delta\big)\cap\rk_\delta(t)$ for $t\in T_\delta$,
    $\rk_\delta(t)\leq \xi$, and 
\item if $t\in T_\delta$, $\rk_\delta(t)=\xi_0<\xi$ then
  $\dom(p^+_t)\supseteq \dom(p^{\delta,\xi_0}_t)$ and for each $\beta\in
  \dom(p^+_t)\cap \xi_0$ we have
\[p^+_t\rest\beta\forces_{\bbP_\beta}\mbox{`` }p^{\delta,\xi_0}_t(\beta)
\leq p^+_t(\beta)\mbox{ '', and}\]
\item if $t\in T_\delta$, $\rk_\delta(t)=\xi$, $\sup(w_\delta\cap\xi)
  \leq\beta<\xi$, $\beta\in\dom(p^+_t)$, then
\[\begin{array}{ll}
p^+_t\rest\beta\forces_{\bbP_\beta}&\mbox{`` if the set } \{
r_\vare(\beta):\vare<\delta\}\cup\{p(\beta)\}\mbox{ has an upper bound in  
}\name{\bbQ}_\beta,\\  
&\mbox{\quad then $p^+_t(\beta)$ is such an upper bound ''.}
  \end{array}\]
\end{itemize}
Then, like in the successor case, we may find a tree of conditions
$\bar{p}^{*,\xi}= \langle p^{*,\xi}_t:t\in T_\delta\ \&\ \rk_\delta(t)\leq
\xi\rangle$ such that $\bar{p}^+\leq\bar{p}^{*,\xi}$ and 
\begin{itemize}
\item if $\xi<\gamma$, $t\in T_\delta$, $\rk_\delta(t)=\xi$, then either
  $p^{*,\xi}_t \forces\name{\cX}_{\delta,\xi}=\emptyset$ or for some
  non-empty set $X^\delta_t\subseteq P_\delta$ we have $p^{*,\xi}_t\forces 
  \name{\cX}_{\delta,\xi}= X_t^\delta$. 
\end{itemize}
Also like in that case we choose $\bar{p}^{\delta,\xi}=\langle
p^{\delta,\xi}_t: t\in T_\delta\ \&\ \rk_\delta(t)\leq\xi\rangle$. Clearly,
all relevant demands in $(*)_{9}$--$(*)_{19}$ are satisfied.
\smallskip

The last stage of the above construction gives us a tree $T_\delta=\{t\in
T_\delta: \rk_\delta(t)\leq\gamma\}$ and a tree of conditions
$\bar{p}^{\delta,\gamma}=\bar{p}^\delta=\langle p^\delta_t:t\in T_\delta
\rangle$. Since $T_\delta\subseteq N_{\delta+\otp(w_\delta)+1}$, we know
that $|T_\delta|<\lambda$ so we may apply Proposition \ref{oldstuff} to get
a tree of conditions $\bar{q}^\delta=\langle q^\delta_t:t\in T_\delta\rangle
\geq\bar{p}^\delta$ such that $(*)_{16}$ is satisfied. Remembering
$(*)_{15}+(*)_{12}$, we easily find $\bbP_\xi$--names
$\name{\bar{q}}_{\delta,\xi}$ (for $\xi\in w_\delta$) such  that   
\begin{itemize}
\item $\forces_{\bbP_\xi}$`` $\name{\bar{q}}_{\delta,\xi}=\langle
  \name{\bar{q}}_{\delta,\xi}(\eta):\eta\in \name{\cX}_{\delta,\xi}\rangle$
  is a system of conditions in $\name{\bbQ}_\xi$ '',
\item $\forces_{\bbP_\xi}$`` $(\forall\eta\in \name{\cX}_{\delta,\xi})
  (\name{\bar{p}}_{\delta,\xi}(\eta)\leq\name{\bar{q}}_{\delta,\xi}(\eta))$  
  '', and 
\item if $t\in T_\delta$, $\rk_\delta(t)>\xi$, then $q^\delta_{t\rest\xi}
  \forces_{\bbP_\xi}$`` $\name{\cX}_{\delta,\xi}\neq\emptyset\ \Rightarrow\
  q^\delta_t(\xi)=\name{\bar{q}}_{\delta,\xi}\big((t)_\xi\big)$ ''.  
\end{itemize}
So then $(*)_{14}$ is satisfied. Now we define
$r^-_\delta,r_\delta\in\bbP_\gamma$ essentially by $(*)_1$--$(*)_3$ and
$(*)_8$.  
\medskip

After completing all $\lambda$ stages of the construction, for each $\xi\in
N\cap\gamma$ we look at the sequence $\langle\name{f}_{\alpha,\xi},
\name{\cX}_{\alpha,\xi},\name{\bar{p}}_{\alpha,\xi},
\name{\bar{q}}_{\alpha,\xi}: \alpha<\lambda\rangle$. By $(*)_7$, it is a
$\bbP_\xi$--name for a play of $\supgame(r_\beta(\xi),\name{\bbQ}_\xi)$
(where $\xi\in w_{\beta+1}\setminus w_\beta$) in which Generic uses her
winning strategy $\name{\st}_\xi$. Therefore, for every $\xi\in N\cap\gamma$
we may pick a $\bbP_\xi$--name $q(\xi)$ for a condition in $\name{\bbQ}_\xi$
such that    
\begin{itemize}
\item if $\xi\in w_{\beta+1}\setminus w_\beta$, $\beta<\lambda$ (or
  $\xi\in w_0$, $\beta=0$), then   
\[\forces_{\bbP_\xi}\mbox{``}q(\xi)\geq r_{\beta}(\xi)\mbox{ and }
q(\xi) \forces_{\name{\bbQ}_\xi}\big\{\delta\!<\!\lambda\!:\big(\exists 
\eta\in\name{\cX}_{\delta,\xi}\big)\big(\name{\bar{q}}_{\delta,\xi}(\eta)\in 
\Gamma_{\name{\bbQ}_\xi}\big)\big\}\in D^{\bp}[\bbP_{\xi+1}]\mbox{''.}\]  
\end{itemize}
This determines a condition $q\in\bbP_\gamma$ (with $\dom(q)=N\cap \gamma$) 
and easily $(\forall\beta<\lambda)(p\leq r_\beta\leq q)$ (remember
$(*)_2$). For each $\xi\in N\cap\gamma$ fix $\bbP_{\xi+1}$--names 
$\name{C}^\xi_i,\name{g}^\xi_i$ (for $i<\lambda$) such that
\[\begin{array}{ll}
q\rest(\xi+1)\forces_{\bbP_{\xi+1}}&\mbox{`` }\big(\forall i<\lambda\big)
\big(\name{C}^\xi_i\in D\cap \bV\ \&\ \name{g}^\xi_i\in {}^\lambda\lambda
\big)\mbox{ and}\\
&\ \ \ \big(\forall\delta\in \mathop{\triangle}\limits_{i<\lambda}
\name{C}^\xi_i\cap \mathop{\triangle}\limits_{i<\lambda}\set^\bp(
\name{g}^\xi_i)\big)\big(\exists\eta\in\name{\cX}_{\delta,\xi}\big)\big(
\name{\bar{q}}_{\delta,\xi}(\eta)\in\Gamma_{\name{\bbQ}_\xi}\big)\mbox{
  ''.} \end{array}\]
Let $\name{B}$ be a $\bbP_\gamma$--name for the set $\{\delta<\lambda:
\Gamma_{\bbP_\gamma}\cap N_\delta\in N_{\delta+1}\}$. It follows from Lemma
\ref{Dlmodels} that $\forces_{\bbP_\gamma}\name{B}\in D^\bp[\bbP_\gamma]$.

\begin{claim}
\label{cl4}
If  $\alpha_\delta=\delta$ is limit, then 
\[\begin{array}{ll}
q\forces_{\bbP_\gamma}&\mbox{`` if }\delta\in\name{B}\mbox{ and }\big(
\forall\xi\in w_\delta\big)\big(\delta\in\mathop{\triangle}\limits_{i<
  \lambda}\name{C}^\xi_i\cap \mathop{\triangle}\limits_{i<\lambda} \set^\bp
(\name{g}^\xi_i)\big)\\  
&\quad\mbox{ then }\ \big(\exists t\in T_\delta\big)\big(\rk_\delta(t)=
\gamma\ \&\ q^\delta_t\in \Gamma_{\bbP_\gamma}\big)\mbox{ ''.}
\end{array}\] 
\end{claim}

\begin{proof}[Proof of the Claim]
Suppose that $\delta=\alpha_\delta$ is a limit ordinal and a condition
$r\geq q$ forces (in $\bbP_\gamma$) that 
\begin{enumerate}
\item[$(*)^a_{20}$] $\delta\in\name{B}$, and 
\item[$(*)^b_{20}$] $\big(\forall\xi\in w_\delta\big)\big(\delta\in
  \mathop{\triangle}\limits_{i<\lambda}\name{C}^\xi_i\cap
  \mathop{\triangle}\limits_{i<\lambda}\set^\bp(\name{g}^\xi_i)\big)$. 
\end{enumerate}
Passing to a stronger condition if necessary, we may also assume that 
\begin{enumerate}
\item[$(*)^a_{21}$] if $r'\in\bbP_\gamma\cap N_\delta$, then either $r'\leq
  r$ or $r',r$ are incompatible in $\bbP_\gamma$.
\end{enumerate}
Let $H^\delta=\{r'\in\bbP_\gamma\cap N_\delta:r'\leq r\}$. It follows from
$(*)^a_{20}+(*)^a_{21}$ that
\begin{enumerate}
\item[$(*)^b_{21}$] $r\forces$`` $\Gamma_{\bbP_\gamma}\cap N_\delta=
  H^\delta$ '' and $H^\delta\in N_{\delta+1}$. 
\end{enumerate}
By \ref{supergame}(1)(a)+ $(*)_7+(*)_0$ we may choose a sequence
$\name{\bar{\tau}}=\langle\name{\tau}(\xi,\alpha):\xi\in w_\delta\ \&\
\alpha<\delta\rangle\in N_{\delta+1}$ such that
\begin{itemize}
\item $\name{\tau}(\xi,\alpha)$ is a $\bbP_\xi$--name for an element of
  $\name{\bbQ}_\xi$, $\name{\tau}(\xi,\alpha)\in N_\delta$, 
\item $\forces_{\bbP_\xi}\name{\tau}(\xi,\alpha)=\name{f}_{\delta,\xi}(
  \alpha)$. 
\end{itemize}
Next we choose a sequence $t^*=\langle (t^*)_\xi:\xi\in w_\delta\rangle \in
\prod\limits_{\xi\in w_\delta}\big(P_\delta\cup\{*\}\big)$ so that for each
$\xi\in w_\delta$, $(t^*)_\xi$ is the $<^*_\chi$--first member of
$P_\delta\cup \{*\}$ satisfying:
\begin{enumerate}
\item[$(*)_{22}$] if $t=t^*\rest\xi=\langle (t^*)_\vare:\vare\in w_\delta
  \cap \xi\rangle\in T_\delta$ and
  \begin{enumerate}
\item[(i)] for some non-empty set $X\subseteq P_\delta$, $p^{\delta,\xi}_t
    \forces_{\bbP_\xi} X=\name{\cX}_{\delta,\xi}$ (remember $(*)_{11}$) and
    there is $\eta\in X$ such that 
\item[(ii)] $\big(\forall\alpha<\delta\big)\big(\name{\tau}(\xi,
  \eta(\alpha))\in \{r'(\xi):r'\in H^\delta\}\big)$, 
  \end{enumerate}
then $(t^*)_\xi\in X$ and $\big(\forall\alpha<\delta\big)\big(
  \name{\tau}(\xi,(t^*)_\xi(\alpha))\in \{r'(\xi):r'\in H^\delta\}\big)$. 
\end{enumerate}
Note that for every $\xi\in w_\delta\cup\{\gamma\}$ the sequence $t^*\rest
\xi$ is definable (in $\cL(\in,<^*_\chi)$) from $\bp,\name{\bar{\tau}},
H^\delta, w_\delta,\xi$ and $\langle\bar{p}^{\delta,\vare}:\vare\in
w_\delta\cap\xi\rangle$. Consequently, if $\xi\in w_\delta\cup\{\gamma\}$
and $\zeta=\otp(w_\delta\cap\xi)$, then $t^*\rest\xi\in
N_{\delta+\zeta+1}$. Now, by induction on $\xi\in w_\delta\cup\{\gamma\}$ we
are going to show that $t^*\rest\xi\in T_\delta$ and choose conditions
$r^*_\xi,r^{**}_\xi\in\bbP_\xi$ such that 
\begin{enumerate}
\item[$(*)_{23}^a$] $q^\delta_{t^*\rest\xi}\leq r^*_\xi$, $r\rest\xi \leq
  r^*_\xi$ and if $\vare\in w_\delta\cap\xi$ then $r^{**}_\vare\leq
  r^*_\xi$, and 
\item[$(*)_{23}^b$] $\dom(r^*_\xi)=\dom(r^{**}_\xi)$ and $r^*_\xi\leq
  r^{**}_\xi$ and for every $\beta\in\dom(r^{**}_\xi)$
\[\begin{array}{ll}
r^{**}_\xi\forces_{\bbP_\beta}&\mbox{`` }\langle r^*_\vare(\beta),
r^{**}_\vare(\beta):\vare\in w_\delta\cap(\xi+1)\rangle\mbox{ is a partial
  play of }\Game^\lambda_0(\name{\bbQ}_\beta,
\name{\emptyset}_{\name{\bbQ}_\beta})\\ 
&\mbox{ in which Complete uses her winning strategy }\name{\st}^0_\beta
\mbox{ ''.}
\end{array}\]
\end{enumerate}

Suppose that $\otp(w_\delta\cap\xi)$ is a limit ordinal and for $\vare\in
w_\delta\cap\xi$ we know that $t^*\rest\vare\in T_\delta$  and we have
defined $r^*_\vare,r^{**}_\vare$. It follows from $(*)_{13}$ that $t^*\rest
\xi\in T_\delta$. Let $\beta=\sup(w_\delta\cap\xi)\leq\xi$. It follows from
$(*)^b_{23}$ that we may find a condition $r^*_\xi\in\bbP_\xi$ such that
$r^*_\xi\rest\beta$ is stronger than all $r^{**}_\vare$ (for $\vare\in
w_\delta\cap\xi$) and $r^*_\xi\rest [\beta,\xi)=r\rest [\beta,\xi)$. Clearly
$r\rest\xi\leq r^*_\xi$ and also $q^\delta_{t^*\rest\xi}\rest\beta\leq
r^*_\xi\rest\beta$ (remember $(*)^a_{23}$ for $\vare\in
w_\delta\cap\xi$). Now by induction on $\alpha\in [\beta,\xi)$ we argue that
$q^\delta_{t^*\rest\xi}\rest\alpha\leq r^*_\xi\rest\alpha$. So suppose that 
$\beta\leq \alpha<\xi$ and we know already that $q^\delta_{t^*\rest\xi}
\rest\alpha\leq r^*_\xi\rest\alpha$. It follows from $(*)_3+(*)_5+(*)_6$
that $r^*_\xi\rest\alpha\forces_{\bbP_\alpha}\big(\forall i<\delta\big)\big(
r_i(\alpha)\leq p^\delta_{t^*\rest\xi}(\alpha)\leq q^\delta_{t^*\rest\xi}
(\alpha)\big)$ and therefore we may use $(*)_8$ to conclude that
\[r^*_\xi\rest\alpha\forces_{\bbP_\alpha} q^\delta_{t^*\rest\xi}(\alpha)
\leq r_\delta(\alpha)\leq q(\alpha)\leq r(\alpha)=r^*_\xi(\alpha),\]
as desired. Finally we define $r^{**}_\xi\in\bbP_\xi$ essentially by
$(*)^b_{23}$.  
\smallskip

Now suppose that $\otp(w_\delta\cap\xi)$ is a successor ordinal and let
$\xi_0=\max(w_\delta\cap\xi)$. Assume we know that $t^*\rest \xi_0\in 
T_\delta$ and that we have already defined $r^*_{\xi_0},r^{**}_{\xi_0} \in
\bbP_{\xi_0}$. It follows from the choice of $q$ and from $(*)^b_{20}$ that
\[r^{**}_{\xi_0}\forces_{\bbP_{\xi_0}}\mbox{`` }r(\xi_0)\forces_{
  \name{\bbQ}_{\xi_0}} \big(\exists \eta\in\name{\cX}_{\delta,\xi_0}\big)
\big(\name{\bar{q}}_{\delta,{\xi_0}}(\eta)\in\Gamma_{\name{\bbQ}_{\xi_0}}
\big)\mbox{ ''.}\]  
Thus we may choose $r^*\in\bbP_{\xi_0+1}$ and $\eta\in P_\delta$ such that
$r^{**}_{\xi_0}\leq r^*\rest\xi_0$, $r^*\rest\xi_0\forces r(\xi_0)\leq
r^*(\xi_0)$ and $r^*\rest\xi_0\forces_{\bbP_{\xi_0}}$`` $\eta\in
\name{\cX}_{\delta,\xi_0}\ \&\ \name{\bar{q}}_{\delta,\xi_0}(\eta)\leq
r^*(\xi_0)$ ''. Then $r\rest (\xi_0+1)\leq r^*$ and (by $(*)_7$,
\ref{supergame}(1)(c,d)) $r^*\rest\xi_0\forces_{\bbP_{\xi_0}}\big(\forall
\alpha<\delta\big)\big(\name{\tau}(\xi_0,\eta(\alpha))\leq r^*(\xi_0) \big)$ 
and  hence (by $(*)^a_{21}$) $r\rest\xi_0\forces \name{\tau}(\xi_0,
\eta(\alpha)) \leq r(\xi_0)$  for all $\alpha<\delta$. Therefore 
\begin{enumerate}
\item[$(*)^a_{24}$] for all $\alpha<\delta$, $\name{\tau}(\xi_0,
  \eta(\alpha)) \in \{r'(\xi_0):r'\in H^\delta\}$.
\end{enumerate}
Since $p^{\delta,\xi_0}_{t^*\rest\xi_0}\leq r^*\rest\xi_0$ (remember
$(*)^a_{23}$ for $\xi_0$ and $(*)_5$) we may use $(*)_{11}$ to conclude that
for some non-empty set $X\subseteq P_\delta$ we got
$p^{\delta,\xi_0}_{t^*\rest\xi_0}\forces_{\bbP_{\xi_0}}X=\name{\cX}_{\delta,\xi_0}$
and $\eta\in X$ satisfies (ii) of $(*)_{22}$. Hence $(t^*)_{\xi_0}\in X$ is
such that 
\begin{enumerate}
\item[$(*)^b_{24}$] for all $\alpha<\delta$, $\name{\tau}(\xi_0,
  (t^*)_{\xi_0}(\alpha)) \in \{r'(\xi_0):r'\in H^\delta\}$,
\end{enumerate}
and in particular $t^*\rest\xi\in T_\delta$ (remember $(*)_{11}$). We claim
that $(t^*)_{\xi_0}=\eta$. If not, then by \ref{supergame}(1)(b) we have 
\[r^*\rest\xi_0\forces_{\bbP_{\xi_0}}\big(\exists\alpha<\delta\big)\big( 
\name{\tau}(\xi_0,(t^*)_{\xi_0}(\alpha)),\name{\tau}(\xi_0,\eta(\alpha))
\mbox{ are incompatible in }\name{\bbQ}_{\xi_0}\big),\] 
so we may pick $\alpha<\delta$ and a condition $r^+\in \bbP_{\xi_0}$ such
that $r^*\rest\xi_0\leq r^+$ and
\[r^+\forces_{\bbP_{\xi_0}}\mbox{`` } \name{\tau}(\xi_0,(t^*)_{\xi_0}
(\alpha)), \name{\tau}(\xi_0,\eta(\alpha))\mbox{ are incompatible in
}\name{\bbQ}_{\xi_0}\mbox{ ''.}\] 
However, $r^+\forces_{\bbP_{\xi_0}}\mbox{`` }\name{\tau}(\xi_0,
(t^*)_{\xi_0}(\alpha)) \leq r(\xi_0)\ \&\  \name{\tau}(\xi_0,\eta(\alpha))
\leq r(\xi_0)\mbox{ ''}$ (by $(*)^a_{24}+(*)^b_{24}$), a contradiction. 

Now we define $r^*_\xi\in\bbP_\xi$ so that $r^*_\xi\rest (\xi_0+1)=r^*$ and
$r^*_\xi \rest (\xi_0,\xi)=r\rest (\xi_0,\xi)$. By the above considerations
and $(*)_{14}$ we know that $q^\delta_{t^*\rest\xi}\rest (\xi_0+1)\leq r^*=
r^*_\xi \rest (\xi_0+1)$. Exactly like in the case of limit $\otp( w_\delta
\cap \xi)$ we argue that $q^\delta_{t^*\rest\xi}\leq r^*_\xi$. Finally, we
choose $r^{**}_\xi\in \bbP_\xi$ by $(*)^b_{23}$. 
\smallskip

The last stage $\gamma$ of the inductive process described above shows that
$t^*\in T_\delta$ and $q^\delta_{t^*}\leq r^*_\gamma$, $r\leq
r^*_\gamma$. Now the claim readily follows. 
\end{proof}

We finish the proof of part (a) of the theorem exactly like in the proof of
\ref{lordsthm}. 
\medskip

\noindent(b)\quad Included in the proof of the first part.
\end{proof}

\section{Examples and counterexamples}
Let us note that our canonical test forcing $\fEE$ is B--noble: 

\begin{proposition}
\label{Qisnoble}
  Assume that $\bar{E},E$ are as in \ref{exBp} and $\bp=(\bar{P},S,D)$ is a
  $\Dl$--parameter on $\lambda$ such that $\lambda\setminus S\in E$. Then
  the forcing $\fEE$ is B--noble over $\bp$.
\end{proposition}

\begin{proof}
  The proof is a small modification of that of \ref{Qisgood}(2). First we
  fix an enumeration $\langle\nu_\alpha:\alpha<\lambda\rangle=
  {}^{<\lambda}\lambda$ (remember $\lambda^{<\lambda}=\lambda$), and for
  $\alpha<\lambda$ let $f(\alpha)\in\fEE$ be a condition such that
  $\mrot\big(f(\alpha)\big)=\nu_\alpha$ and 
\[\big(\forall\nu\in f(\alpha)\big)\big(\nu_\alpha\trianglelefteq\nu\
\Rightarrow\ \suc_{f(\alpha)}(\nu)=\lambda\big).\]
Let $p\in \fEE$. Consider the following strategy $\st$ of Generic in
$\supgame(p,\fEE)$. In the course of the play Generic is instructed to
build aside a sequence $\langle T_\xi:\xi<\lambda\rangle$ so that if
$\langle f_\xi,\cX_\xi,\bar{p}_\xi,\bar{q}_\xi:\xi<\lambda\rangle$ is the
sequence of the innings of the two players, then the following conditions
(a)--(d) are satisfied.
\begin{enumerate}
\item[(a)]  $T_\xi\in \fEE$ and if $\xi<\zeta<\lambda$ then $p=T_0 \supseteq
  T_\xi\supseteq T_\zeta$ and $T_\zeta\cap {}^\xi\lambda= T_\xi\cap {}^\xi
  \lambda$,  
\item[(b)] if $\xi<\lambda$ is limit, then $T_\xi=\bigcap\limits_{\zeta<\xi}
  T_\zeta$, 
\item[(c)] if $\xi\in S$, then 
\item[$\bullet$] $f_\xi=f\rest\xi$ and $\cX_\xi\subseteq P_\xi$ is a maximal set
    (possibly empty) such that 
    \begin{enumerate}
\item[$(\alpha)$] for each $\eta\in\cX_\xi$ the family $\big\{f\big(\eta
(\alpha)\big):\alpha<\xi\big\}\cup \{T_\xi\}$ has an upper bound in $\fEE$
and $\lh\big(\bigcup\{\nu_{\eta(\alpha)}:\alpha<\xi\}\big)=\xi$, 
\item[$(\beta)$] if $\eta_0,\eta_1\in\cX_\xi$ are distinct, then for some
  $\alpha<\xi$ the conditions $f\big(\eta_0(\alpha)\big)$,
  $f\big(\eta_1(\alpha)\big)$ are incompatible,
    \end{enumerate}
\item[$\bullet$] for $\eta\in \cX_\xi$ the condition $p^\xi_\eta$ is an upper bound to
  $\big\{f\big(\eta(\alpha)\big):\alpha<\xi\big\}\cup\{T_\xi\}$, 
\item[$\bullet$] $T_{\xi+1}=\bigcup\{q^\xi_\eta:\eta\in\cX_\xi\}\cup \bigcup\big\{
  (T_\xi)_\nu: \nu\in {}^\xi\lambda\cap T_\xi\ \mbox{ and } \nu\notin
  p^\xi_\eta\mbox{ for }\eta\in\cX_\xi\big\}$,
\item[(d)] if $\xi\notin S$, then $\cX_\xi=\emptyset$, $f_\xi=f\rest \xi$
  and $T_{\xi+1}=T_\xi$.  
\end{enumerate}
After the play is over, Generic puts $p^*=\bigcap\limits_{\xi<\lambda} T_\xi
\subseteq {}^{<\lambda}\lambda$. Almost exactly as in the proof of
\ref{Qisgood}(2), one checks that $p^*\in\fEE$ is a condition witnessing
$(\circledast)^{B+}_{\bp}$ of \ref{supergame}(1).
\end{proof}

\begin{definition}
\label{defpEE}
  Let $\bar{E},E$ be as in \ref{exBp}. We define a forcing notion $\pEE$ as
  follows.

\noindent {\bf A condition $p$ in $\pEE$} is a complete $\lambda$--tree
$p\subseteq {}^{<\lambda}\lambda$ such that 
\begin{itemize}
\item for every $\nu\in p$, either $|\suc_p(\nu)|=1$ or $\suc_p(\nu)\in
  E_\nu$, and
\item for some set $A\in E$ we have 
\[\big(\forall\nu\in p\big)\big(\lh(\nu)\in A\ \Rightarrow\ \suc_p(\nu)\in
E_\nu \big).\]
\end{itemize}

\noindent {\bf The order $\leq=\leq_{\pEE}$} is the reverse inclusion:
$p\leq q$ if and only if ($p,q\in\pEE$ and ) $q\subseteq p$.
\end{definition}

\begin{proposition}
\label{Pisnoble}
  Assume that $\bar{E},E$ are as in \ref{exBp} and $\bp=(\bar{P},S,D)$ is a
  $\Dl$--parameter on $\lambda$ such that $\lambda\setminus S\in E$. Then
  $\pEE$ is a $({<}\lambda)$-complete forcing notion of size $2^\lambda$
  which is B-noble over $\bp$.
\end{proposition}

\begin{proof}
  The arguments of \ref{Qisnoble} can be repeated here with almost no
  changes (a slight modification is needed for the justification that 
  $p^*\in\pEE$).  
\end{proof}

We may use the forcing $\pEE$ to substantially improve \cite[Corollary
5.1]{RoSh:860}. First, let us recall the following definition.

\begin{definition}
Let $\cF$ be a filter on $\lambda$ including all co-bounded subsets of
$\lambda$, $\emptyset\notin\cF$. 
\begin{enumerate}
\item We say that a family $F\subseteq{}^\lambda\lambda$ is {\em
    $\cF$--dominating\/}  whenever 
\[\big(\forall g\in {}^\lambda\lambda\big)\big(\exists f\in F\big)
\big(\{\alpha<\lambda: g(\alpha)<f(\alpha)\}\in\cF\big).\]
The $\cF$--dominating number $\gd_\cF$ is the minimal size of an
  $\cF$--dominating family in ${}^\lambda\lambda$.
\item We say that a family $F\subseteq{}^\lambda\lambda$ is {\em
    $\cF$--unbounded\/}  whenever 
\[\big(\forall g\in {}^\lambda\lambda\big)\big(\exists f\in F\big)
\big(\{\alpha<\lambda: g(\alpha)<f(\alpha)\}\in(\cF)^+\big).\]
The $\cF$--unbounded number $\gb_\cF$ is the minimal size of an
  $\cF$--unbounded family in ${}^\lambda\lambda$.
\item If $\cF$ is the filter of co-bounded subsets of $\lambda$, then the
  corresponding dominating/unbounded numbers are also denoted by
  $\gd_\lambda,\gb_\lambda$. If $\cF$ is the filter generated by club
  subsets of $\lambda$, then the corresponding numbers are called $\gd_{\rm  
    cl},\gb_{\rm cl}$.   
\end{enumerate}
\end{definition}

\begin{corollary}
\label{improve860}
  Assume $\lambda=\lambda^{<\lambda}$, $2^\lambda=\lambda^+$. Suppose that
  $\bp=(\bar{P},S,D)$ is a $\Dl$--parameter on $\lambda$, and $E$ is a
  normal filter on $\lambda$ such that $\lambda\setminus S\in E$. Then there
  is a $\lambda^{++}$--cc $\lambda$--proper forcing notion $\bbP$ such that
\[\forces_{\bbP}\mbox{`` }2^\lambda=\lambda^{++}=\gb_{E^\bbP}=\gd_{E^\bbP}=
\gd_\lambda\ \&\ \gb_\lambda=\gb_{D[\bbP]}=\gd_{D[\bbP]}=\lambda^+ \mbox{ ''.}\] 
\end{corollary}

\begin{proof}
For $\nu\in {}^{<\lambda}\lambda$ let $E_\nu$ be the filter generated by
clubs of $\lambda$ and let $\bar{E}=\langle E_\nu:\nu\in {}^{<\lambda}
\lambda\rangle$. Let $\bar{\bbQ}=\langle\bbP_\xi,
\name{\bbQ}_\xi:\xi<\lambda^{++}\rangle$ be a $\lambda$--support iteration
such that for every $\xi<\lambda^{++}$,  $\forces_{\bbP_\xi}\mbox{``
}\name{\bbQ}_\xi=\pEE$ ''. (Remember, we use the convention that in
$\bV^{\bbP_\xi}$ the normal filter generated by $E$ is also denoted by $E$
etc.) Let $\bbP=\bbP_{\lambda^{++}}=\lim(\bar{\bbQ})$.  

It follows from \ref{nobleiteration}(a)+\ref{Pisnoble} that $\bbP$ is 
$\lambda$--proper. Using \cite[Theorem 2.2]{RoSh:890} (see also Eisworth
\cite[\S 3]{Ei03}) we see that $\bbP$ satisfies the $\lambda^{++}$--cc,
$\forces_{\bbP}2^\lambda=\lambda^{++}$ and $\bbP$ is
$({<}\lambda)$--complete. Thus, the forcing with $\bbP$ does  
not collapse cardinals and it also follows from \ref{nobleiteration}(b) that     
\[\forces_{\bbP}\mbox{`` ${}^\lambda\lambda\cap\bV$ is $D[\bbP]$--dominating
  in ${}^\lambda\lambda$ ''.}\]  
It is also easy to check, that for each $\xi<\lambda^{++}$
\[\forces_{\bbP}\mbox{`` ${}^\lambda\lambda\cap\bV^{\bbP_\xi}$ is not  
  $E$--unbounded in ${}^\lambda\lambda$ ''}\]
and hence  we may easily conclude that $\forces_{\bbP}$`` $\gb_{E^\bbP}= 
2^\lambda$ ''. 
\end{proof}

\begin{definition}
\label{boundedPQ}
Assume that 
\begin{itemize}
\item $\lambda$ is weakly inaccessible, $\lambda^{<\lambda}=\lambda$,
\item $\bH:\lambda\longrightarrow\lambda$ is such that $|\alpha|^+\leq
  |\bH(\alpha)|$ for each $\alpha<\lambda$,
\item $F$ is a normal filter on $\lambda$, $\bar{F}=\langle F_\nu:\nu\in
  \bigcup\limits_{\alpha<\lambda}\prod\limits_{\xi<\alpha} \bH(\xi)\rangle$
  where $F_\nu$ is a $\big({<}|\alpha|^+\big)$--complete filter on
  $\bH(\alpha)$ whenever $\nu\in \prod\limits_{\xi<\alpha}\bH(\xi)$,
  $\alpha<\lambda$. 
\end{itemize}
We define forcing notions $\qFF$ and $\pFF$ as follows.
\begin{enumerate}
\item {\bf A condition $p$ in $\qFF$} is a complete $\lambda$--tree
  $p\subseteq \bigcup\limits_{\alpha<\lambda} \prod\limits_{\xi<\alpha}
  \bH(\xi)$ such that  
\begin{enumerate}
\item[(a)] for every $\nu\in p$, either $|\suc_p(\nu)|=1$ or $\suc_p(\nu)\in
  F_\nu$, and
\item[(b)] for every $\eta\in \lim_\lambda(p)$ the set $\{\alpha<\lambda:
  \suc_p(\eta\rest\alpha)\in F_{\eta\rest\alpha}\}$ belongs to $F$.
\end{enumerate}
\noindent {\bf The order of $\qFF$} is the reverse inclusion. 
\item {\bf A condition $p$ in $\pFF$} is a complete $\lambda$--tree
  $p\subseteq \bigcup\limits_{\alpha<\lambda} \prod\limits_{\xi<\alpha}
  \bH(\xi)$ satisfying (a) above and 
\begin{enumerate}
\item[(b)$^+$] for some set $A\in F$ we have 
\[\big(\forall\nu\in p\big)\big(\lh(\nu)\in A\ \Rightarrow\ \suc_p(\nu)\in
F_\nu \big).\]
\end{enumerate}
\noindent {\bf The order of $\pFF$} is the reverse inclusion. 
\end{enumerate}
\end{definition}

\begin{proposition}
\label{boundarenobel}
Assume that $\lambda,\bH,\bar{F},F$ are as in \ref{boundedPQ}. Let
$\bp=(\bar{P},S,D)$ be a $\Dl$--parameter such that $\lambda\setminus S\in
F$. Then both $\qFF$ and $\pFF$ are strategically $({<}\lambda)$--complete
forcing notions of size $2^\lambda$ which are also B--noble over $\bp$.
\end{proposition}

\begin{proof}
Like \ref{Qisgood}(2), \ref{Qisnoble}, \ref{Pisnoble}.
\end{proof}

The property of being B--noble seems to be a relative of {\em properness for
  $D$--semi diamonds\/} introduced in \cite{RoSh:655} and even more so of
{\em properness over $D$--diamonds\/} studied in Eisworth
\cite{Ei03}. However, technical differences make it difficult to see what
are possible dependencies between these notions (see Problem
\ref{prob1}). In this context, let us note that there are forcing notions
which are proper over semi diamonds, but are not B--noble over any
$\Dl$--parameter $\bp$. Let us consider, for example, a forcing notion
$\bbP^*$ defined as follows: 
\smallskip

\noindent{\bf a condition in $\bbP^*$}\quad is a function $p$ such that 
\begin{enumerate}
\item[(a)] $\dom(p)\subseteq\lambda^+$, $\rng(p)\subseteq\lambda^+$,
$|\dom(p)|<\lambda$, and 
\item[(b)] if $\alpha_1<\alpha_2$ are both from $\dom(p)$, then
$p(\alpha_1)<\alpha_2$; 
\end{enumerate}
\noindent{\bf the order $\leq$ of $\bbP^*$}\quad is the inclusion
$\subseteq$. 
\smallskip

\begin{proposition}
[See {\cite[Prop. 4.1, 4.2]{RoSh:655}}]
$\bbP^*$ is $({<}\lambda)$-complete forcing notion which is proper over all
semi diamonds.
\end{proposition}

\begin{proposition}
 $\bbP^*$ is not B--noble over any $\Dl$--parameter $\bp$ on $\lambda$. 
\end{proposition}

\begin{proof}
Let $q\in\bbP^*$ be such that that $\lambda\in\dom(q)$ and let $\name{W}_0$ 
be a $\bbP^*$--name such that $\forces_{\bbP^*}$`` $\name{W}_0=\bigcup
\Gamma_{\bbP^*}\rest \lambda$ ''. Clearly 
\[q\forces_{\bbP^*}\mbox{`` }\name{W}_0\mbox{ is a function with }
\dom(\name{W}_0)\subseteq\lambda\mbox{ and }\rng(\name{W}_0) \subseteq
  \lambda \mbox{ ''.}\]
Let $\name{W}$ be a $\bbP^*$--name for a member of ${}^\lambda\lambda$ such
that  
\[q\forces_{\bbP^*}\mbox{`` }\name{W}_0\subseteq\name{W}\mbox{ and }
\big(\forall \alpha\in\lambda\setminus\dom(\name{W}_0)\big) \big(
  \name{W}(\alpha)=\alpha\big)\mbox{ ''.}\]
Now suppose that $\bp=(\bar{P},S,D)$ is a $\Dl$--parameter and $\langle
A_\alpha:\alpha<\lambda\rangle$ is a sequence of subsets of $\lambda$ such
that $|A_\alpha|<\lambda$ for $\alpha<\lambda$.  The following claim implies
that $\bbP^*$ (above the condition $q$) is not B--noble over $\bp$ (remember
\ref{nobleiteration}(b)). 

\begin{claim}
\label{cl1}
$\forces_{\bbP^*}\mbox{`` }\big\{\alpha<\lambda:\name{W}(\alpha)\in
A_\alpha \big\}\notin D[\bbP^*]\mbox{ ''.}$
\end{claim}

\begin{proof}[Proof of the Claim]
Suppose that $p\geq q$ and $\name{B}_i,\name{f}_i$ (for $i<\lambda$) are
$\bbP^*$--names for members of $D\cap\bV$ and members of
${}^\lambda\lambda$, respectively, such that
\[p\forces_{\bbP^*}\mbox{`` }
\mathop{\triangle}\limits_{i<\lambda} \name{B}_i\cap
\mathop{\triangle}\limits_{i<\lambda} \set^\bp(\name{f}_i)\subseteq
\{\alpha<\lambda:\name{W}(\alpha)\in A_\alpha\}\mbox{ ''}.\]
Build inductively a sequence $\langle p_i,B_i,f_i:i<\lambda\rangle$ such
that for each $i<\lambda$: 
\begin{enumerate}
\item[(i)]  $p_i\in\bbP^*$, $p\leq p_0\leq p_j\leq p_i$ for $j<i$, 
\item[(ii)]  $B_i\in D\cap\bV$, $f_i\in {}^\lambda\lambda$ and 
\item[(iii)] $p_i\forces_{\bbP^*}$`` $\name{B}_i=B_i$ and $\name{f}_j \rest
  i=f_j\rest i$ for all $j\leq i$'', and 
\item[(iv)]  $i<\sup\big(\dom(p_i)\cap \lambda\big)$.
\end{enumerate}
Since $B=\mathop{\triangle}\limits_{i< \lambda} B_i\cap
\mathop{\triangle}\limits_{i< \lambda}\set^\bp(f_i)\in D$, we may pick a
limit ordinal $\delta\in B$ such that $\big(\forall i<\delta\big)\big(
\sup\big((\dom(p_i)\cup\rng(p_i))\cap\lambda\big)<\delta\big)$. Put
$\alpha=\sup(A_\delta)+888$ and $p^+=\bigcup\limits_{i<\delta}
p_i\cup\{(\delta,\alpha)\}$. Then $p^+\in\bbP^*$ 
is a condition stronger than all $p_i$ for $i<\delta$ and 
$p^+\forces_{\bbP^*}\mbox{`` } \delta\in \mathop{\triangle}\limits_{i<
  \lambda} \name{B}_i\cap\mathop{\triangle}\limits_{i< \lambda}
\set^\bp(\name{f}_i) \mbox{ and }\name{W}(\delta)=\alpha\notin
A_\delta$ '', a contradiction.
\end{proof}
A similar construction can be carried out above any condition $q$ such that
for some $\alpha\in\dom(q)$ we have $\cf(\alpha)=\lambda$ (the set of such
conditions is dense in $\bbP^*$).   
\end{proof}

\section{$\fEE$ vs $\pEE$ and Cohen $\lambda$--reals}
The forcing notions $\fEE$ and $\pEE$ (introduced in \ref{exBp} and
\ref{defpEE}, respectively) may appear to be {\em almost\/} the
same. However, at least under some reasonable assumptions on $\bar{E},E$
they do have different properties. 

Suppose that $\bV\subseteq\bV^*$ are transitive universes of ZFC (with the
same ordinals) such that ${}^{<\lambda}\lambda\cap \bV= {}^{<\lambda}
\lambda\cap\bV^*$. We say that a function $c\in {}^\lambda 2\cap\bV^*$ is a  
$\lambda$--Cohen over $\bV$ if for every open dense set $U\subseteq
{}^{<\lambda}2$ (where ${}^{<\lambda}2$ is equipped with the partial order
of the extension of sequences), $U\in\bV$, there is $\alpha<\lambda$ such
that $c\rest \alpha\in U$. 

\begin{proposition}
  \label{PaddsCohen}
Assume that
\begin{enumerate}
\item[(a)] $\lambda$ is a strongly inaccessible cardinal,
\item[(b)] $S$ is the set of all strong limit cardinals $\kappa<\lambda$ of
  countable cofinality,
\item[(c)]  $E$ is a normal filter on $\lambda$ such that $S\in E$,
\item[(d)] $\bar{E}=\langle E_\nu:\nu\in {}^{<\lambda}\lambda\rangle$ is a
  system of $({<}\lambda)$--complete non-principal filters on $\lambda$.
\end{enumerate}
Then the forcing notion $\pEE$ adds a $\lambda$--Cohen over $\bV$.
\end{proposition}

\begin{proof}
Let $\kappa\in S$. We will say that a tree $T\subseteq {}^{<\kappa}\kappa$
is {\em $\kappa$--interesting\/} if 
\begin{enumerate}
\item[$(\odot)_0$]  $|\lim_\kappa(T)|=2^\kappa$ and for some increasing
  cofinal in $\kappa$ sequence $\langle\delta_n:n<\omega\rangle
  \subseteq\kappa$ we have  
\[\big(\forall n<\omega\big)\big(\forall\nu\in T\big)\big(\lh(\nu)\leq
\delta_n\ \Rightarrow\ \sup\big(\rng(\nu)\big)<\delta_{n+1}\big).\]   
\end{enumerate}
Note that there are only $2^\kappa$ many $\kappa$--interesting trees (for
$\kappa\in S$). Therefore we may fix a well ordering of the family of
$\kappa$--interesting trees of length $2^\kappa$ and choose by induction a
function $f_\kappa:{}^\kappa\kappa\longrightarrow {}^{<\kappa}2\setminus
\{\langle\rangle\}$ such that   
\begin{enumerate}
\item[$(\odot)_1$]  if $T\subseteq {}^{<\kappa}\kappa$ is a
  $\kappa$--interesting tree, then 
\[\big(\forall\sigma\in {}^{<\kappa}2\setminus\{\langle\rangle\}\big) \big( 
\exists\eta\in{\lim}_\kappa(T)\big)\big(f_\kappa( \eta)=\sigma\big).\]   
\end{enumerate}
Let $\name{W}$ be a $\pEE$--name such that $\forces_{\pEE}\name{W} = \bigcup
\big\{\mrot(p):p\in\Gamma_{\pEE}\big\}$. Plainly, $\forces\name{W}\in
{}^\lambda\lambda$. Next, let $\name{C}$ be a $\pEE$--name such that
$\forces_{\pEE}\name{C}=\{\kappa\in S:\name{W}\rest\kappa\in {}^\kappa
\kappa\}$ and let $\name{\tau}$ be a $\pEE$--name such that
\[\begin{array}{ll}
\forces_{\pEE}&\mbox{`` $\name{\tau}$ is the concatenation of all elements
  of the sequence }\\
&\ \ \langle f_\kappa(\name{W}\rest\kappa):\kappa\in\name{C}\rangle, \mbox{
  i.e., }\name{\tau}=\langle\ldots \conc f_\kappa(\name{W}\rest\kappa)\conc
\ldots\rangle_{\kappa\in\name{C}}  \mbox{ ''.} 
\end{array}\]
Plainly, $\forces\name{\tau}\in {}^\lambda 2$ (remember, $\langle\rangle
\notin \rng(f_\kappa)$ for $\kappa\in S$). 

We are going to argue that 
\[\forces_{\pEE}\mbox{`` $\name{\tau}$ is $\lambda$-Cohen over $\bV$ ''.}\] 
To this end suppose that $U\subseteq {}^{<\lambda}2$ is an open dense set
and $p\in\pEE$. Let 
\[B=\big\{\kappa\in S:\big(\forall\eta\in p\cap {}^\kappa\lambda\big) \big(  
\suc_p(\eta)\in E_\eta\big)\big\}\]  
(so $B\in E$). By induction on $n<\omega$ choose $\delta_n,T_n$ so that
\begin{enumerate}
\item[$(\odot)_2$] $\delta_n\in B$, $\delta_n<\delta_{n+1}$, $T_n\subseteq
  {}^{\leq\delta_n}\lambda$ is a complete tree (thus every chain in $T_n$
  has a $\trianglelefteq$--bound in $T_n$), 
\item[$(\odot)_3$] $T_n\subseteq T_{n+1}\subseteq p$, $T_{n+1}\cap
  {}^{\delta_n}\lambda=T_n\cap {}^{\delta_n}\lambda$, $T_0\subseteq {}^{\leq
    \delta_0}\delta_0$, $|T_0\cap {}^{\delta_0}\delta_0|=1$,
\item[$(\odot)_4$] if $\nu\in T_n\cap {}^{<\delta_n}\lambda$, then
  $\suc_{T_n}(\nu)\neq\emptyset$ and 
\[2\leq |\suc_{T_n}(\nu)|\ \Rightarrow\ \lh(\nu)\in\{\delta_0,\ldots,
\delta_{n-1}\},\]
\item[$(\odot)_5$]  if $\nu\in T_{n+1}\cap {}^{\delta_n}\lambda$, then
  $|\suc_{T_{n+1}}(\nu)|=\delta_n$ and 
\[\sup\big(\rng(\nu)\big)<\delta_{n+1}<\min\big(\suc_{T_{n+1}}(\nu) \big).\] 
 \end{enumerate}
Next put $\kappa=\sup(\delta_n:n<\omega)$ and $T=\bigcup\limits_{n<\omega} 
T_n$. Clearly $\kappa\in S$ and $T\subseteq {}^{<\kappa}\kappa$ is a
$\kappa$--interesting tree such that 
\[\big(\forall\eta\in{\lim}_\kappa(T)\big)\big(\kappa=\min\{\delta\in S:
\delta_0<\delta\ \&\ \eta\rest\delta\in {}^\delta\delta\}\big).\] 
Let $T_0\cap {}^{\delta_0}\delta_0=\{\eta_0\}$ and $C(\eta_0)=\{\delta\in 
S\cap (\delta_0+1): \eta_0\rest\delta\in {}^\delta\delta\}$, and let
$\tau_0$ be the concatenation of all elements of the sequence $\langle 
f_\delta(\eta_0\rest\delta): \delta\in C(\eta_0)\rangle$, i.e.,
$\tau_0=\langle \ldots\conc f_\delta(\eta_0\rest\delta)\conc\ldots 
\rangle_{\delta\in C(\eta_0)}\in {}^{<\lambda}2$. Pick $\sigma\in
{}^{<\lambda}2$ such that $\tau_0\vtl\sigma\in U$. It follows from
$(\odot)_1$ that we may find $\eta\in\lim_\kappa(T)$ such that
$\sigma=\tau_0\conc f_\kappa(\eta)$. Now note that
$(p)_\eta\forces_{\pEE}\sigma\vtl \name{\tau}$.  
\end{proof}

\begin{proposition}
  \label{QnoCohen}
Assume that
\begin{enumerate}
\item[(a)] $\lambda$ is a measurable cardinal,
\item[(b)] $\bar{E}=\langle E_\nu:\nu\in {}^{<\lambda}\lambda\rangle$ is a
  system of normal ultrafilters on $\lambda$, 
\item[(c)]  $E$ is a normal filter on $\lambda$,
\item[(d)] $p\in \fEE$ and $\name{\tau}$ is a $\fEE$--name such that
  $p\forces\name{\tau}\in {}^\lambda 2$,
\item[(e)] $\bar{\delta}=\langle\delta_\alpha:\alpha<\lambda\rangle$ is an
  increasing continuous sequence of non-successor ordinals below $\lambda$
  such that $\delta_0=0$ and $2^{2^{|\delta_\alpha|}}<|\delta_{\alpha+1}|$
  for all $\alpha<\lambda$.
\end{enumerate}
Then there are a condition $q\in\fEE$ and a sequence $\langle A_\alpha:
\alpha<\lambda\rangle$ such that 
\begin{enumerate}
\item[(i)]  $q\geq p$ and $A_\alpha\subseteq
  {}^{[\delta_\alpha,\delta_{\alpha+1})} 2$,
  $|A_\alpha|<|\delta_{\alpha+1}|$ for $\alpha<\lambda$, and 
\item[(ii)] $q\forces_{\fEE}\big(\forall\alpha<\lambda\big)\big( \name{\tau}
  \rest [\delta_\alpha,\delta_{\alpha+1})\in A_\alpha\big)$. 
\end{enumerate}
In particular, the forcing notion $\fEE$ does not add any $\lambda$--Cohen
over $\bV$. 
\end{proposition}

\begin{proof}
 Let $\langle\nu_\alpha:\alpha<\lambda\rangle$ be an enumeration of
 ${}^{<\lambda}\lambda$ such that $\nu_\alpha\vtl\nu_\beta$ implies
 $\alpha<\beta$. By induction on $\alpha<\lambda$ we will construct a
 sequence $\langle A_\alpha,p_\alpha,X_\alpha: \alpha<\lambda\rangle$ so
 that for each $\alpha<\lambda$ we have:
\begin{enumerate}
\item[$(\boxtimes)_1$] $A_\alpha\subseteq
  {}^{[\delta_\alpha,\delta_{\alpha+1})} 2$,
  $|A_\alpha|<|\delta_{\alpha+1}|$, $p_\alpha\in\fEE$, $X_\alpha\subseteq
  p_\alpha$, $|X_\alpha|\leq\delta_\alpha$, 
\item[$(\boxtimes)_2$] if $\alpha<\beta<\lambda$, then $p_\alpha\leq
  p_\beta$ and $X_\alpha\subseteq X_\beta$, 
\item[$(\boxtimes)_3$] $p_0=p$, $X_0=\{\mrot(p_0)\}$, and if $\alpha$ is
  limit then $p_\alpha=\bigcap\limits_{\beta<\alpha} p_\beta$ and $X_\alpha=
  \bigcup\limits_{\beta<\alpha} X_\beta$, 
\item[$(\boxtimes)_4$] if $\nu\in X_{\alpha+1}$, then $\suc_{p_\alpha}(\nu)
  \in E_\nu$, 
\item[$(\boxtimes)_5$] if $\alpha$ is limit, $\nu\in X_\alpha$ and
  $\alpha\in \bigcap\limits_{\beta<\alpha}\suc_{p_\beta}(\nu)$, then for
  some $\eta\in X_{\alpha+1}$ we have $\nu\conc \langle\alpha\rangle
  \trianglelefteq\eta$, 
\item[$(\boxtimes)_6$] if $\nu_\alpha\in p_\alpha$, then there is $\eta\in
  X_{\alpha+1}$ such that $\nu_\alpha\trianglelefteq\eta$ and if
  (additionally) $\suc_{p_\alpha}(\nu_\alpha)\in E_{\nu_\alpha}$, then
  $\eta=\nu_\alpha$, 
\item[$(\boxtimes)_7$] if $\alpha$ is limit and $\nu\in {}^\alpha\alpha \cap
  p_\alpha$ and $\suc_{p_\alpha}(\nu)\in E_\nu$, then $\nu\in X_{\alpha+1}$,
\item[$(\boxtimes)_8$] $p_{\alpha+1}\forces\name{\tau}\rest [\delta_\alpha,
  \delta_{\alpha+1})\in A_\alpha$. 
\end{enumerate}
Suppose that we have determined $p_\beta,X_\beta$ for $\beta<\alpha$ and
$A_\beta$ for $\beta+1<\alpha$ so that the relevant instances of
$(\boxtimes)_1$--$(\boxtimes)_8$ are satisfied. If $\alpha$ is limit or 0,
then $p_\alpha,X_\alpha$ are defined by $(\boxtimes)_3$ (and $A_\alpha$ will
be chosen at the next step). One easily verifies that $p_\alpha,X_\alpha$
satisfy the requirements in $(\boxtimes)_1$--$(\boxtimes)_4$. 

So suppose now that $\alpha=\gamma+1$ (and we have defined
$p_\gamma,X_\gamma$ and $A_\beta$ for $\beta<\gamma$). We may easily choose
a set $X_\alpha\subseteq p_\gamma$ such that 
\begin{enumerate}
\item[$(\boxtimes)_9$]   $X_\gamma\subseteq X_\alpha$,
  $|X_\alpha|<\delta_\alpha$ and $X_\alpha$ satisfies
  $(\boxtimes)_4$--$(\boxtimes)_7$ (with $\alpha$ there corresponding to
  $\gamma$ here), and 
\item[$(\boxtimes)_{10}$] if $\eta_0,\eta_1\in X_\alpha$, $\nu=\eta_0\cap
  \eta_1$, then $\nu\in X_\alpha$, and 
\item[$(\boxtimes)_{11}$] if $\langle\eta_\xi:\xi<\zeta\rangle\subseteq
  X_\alpha$ is $\vtl$--increasing, then there is $\eta\in X_\alpha$ such
  that $(\forall\xi<\zeta)(\eta_\xi\trianglelefteq\eta)$.
\end{enumerate}
Next, for each $\eta\in X_\alpha$ choose a function
$\sigma_\eta:[\delta_\gamma,\delta_{\gamma+1})\longrightarrow 2$ and a
condition $q_\eta\in\fEE$ so that 
\begin{enumerate}
\item[$(\boxtimes)_{12}$] $\mrot(q_\eta)=\eta$, $(p_\gamma)_\eta\leq
  q_\eta$, $\big(\forall \nu\in X_\alpha\big)\big(\eta\vtl\nu\ \Rightarrow\  
  \nu(\lh(\eta))\notin\suc_{q_\eta}(\eta)\big)$ and $q_\eta\forces
  \name{\tau}\rest [\delta_\gamma, \delta_{\gamma+1})=\sigma_\eta$.   
\end{enumerate}
(Possible by assumption (b) and \ref{expur}.) Put $p_\alpha=
\bigcup\limits_{\eta\in X_\alpha} q_\eta$ and $A_\gamma=\{\sigma_\eta:
\eta\in X_\alpha\}$. Plainly, $|A_\gamma|\leq |X_\alpha|<\delta_{\gamma+1}$
and $p_\alpha\in \fEE$ (to verify that $p_\alpha$ is a complete
$\lambda$-tree use $(\boxtimes)_{10}+(\boxtimes)_{11}$; the other
requirements easily follow from the fact that $|X_\alpha|<\lambda$). One
also easily checks that $p_\alpha\forces\name{\tau}\rest [\delta_\gamma, 
\delta_{\gamma+1}) \in A_\gamma$. 
\smallskip

After the inductive construction is carried out, we put
$q=\bigcap\limits_{\alpha<\lambda} p_\alpha$. It follows from
$(\boxtimes)_2+(\boxtimes)_5+(\boxtimes)_6$ that $q$ is a complete
$\lambda$--tree, $|\suc_q(\nu)|=1$ or $\suc_q(\nu)\in E_\nu$ for each
$\nu\in q$, and 
\[\bigcup\limits_{\alpha<\lambda}X_\alpha=\{\nu\in q:\suc_q(\nu)\in
E_\nu\}.\] 
Suppose now that $\eta\in\lim_\lambda(q)$. Then for each $\alpha<\lambda$ we
have $\eta\in\lim_\lambda(p_\alpha)$ and hence $B_\alpha \stackrel{\rm
  def}{=} \{\xi<\lambda:\suc_{p_\alpha}(\eta\rest\xi)\in E_{\eta\rest\xi}\}
\in E$. Let 
\[C=\{\delta<\lambda:\delta\mbox{ is limit and }\eta\rest\delta\in {}^\delta
\delta\}\] 
(it is a club of $\lambda$). Since $E$ is a normal filter, $C\cap
\mathop{\triangle}\limits_{\alpha< \lambda} B_\alpha\in E$. Suppose
$\delta\in C\cap\mathop{\triangle}\limits_{\alpha< \lambda} B_\alpha$. Then
by $(\boxtimes)_3+(\boxtimes)_7$ we have $\eta\rest\delta\in X_{\delta+1}$
and thus $\suc_q(\eta\rest\delta)\in E_{\eta\rest\delta}$. Consequently,
$q\in\fEE$.

Finally, it follows from $(\boxtimes)_8$ that $q\forces \big(\forall\alpha<
\lambda\big) \big(\name{\tau}\rest [\delta_\alpha,\delta_{\alpha+1})\in
A_\alpha \big)$. 
\end{proof}

Let us note that forcing notions of the form $\fEE$ may add
$\lambda$--Cohens if the filters $E_\nu$ are far from being ultrafilters.

\begin{proposition}
  Assume that
\begin{enumerate}
\item[(a)] $E$ is a normal filter on $\lambda=\lambda^{<\lambda}$,
\item[(b)] $\bar{E}=\langle E_\nu:\nu\in {}^{<\lambda}\lambda\rangle$ is a
  system of $({<}\lambda)$--complete filters on $\lambda$, 
\item[(c)]  for every $\nu\in {}^{<\lambda}\lambda$ there is a family
  $\{A^\nu_\alpha:\alpha<\lambda\}$ of pairwise disjoint sets from
  $(E_\nu)^+$. 
\end{enumerate}
Then both the forcing notions $\fEE$ and $\pEE$ add $\lambda$--Cohens over
$\bV$. 
\end{proposition}

\begin{proof}
We will sketch the argument for $\fEE$ only (no changes are needed for the
case of $\pEE$).  

For each $\nu\in {}^{<\lambda}\lambda$ choose a function $h_\nu:\lambda
  \longrightarrow {}^{<\lambda}2\setminus\{\langle\rangle\}$ such that 
\[\big(\forall\sigma\in {}^{<\lambda}2\big)\big(\exists \alpha<\lambda
\big)\big(\forall\xi\in A^\nu_\alpha\big)\big(h_\nu(\xi)=\sigma\big).\]
Let $\name{W}$ be a $\fEE$--name such that $\forces_{\fEE}\name{W} = \bigcup  
\big\{\mrot(p):p\in\Gamma_{\fEE}\big\}$ and let $\name{\tau}$ be a
$\fEE$--name such that 
\[\begin{array}{ll}
\forces_{\fEE}&\mbox{`` $\name{\tau}$ is the concatenation of all elements
  of the sequence }\\
&\ \ \big\langle h_{\name{W}\rest\xi}\big(\name{W}(\xi)\big):\xi<\lambda
\big\rangle, \mbox{  i.e., }\\
&\ \ \name{\tau}=\big\langle h_{\langle\rangle}\big(\name{W}(0)\big) \conc
h_{\langle\name{W}(0)\rangle}\big(\name{W}(1)\big)\conc \ldots \conc
h_{\name{W}\rest\xi} \big(\name{W}(\xi)\big)\conc
\ldots\rangle_{\xi<\lambda} \mbox{ ''.}   
\end{array}\]
One easily verifies that $\forces$`` $\name{\tau}\in {}^\lambda 2$ is a
$\lambda$--Cohen over $\bV$ ''.
\end{proof}

The result in \ref{QnoCohen} would be specially interesting if we only knew
that it is preserved in $\lambda$--support iterations. Unfortunately, at the
moment we do not know if this is true (see Problem \ref{prob2}(1)). However,
we may consider properties stronger than adding $\lambda$--Cohens and then
our earlier results give some input.

\begin{definition}
\label{strongdef}
Suppose that $\bV\subseteq\bV^*$ are transitive universes of ZFC (with the
same ordinals) such that ${}^{<\lambda}\lambda\cap \bV= {}^{<\lambda}
\lambda\cap\bV^*$. We say that a function $c\in {}^\lambda 2\cap\bV^*$ is a  
\begin{enumerate}
\item {\em strongly$^\ominus$ $\lambda$--Cohen over $\bV$\/} if it is a 
  $\lambda$--Cohen (i.e.,  for every open dense set $U\subseteq
{}^{<\lambda}2$ from $\bV$ there is $\alpha<\lambda$ such that $c\rest
\alpha\in U$) and
\begin{enumerate}
\item[$(\ominus)$]  if $\langle\eta_\alpha,\beta_\alpha:\alpha<\lambda
  \rangle\in\bV$ is such that $\alpha<\beta_\alpha<\lambda$ and
  $\eta_\alpha \in {}^{[\alpha,\beta_\alpha)}2 $ for $\alpha<\lambda$, then 
\[\bV^*\models\{\alpha<\lambda:\eta_\alpha \nsubseteq c\}\mbox{ is
  stationary};\]
\end{enumerate}
\item {\em strongly$^\oplus$ $\lambda$--Cohen over $\bV$\/} if it is a 
  $\lambda$--Cohen and
\begin{enumerate}
\item[$(\oplus)$]  if $\langle\eta_\alpha,\beta_\alpha:\alpha<\lambda
  \rangle\in\bV$ is such that $\alpha<\beta_\alpha<\lambda$ and
  $\eta_\alpha \in {}^{[\alpha,\beta_\alpha)} 2$ for $\alpha<\lambda$, then  
\[\bV^*\models\{\alpha<\lambda:\eta_\alpha\subseteq c\}\mbox{ is
  stationary}.\]
\end{enumerate}
\item More generally, if $D$ is a normal filter on $\lambda$, $D\in \bV$ then 
  we say that $c\in {}^\lambda2\cap \bV^*$ is {\em $D$--strongly$^\oplus$
    $\lambda$--Cohen over $\bV$\/} if in $(\oplus)$ we replace
  ``stationary'' by ``$\in (D^{\bV^*})^+$'' (where $D^{\bV^*}$ is the normal
  filter generated by $D$ in $\bV^*$). Similarly for {\em
    strongly$^\ominus$}. 
\end{enumerate}
\end{definition}

\begin{remark}
\begin{enumerate}
\item To explain our motivation for \ref{strongdef}, let us recall that  
if $c\in {}^\lambda 2$ is  $\lambda$--Cohen over $\bV$ and $\langle
\eta_\alpha,\beta_\alpha:\alpha<\lambda\rangle\in\bV$ is such that
$\alpha<\beta_\alpha<\lambda$ and $\eta_\alpha \in
{}^{[\alpha,\beta_\alpha)} 2$ for $\alpha<\lambda$, then
\[\bV^*\models\mbox{`` both }\{\alpha<\lambda:\eta_\alpha\subseteq c\}
\mbox{ and } \{\alpha<\lambda:\eta_\alpha\nsubseteq c\}\mbox{ are unbounded
  in $\lambda$ ''}.\] 
\item Let $\langle\eta_\alpha,\beta_\alpha:\alpha<\lambda\rangle \in\bV$ be 
  such that $\alpha<\beta_\alpha<\lambda$ and $\eta_\alpha \in
  {}^{[\alpha,\beta_\alpha)} 2$ for $\alpha<\lambda$. Let
  $\bbC=({}^{<\lambda}2,\vtl)$ (so this is the $\lambda$--Cohen forcing
  notion) and let $\name{c}$ be the canonical $\bbC$--name for  the generic
  $\lambda$--real (i.e., $\forces_{\bbC}\name{c}=\bigcup\Gamma_{\bbC}$). Let
  $\name{\bbQ}$ be a $\bbC$--name for a forcing notion in which conditions
  are closed bounded sets $d\subseteq\lambda$ such that $(\forall\alpha\in
  d) (\eta_\alpha\subseteq \name{c})$ ordered by end extension. 
Then $\bbC*\name{\bbQ}$ is essentially the $\lambda$-Cohen forcing and 
\[\forces_{\bbC*\name{\bbQ}}\mbox{`` $\name{c}$ is not strongly$^\ominus$
  $\lambda$--Cohen over $\bV$ ''.}\]
Hence, if we add a $\lambda$--Cohen then we also add a non-strong$^\ominus$ 
$\lambda$--Cohen.
\item Note that {\em strongly$^\oplus$ $\lambda$--Cohen\/} implies
  {\em strongly$^\ominus$ $\lambda$--Cohen\/}. (Simply, for a sequence
  $\langle \eta_\alpha,\beta_\alpha:\alpha<\lambda\rangle$ consider $\langle
  1-\eta_\alpha, \beta_\alpha:\alpha<\lambda\rangle$.)
 \end{enumerate}
\end{remark}

\begin{proposition}
Assume that $\lambda$ is a strongly inaccessible cardinal and $\bp=
(\bar{P},\lambda,D_\lambda)$ is a $\Dl$--parameter such that
$P_\delta={}^\delta\delta$ and $D_\lambda$ is the filter generated by club
subsets of $\lambda$. 
\begin{enumerate}
\item If a forcing notion $\bbQ$ is reasonably B--bounding
over $\bp$, then 
\[\forces_{\bbQ}\mbox{`` there is no strongly$^\oplus$ $\lambda$--Cohen
over $\bV$ ''}.\]
\item If $\bar{\bbQ}=\langle\bbP_\alpha,\name{\bbQ}_\alpha:\alpha<\gamma 
  \rangle$ is a $\lambda$--support iteration such that for every
  $\alpha<\lambda$, 
\[\forces_{\bbP_\alpha}\mbox{`` }\name{\bbQ}_\alpha\mbox{ is reasonably
  B--bounding over $\bp[\bbP_\alpha]$ '',}\]   
then 
\[\forces_{\bbP_\gamma}\mbox{`` there is no strongly$^\oplus$
  $\lambda$--Cohen over $\bV$ ''}.\]
\end{enumerate}
\end{proposition}

\begin{proof}
(1)\quad Note that, in $\bV^\bbQ$,  $D_\lambda[\bbQ]$ is the filter
generated by clubs of $\lambda$. 

Let $p \in\bbQ$ and let $\name{\eta}$ be a $\bbQ$--name such that
$p\forces\name{\eta}\in {}^\lambda 2$. Let $\st$ be a winning strategy of 
Generic in the game $\supgame(p,\bbQ)$. 

Let us consider a play of $\supgame(p,\bbQ)$ in which Generic follows the
instructions of $\st$ while Antigeneric plays as follows. In the course of
the play, in addition to his innings $q^\alpha_t$, Antigeneric constructs
aside a sequence $\big\langle \kappa_\alpha,\langle\eta^\alpha_t:t\in
I_\alpha\rangle:\alpha<\lambda\big\rangle$ such that if  $\big\langle
I_\alpha, \langle p^\alpha_t,q^\alpha_t:t\in I_\alpha\rangle: \alpha<\lambda 
\big\rangle$ is the sequence of the innings of the two players then the
following two demands are satisfied.
\begin{enumerate}
\item[$(\boxdot)_1$] $\kappa_\alpha$ is a cardinal such that
  $2^{|I_\alpha|+\alpha+\aleph_0}<\kappa_\alpha$ and $\eta^\alpha_t\in
  {}^{\kappa_\alpha}2$ (for $t\in I_\alpha$),
\item[$(\boxdot)_2$] $q^\alpha_t\forces_{\bbQ}\name{\eta}\rest\kappa_\alpha
  =\eta^\alpha_t$ for each $t\in I_\alpha$. 
\end{enumerate}
Since the play is won by Generic, there is a condition $q\geq p$ such that
\[q\forces_{\bbQ}\mbox{`` }\{\alpha<\lambda: (\exists t\in
I_\alpha)(q^\alpha_t\in \Gamma_\bbQ)\}\mbox{ contains a club of $\lambda$
  ''.}\] 
It follows from $(\boxdot)_1$ that for each $\alpha<\lambda$ we may choose
$\vare^0_\alpha<\vare^1_\alpha$ from the interval $(\alpha,\kappa_\alpha)$
such that $\big(\forall t\in I_\alpha\big)\big(\eta^\alpha_t(\vare^0_\alpha
)=\eta^\alpha_t(\vare^1_\alpha)\big)$. For each $\alpha<\lambda$ choose
$\nu_\alpha:[\alpha,\kappa_\alpha)\longrightarrow 2$ so that
$\nu_\alpha(\vare^0_\alpha)=0$ and $\nu_\alpha(\vare^1_\alpha)=1$. Then 
\[q\forces_{\bbQ}\mbox{`` }\{\alpha<\lambda:\nu_\alpha\nsubseteq \name{\eta}
\}\mbox{ contains a club of $\lambda$ ''.}\] 
\smallskip

\noindent (2)\quad Similar, but we have to work with trees of conditions as
in the proof of \ref{like860}.
\end{proof}

\section{Marrying {\em B--bounding\/} with {\em fuzzy proper}}
In this section we introduce a property of forcing notions which, in a
sense, marries the B--bounding forcing notions of \cite[Definition
3.1(5)]{RoSh:860} with the fuzzy proper forcings introduced in \cite[\S
A.3]{RoSh:777}. This property, defined in the language of games, is based on
two games: {\em the servant game\/} $\sergame$ which is the part coming from
the fuzzy properness and {\em the master game\/} $\masgame$ which is
related to the reasonable boundedness property. Later in this section we
will even formulate a true preservation theorem for a slightly modified game. 

In this section we assume the following:

\begin{context}
\label{consec6}
\begin{enumerate}
\item $\lambda$ is a strongly inaccessible cardinal,
\item $D$ is a normal filter on $\lambda$,
\item $S\in D$, $0\notin S$, all successor ordinals below $\lambda$ belong
  to $S$, $\lambda\setminus S$ is unbounded.
\end{enumerate}
\end{context}

\begin{definition}
  \label{servant}
Let $\bbQ$ be a forcing notion.
\begin{enumerate}
\item {\em A $\bbQ$--servant over $S$\/} is a sequence $\bar{q}=\langle
  q^\delta_t:\delta\in S\ \&\ t\in I_\delta\rangle$ such that
  $|I_\delta|<\lambda$ (for $\delta\in S$) and $q^\delta_t\in\bbQ$ (for
  $\delta\in S$, $t\in I_\delta$).
\item Let $\bar{q}$ be a $\bbQ$--servant over $S$ and $q\in\bbQ$. We define
  a game $\sergame(\bar{q},q,\bbQ)$ as follows. A play of
  $\sergame(\bar{q},q,\bbQ)$ lasts at most $\lambda$ steps during which the
  players, COM and INC, attempt to construct a sequence $\langle r_\alpha,
  A_\alpha:\alpha<\lambda\rangle$ such that 
  \begin{itemize}
\item $r_\alpha\in\bbQ$, $q\leq r_\alpha$, $A_\alpha\in D$ and $\alpha<\beta
  <\lambda\ \Rightarrow\ r_\alpha\leq r_\beta$.
  \end{itemize}
The terms $r_\alpha,A_\alpha$ are chosen successively by the two players so
that 
\begin{itemize}
\item if $\alpha\notin S$, then INC picks $r_\alpha,A_\alpha$, and 
\item if $\alpha\in S$, then COM chooses $r_\alpha,A_\alpha$.
\end{itemize}
If at some moment of the play one of the players has no legal move, then INC
wins; otherwise, if both players always had legal moves and the sequence
$\langle r_\alpha,A_\alpha:\alpha<\lambda\rangle$ has been constructed, then
COM wins if and only if 
\begin{enumerate}
\item[$(\odot)$] \quad $\big(\forall\delta\in S\big)\big([\delta\in
  \bigcap\limits_{\alpha<\delta} A_\alpha\ \&\ \delta\mbox{ is limit }]\
  \Rightarrow\ (\exists t\in I_\delta)(q^\delta_t\leq r_\delta)\big)$. 
\end{enumerate}
\item If COM has a winning strategy in the game $\sergame(\bar{q},q,\bbQ)$,
  then we will say that {\em $q$ is an $(S,D)$--knighting condition for the
    servant $\bar{q}$}. 
\end{enumerate}
\end{definition}

\begin{definition}
\label{master}
Let $\bbQ$ be a strategically $({<}\lambda)$--complete forcing notion.
\begin{enumerate}
\item For a condition $p\in\bbQ$ we define a game $\masgame(q,\bbQ)$ between
  two players, Generic and Antigeneric. The game is a small modification of
  $\bgame(p,\bbQ)$ (see \ref{Bplus}) --- the main difference is in the
  winning condition. A play of $\masgame(p,\bbQ)$ lasts $\lambda$ steps and
  during a play a sequence    
\[\Big\langle I_\alpha,\langle p^\alpha_t,q^\alpha_t:t\in I_\alpha\rangle:
\alpha<\lambda\Big\rangle\]   
is constructed. Suppose that the players have arrived to a stage
$\alpha<\lambda$ of the game. Now,
\begin{enumerate}
\item[$(\aleph)_\alpha$]  first Generic chooses a non-empty set $I_\alpha$
  of cardinality $<\lambda$ and a system $\langle p^\alpha_t:t\in I_\alpha
  \rangle$ of conditions from $\bbQ$, 
\item[$(\beth)_\alpha$]  then Antigeneric answers by picking a system
  $\langle q^\alpha_t:t\in I_\alpha\rangle$ of conditions from $\bbQ$ such
  that $(\forall t\in I_\alpha)(p^\alpha_t\leq q^\alpha_t)$.
\end{enumerate}
At the end, Generic wins the play $\big\langle I_\alpha, \langle p^\alpha_t,
q^\alpha_t:t\in I_\alpha\rangle:\alpha<\lambda\big\rangle$ of
$\masgame(p,\bbQ)$ if and only if letting $\bar{q}=\langle q^\alpha_t:
\alpha\in S\ \&\ t\in I_\alpha\rangle$ (it is a $\bbQ$--servant over $S$) we
have  
\begin{enumerate}
\item[$(\circledast)^{D,S}_{\rm master}$] there exists an $(S,D)$--knighting
  condition $q\geq p$ for the servant $\bar{q}$. 
\end{enumerate}
\item A forcing notion $\bbQ$ is {\em reasonably merry over $(S,D)$\/} if
  (it is strategically $({<}\lambda)$--complete and) Generic has a winning
  strategy in the game $\masgame(p,\bbQ)$ for any $p\in\bbQ$.
\end{enumerate}
\end{definition}

\begin{theorem}
 \label{masiter}
Assume that $\lambda,S,D$ are as in \ref{consec6}. Let $\bar{\bbQ}=\langle
\bbP_\alpha,\name{\bbQ}_\alpha:\alpha<\gamma\rangle$ be a $\lambda$--support
iteration such that for each $\alpha<\gamma$:
\[\forces_{\bbP_\alpha}\mbox{`` }\name{\bbQ}_\alpha\mbox{ is reasonably
  merry over $(S,D)$''.}\]
Then 
\begin{enumerate}
\item[(a)] $\bbP_\gamma=\lim(\bar{\bbQ})$ is $\lambda$--proper, and 
\item[(b)] for every $\bbP_\gamma$--name $\name{\tau}$ for a function
  from $\lambda$ to $\bV$ and a condition $p\in\bbP_\gamma$, there are
  $q\geq p$ and  $\langle A_\xi:\xi<\lambda\rangle$ such that
  $(\forall\xi<\lambda)(|A_\xi|<\lambda)$ and $q\forces\mbox{``
  }\{\xi<\lambda:\name{\tau}(\xi) \in A_\xi\}\in \left(D^{\bbP_\gamma}
  \right)^+\mbox{ ''}$.  
\end{enumerate}
\end{theorem}

\begin{proof}
(a)\quad  The proof starts with arguments very much like those in
\cite[Theorems 3.1, 3.2]{RoSh:860}, so we will state only what should be
done (without actually describing how the construction can be carried
out). The major difference comes later, in arguments that the chosen
condition is suitably generic.

Suppose that $N\prec (\cH(\chi),\in,<^*_\chi)$ is such that 
\[{}^{<\lambda}N\subseteq N,\quad |N|=\lambda\ \mbox{ and }\
\bar{\bbQ},S,D,\ldots\in N.\] 
Let $p\in N\cap \bbP_\gamma$ and $\langle\name{\tau}_\alpha: \alpha<
\lambda\rangle$ list all $\bbP_\gamma$--names for ordinals from $N$. For
each $\xi\in N\cap \gamma$ fix a $\bbP_\xi$--name $\name{\st}^0_\xi\in N$
for a winning strategy of Complete in $\Game^\lambda_0(\name{\bbQ}_\xi,
\name{\emptyset}_{\name{\bbQ}_\xi})$ such that it instructs Complete to play
$\name{\emptyset}_{\name{\bbQ}_\xi}$ as long as her opponent plays
$\name{\emptyset}_{\name{\bbQ}_\xi}$.

Let us pick an increasing continuous sequence $\langle w_\delta:\delta<
\lambda\rangle$ of subsets of $\gamma$ such that
$\bigcup\limits_{\delta<\lambda} w_\delta=N\cap\gamma$, $w_0=\{0\}$ and
$|w_\delta|<\lambda$. 

By induction on $\delta<\lambda$ choose
\begin{enumerate}
\item[$(\boxtimes)_\delta$] \qquad $\cT_\delta,\bar{p}^\delta,
  \bar{q}^\delta, r^-_\delta, r_\delta$, $\langle\name{\vare}_{\delta,\xi},
  \name{\bar{p}}_{\delta,\xi},\name{\bar{q}}_{\delta,\xi}:\xi\in
  w_\delta\rangle$, and $\name{\st}_\xi$ for $\xi\in w_{\delta+1}\setminus w_\delta$  
\end{enumerate}
so that if the following conditions $(*)_0$--$(*)_{11}$ are satisfied (for each 
$\delta<\lambda$).  
\begin{enumerate}
\item[$(*)_0$] All objects listed in $(\boxtimes)_\delta$ belong to $N$ and
  they are known after stage $\delta$ of the construction. 
\item[$(*)_1$] $r^-_\delta,r_\delta\in \bbP_\gamma$, $r_0^-(0)=r_0(0)=p(0)$,
  and for each $\alpha<\delta<\lambda$ we have $(\forall\xi\in w_{\alpha+1})
  (r_\alpha(\xi)=r^-_\delta(\xi)=r_\delta(\xi))$ and $p\leq r_\alpha^-\leq
  r_\alpha\leq r^-_\delta\leq r_\delta$. 
\item[$(*)_2$] If $\xi\in\dom(r_\delta)\setminus w_\delta$, then 
\[\begin{array}{ll}
r_\delta\rest\xi\forces&\mbox{`` the sequence }\langle r^-_\alpha(\xi),
r_\alpha(\xi):\alpha\leq\delta\rangle\mbox{ is a legal partial play of }\\
&\quad\Game_0^\lambda\big(\name{\bbQ}_\xi,\name{\emptyset}_{\name{\bbQ}_\xi} 
\big)\mbox{ in which Complete follows }\name{\st}^0_\xi\mbox{ ''} 
\end{array}\]
and if $\xi\in w_{\delta+1}\setminus w_\delta$, then $\name{\st}_\xi$ is a
$\bbP_\xi$--name  for a winning strategy of Generic in
$\masgame(r_\delta(\xi),\name{\bbQ}_\xi)$ such that if $\langle p^\alpha_t:
t\in I_\alpha\rangle$ is given by that strategy to Generic at stage
$\alpha$, then $I_\alpha$ is an ordinal below $\lambda$. Also $\st_0$ is 
a suitable winning strategy of Generic in $\masgame(p(0),\bbQ_0)$.   
\item[$(*)_3$] $\cT_\delta=(T_\delta,\rk_\delta)$ is a standard $(w_\delta,
1)^\gamma$--tree, $|T_\delta|<\lambda$.   
\item[$(*)_4$] $\bar{p}^\delta=\langle p^\delta_t:t\in T_\delta\rangle$ and
  $\bar{q}^\delta=\langle q^\delta_t:t\in T_\delta\rangle$ are standard
  trees of conditions in $\bar{\bbQ}$, $\bar{p}^\delta\leq\bar{q}^\delta$. 
\item[$(*)_5$] If $t\in T_\delta$, $\rk_\delta(t)=\gamma$, then the
  condition $q^\delta_t$ decides the values of all names $\langle
  \name{\tau}_\alpha:\alpha\leq\delta\rangle$.  
\item[$(*)_6$] For $t\in T_\delta$ we have $\big(\dom(p)\cup
  \bigcup\limits_{\alpha<\delta}\dom(r_\alpha)\cup w_\delta\big) \cap
  \rk_\delta(t)\subseteq\dom(p^\delta_t)$ and for each $\xi\in
  \dom(p^\delta_t)\setminus w_\delta$: 
\[\begin{array}{ll}
p^\delta_t\rest\xi\forces_{\bbP_\xi}&\mbox{`` if the set } \{r_\alpha(\xi):
\alpha<\delta\}\cup\{p(\xi)\}\mbox{ has an upper bound in
}\name{\bbQ}_\xi,\\  
&\mbox{\quad then $p^\delta_t(\xi)$ is such an upper bound ''.} 
  \end{array}\]
\item[$(*)_7$] If $\xi\in w_\delta$, then $\name{\vare}_{\delta,\xi}$ is a
  $\bbP_\xi$--name for an ordinal below $\lambda$,
  $\name{\bar{p}}_{\delta,\xi},\name{\bar{q}}_{\delta,\xi}$ are
  $\bbP_\xi$--names for $\name{\vare}_{\delta,\xi}$--sequences of conditions
  in $\name{\bbQ}_\xi$.    
\item[$(*)_8$] If $\xi\in w_{\delta+1}\setminus w_\delta$, then 
\[\begin{array}{r}
\forces_{\bbP_\xi}\mbox{`` }\langle\name{\vare}_{\alpha,\xi},
\name{\bar{p}}_{\alpha,\xi},\name{\bar{q}}_{\alpha,\xi}:\alpha<\lambda
\rangle\mbox{ is a play of }\masgame(r_\delta(\xi),\name{\bbQ}_\xi)\\   
\mbox{ in which Generic uses $\name{\st}_\xi$ ''.}
  \end{array}\]
\item[$(*)_9$] If $t\in T_\delta$, $\rk_\delta(t)=\xi<\gamma$, then the
  condition $p^\delta_t$ decides the value of $\name{\vare}_{\delta, \xi}$,
  say $p^\delta_t\forces$``$\name{\vare}_{\delta,\xi}=
  \vare^t_{\delta,\xi}$'', and $\{(s)_\xi:t\vtl s\in T_\delta\}=
  \vare^t_{\delta,\xi}$ and   
\[q^\delta_t\forces_{\bbP_\xi}\mbox{`` } \name{\bar{p}}_{\delta,\xi}
(\vare)\leq p^\delta_{t\conc\langle\vare\rangle}(\xi)\mbox{ and }
\name{\bar{q}}_{\delta,\xi}(\vare)=q^\delta_{t\conc\langle\vare\rangle}(\xi)
\mbox{ for }\vare<\vare^t_{\delta,\xi}\mbox{ ''.}\]
\item[$(*)_{10}$] If $t_0,t_1\in T_\delta$, $\rk_\delta(t_0)=
  \rk_\delta(t_1)$ and $\xi\in w_\delta\cap\rk_\delta(t_0)$, $t_0\rest
  \xi=t_1\rest\xi$ but $\big(t_0\big)_\xi\neq \big(t_1\big)_\xi$, then  
\[p^\delta_{t_0\rest\xi}\forces_{\bbP_\xi}\mbox{`` the conditions  
$p^\delta_{t_0}(\xi),p^\delta_{t_1}(\xi)$ are incompatible ''.}\] 
\item[$(*)_{11}$] $\dom(r^-_\delta)=\dom(r_\delta)= \bigcup\limits_{t\in
    T_\delta}\dom(q^\delta_t)\cup\dom(p)$ and if $t\in T_\delta$, $\xi\in
\dom(r_\delta)\cap \rk_\delta(t)\setminus w_\delta$, and $q^\delta_t\rest
\xi\leq q\in\bbP_\xi$, $r_\delta\rest\xi\leq q$, then    
\[\begin{array}{ll}
q\forces_{\bbP_\xi}&\mbox{`` if the set }\{r_\alpha(\xi):\alpha<\delta\}
\cup\{q^\delta_t(\xi), p(\xi)\}\mbox{ has an upper bound in }  
\name{\bbQ}_\xi,\\ 
&\mbox{\quad then $r_\delta^-(\xi)$ is such an upper bound ''.}
  \end{array}\]
\end{enumerate}
After the construction is carried out we define a condition
$r\in\bbP_\gamma$ as follows. We let $\dom(r)=N\cap\gamma$ and for
$\xi\in\dom(r)$ we let $r(\xi)$ be a $\bbP_\xi$--name for a condition in
$\name{\bbQ}_\xi$ such that if $\xi\in w_{\alpha+1}\setminus w_\alpha$,
$\alpha<\lambda$ (or $\xi=0=\alpha$), then
\[\begin{array}{r}
\forces_{\bbP_\xi}\mbox{``}r(\xi)\geq r_\alpha(\xi)\mbox{ is an
  $(S,D)$--knighting condition for the servant }\ \ \\
\name{\bar{q}}^\xi\stackrel{\rm def}{=}\langle \name{\bar{q}}_{\delta,\xi}
(\vare):\delta\in S\ \&\ \vare<\name{\vare}_{\delta,\xi}\rangle\mbox{ ''.}  
\end{array}\] 
Clearly $r$ is well defined (remember $(*)_8$). Note also that $r_\delta\leq
r$ for all $\delta<\lambda$ and $p\leq r$. We will argue that $r$ is an
$(N,\bbP_\gamma)$--generic condition. To this end suppose towards
contradiction that $r^*\geq r$, $\alpha^*<\lambda$ and $r^*\forces
\name{\tau}_{\alpha^*}\notin N$. 

For each $\xi\in N\cap\gamma$ fix a $\bbP_\xi$--name $\name{\st}^*_\xi$ for
a winning strategy of COM in the game $\sergame(\name{\bar{q}}^\xi,r(\xi),
\name{\bbQ}_\xi)$.   Moreover, for each $\xi<\gamma$ fix a $\bbP_\xi$--name 
$\name{\st}^0_\xi$ for a winning strategy of Complete in $\Game^\lambda_0(
\name{\bbQ}_\xi,\name{\emptyset}_{\name{\bbQ}_\xi})$ such that it instructs
Complete to play $\name{\emptyset}_{\name{\bbQ}_\xi}$ as long as her
opponent plays $\name{\emptyset}_{\name{\bbQ}_\xi}$.

By induction on $\delta<\lambda$ we will build a sequence 
\[\big\langle r^*_\delta,r^+_\delta,\langle \name{A}^\xi_{\delta,i},
A^\xi_{\delta,i}: i<\lambda\ \&\ \xi\in N\cap\gamma\rangle:\delta< \lambda
\big\rangle\] 
such that the following demands $(*)_{12}$--$(*)_{16}$ are satisfied:
\begin{enumerate}
\item[$(*)_{12}$] $r^*_\alpha\in\bbP_\gamma$, $r^*\leq r^*_\alpha\leq
  r^+_\alpha \leq r^*_\delta$ for $\alpha<\delta<\lambda$,   
\item[$(*)_{13}$] $\name{A}^\xi_{\delta,i}$ is a $\bbP_\xi$--name for an
  element of $D\cap\bV$ (for $\xi\in N\cap\gamma$, $i<\lambda$), 
\item[$(*)_{14}$] if $\delta\in\lambda\setminus S$ and $\xi\in w_\delta$,
  then $r^*_\delta\rest\xi\forces_{\bbP_\xi}(\forall\alpha<\delta) (\forall
  i<\delta)(\name{A}^\xi_{\alpha,i}=A^\xi_{\alpha,i})$,   
\item[$(*)_{15}$] if $\beta<\delta<\lambda$ and $\xi\in w_{\beta+1}\setminus
  w_\beta$, then for some $\bbP_\xi$--names $\langle \name{s}_\alpha^\xi:
  \alpha\leq\beta\rangle$ we have    
\[\begin{array}{ll}
r^*_\delta\rest\xi\forces&\mbox{`` the sequence }
\langle \name{s}_\alpha^\xi,\mathop{\triangle}\limits_{i< \lambda} 
\name{A}^\xi_{\alpha,i}:\alpha\leq\beta\rangle\conc\langle r^*_\alpha(\xi), 
\mathop{\triangle}\limits_{i< \lambda}\name{A}^\xi_{\alpha,i}:\beta<\alpha
\leq \delta\rangle\\ 
&\ \mbox{ is a legal partial play of }\sergame(\name{\bar{q}}^\xi,r(\xi),
\name{\bbQ}_\xi)\\ 
&\ \mbox{ in which Generic follows } \name{\st}^*_\xi\mbox{ '',}  
\end{array}\] 
\item[$(*)_{16}$] $\dom(r^+_\delta)=\dom(r^*_\delta)$, $r^+_\delta\rest
  w_\delta= r^*_\delta\rest w_\delta$ and for each $\xi\in\dom(r^+_\delta)
  \setminus w_\delta$ we have  
\[\begin{array}{ll}
r^+_\delta\rest\xi\forces&\mbox{`` the sequence }\langle r^*_\alpha(\xi),
r^+_\alpha(\xi):\alpha\leq\delta\rangle\mbox{ is a legal partial play of }\\
&\quad\Game_0^\lambda\big(\name{\bbQ}_\xi,\name{\emptyset}_{\name{\bbQ}_\xi} 
\big)\mbox{ in which Complete follows }\name{\st}^0_\xi\mbox{ ''.}
\end{array}\]  
\end{enumerate}
So suppose that we have arrived to a stage $\delta<\lambda$ of the
construction and 
\begin{itemize}
\item $r^*_\alpha,r^+_\alpha$ for $\alpha<\delta$,
\item $\name{A}^\xi_{\alpha,i}$ for $i<\lambda$, $\alpha<\delta$ and $\xi\in
  \bigcup\limits_{\beta<\delta} w_\beta$,
\item $A^\xi_{\alpha,i}$ for $\alpha,i<\sup(\delta\setminus S)$ and $\xi\in
  \bigcup\limits_{\beta<\sup(\delta\setminus S)} w_\beta$
\end{itemize}
have been determined. 

\noindent {\bf Case 1:}\quad $\delta\notin S$.\\
Note that by our assumption on $S$ (in \ref{consec6}), $\delta$ is not a
successor ordinal, so $w_\delta=\bigcup\limits_{\alpha<\delta} w_\alpha$ (or
$\delta=0$ and $w_0=\{0\}$).  By $(*)_{16}+(*)_{15}$ we may choose a
condition $r^*_\delta$ stronger than all $r^+_\alpha$ (for $\alpha<\delta$)
and stronger than $r^*$ and such that for each $\xi\in w_\delta$ 
\begin{itemize}
\item if $\alpha<\delta$ and $i<\delta$, then $r^*_\delta \rest\xi$ forces a
  value to $\name{A}^\xi_{\alpha,i}$, say
\[r^*_\delta\rest\xi\forces_{\bbP_\xi}\name{A}^\xi_{\alpha,i}=
A^\xi_{\alpha,i}.\]  
\end{itemize}
For $\xi\in w_\delta$ and $i<\lambda$ we also let $\name{A}^\xi_{\delta,i}$
be a $\bbP_\xi$--name for the interval $(\delta,\lambda)$. The condition
$r^+_\delta$ is fully determined by $(*)_{16}$.

\noindent {\bf Case 2:}\quad $\delta\in S$ is a successor ordinal, say
$\delta=\beta+1$.\\ 
First, for each $\xi\in w_\delta\setminus w_\beta$ we pick $\bbP_\xi$--names
$\name{s}^\xi_\alpha$ and $\name{A}^\xi_{\alpha,i}$ (for $\alpha<\delta$,
$i<\lambda$) such that $r^+_\beta\rest\xi\forces_{\bbP_\xi}
r^+_\beta(\xi)=\name{s}^\xi_0$ and  
\[\begin{array}{ll}
r^*_\delta\rest\xi\forces&\mbox{`` the sequence }
\langle \name{s}^\xi_\alpha,\mathop{\triangle}\limits_{i< \lambda}
\name{A}^\xi_{\alpha,i}:\alpha\leq\beta\rangle\mbox{ is a legal partial play
  of }\\
&\quad \sergame(\name{\bar{q}}^\xi,r(\xi), \name{\bbQ}_\xi)\mbox{ in which
  Generic follows } \name{\st}^*_\xi\mbox{ ''.}   
\end{array}\] 
Next, we let $\dom(r^*_\delta)=\dom(r^+_\beta)$ and for each $\xi\in
w_\delta$ we choose $\bbP_\xi$--names $r^*_\delta(\xi)$ and
$\name{A}^\xi_{\delta,i}$ (for $i<\lambda$) such that  
\[r^+_\beta\rest\xi\forces_{\bbP_\xi} r^*_\delta(\xi),
\mathop{\triangle}\limits_{i< \lambda}\name{A}^\xi_{\delta,i}\mbox{ is the
  answer to the partial game as in $(*)_{15}$ given by $\st^*_\xi$.}\] 
For $\xi\in\dom(r^*_\delta)\setminus w_\delta$ we put $r^*_\delta(\xi)=
r^+_\beta(\xi)$. Then we define condition $r^+_\delta$ by $(*)_{16}$.  

\noindent {\bf Case 3:}\quad $\delta\in S$ is a limit ordinal.\\ 
We let $\dom(r^*_\delta)=\bigcup\limits_{\alpha<\delta}\dom(r^+_\alpha)$ and
by induction on $\xi\in\dom(r^*_\delta)$ we define $r^*_\delta(\xi)$ so that  
\begin{itemize}
\item if $\xi\notin w_\delta$, then $r^*_\delta\rest\xi\forces (\forall
  \alpha<\delta)(r^+_\alpha(\xi)\leq r^*_\delta(\xi))$ (exists by
  $(*)_{16}$), 
\item if $\xi\in w_\delta$ then for some $\bbP_\xi$--names
  $\name{A}^\xi_{\delta,i}$ for members of $D\cap\bV$
\[r^*_\delta\rest\xi\forces_{\bbP_\xi} r^*_\delta(\xi),
\mathop{\triangle}\limits_{i< \lambda}\name{A}^\xi_{\delta,i}\mbox{ is the
  answer to the partial game as in $(*)_{15}$ given by $\st^*_\xi$.}\]
\end{itemize}
The condition $r^+_\delta$ is given by $(*)_{16}$.
\medskip

After the above construction is carried out we note that 
\[\big\{\delta<\lambda:(\forall\xi\in w_\delta)(\forall\alpha,i<\delta) 
(\delta\in A^\xi_{\alpha,i})\big\}\in D,\]
so we may choose an ordinal $\delta\in S\setminus (\alpha^*+1)$ which is
a limit of points from $\lambda\setminus S$ and such that $\delta\in
\bigcap\limits_{\alpha<\delta}\mathop{\triangle}\limits_{i< \lambda}
A^\xi_{\alpha,i}$ for all $\xi\in w_\delta$.  The following claim provides
the desired contradiction (remember $(*)_0+(*)_5$). 

\begin{claim}
\label{cl5}
For some $t\in T_\delta$ such that $\rk_\delta(t)=\gamma$ the conditions
$q^\delta_t$ and $r^*_\delta$ are compatible
\end{claim}

\begin{proof}[Proof of the Claim]
The proof is very much like that of Claim \ref{cl3}. Let $\langle
\vare_\beta: \beta\leq\beta^*\rangle=w_\delta\cup\{\gamma\}$ be 
the increasing enumeration. For each $\xi<\gamma$ fix a $\bbP_\xi$--name  
$\name{\st}^*_\xi$ for a winning strategy of Complete in $\Game^\lambda_0(
\name{\bbQ}_\xi,\name{\emptyset}_{\name{\bbQ}_\xi})$ such that it instructs
Complete to play $\name{\emptyset}_{\name{\bbQ}_\xi}$ as long as her
opponent plays $\name{\emptyset}_{\name{\bbQ}_\xi}$. 

By induction on $\beta\leq\beta^*$ we will choose conditions $s_\beta,
s^*_\beta\in\bbP_{\vare_\beta}$ and $t=\langle (t)_{\vare_\beta}:\beta<
\beta^*\rangle\in T_\delta$ such that letting $t^\beta_\circ=\langle
(t)_{\vare_\beta'}:\beta'<\beta\rangle\in T_\delta$ we have
\begin{enumerate}
\item[$(\boxdot)_a$] $q^\delta_{t^\beta_\circ}\leq s_\beta$ and 
$r^*_\delta\rest\vare_\beta\leq s_\beta$,  
\item[$(\boxdot)_b$] $\dom(s_\beta)=\dom(s^*_\beta)$ and for every $\zeta<
  \vare_\beta$,  
\[\begin{array}{ll}
s^*_\beta\rest\zeta\forces_{\bbP_\zeta}&\mbox{`` }\langle s_{\beta'}(\zeta),
s^*_{\beta'}(\zeta):\beta'<\beta\rangle\mbox{ is a partial legal play of
}\Game^\lambda_0(\name{\bbQ}_\zeta,\name{\emptyset}_{\name{\bbQ}_\zeta})\\ 
&\ \ \mbox{ in which Complete uses her winning strategy $\name{\st}^*_\zeta$   
  ''.}\end{array}\]   
\end{enumerate}
Suppose that $\beta\leq\beta^*$ is a limit ordinal and we have already
defined $t^\beta_\circ=\langle (t)_{\vare_{\beta'}}: \beta'<\beta\rangle$ and 
$\langle s_{\beta'}, s^*_{\beta'}: \beta'<\beta\rangle$. Let $\xi=\sup(
\vare_{\beta'}:\beta'<\beta)$. It follows  from $(\boxdot)_b$ that we may
find a condition $s_\beta\in\bbP_{\vare_\beta}$ such that $s_\beta\rest\xi$
is stronger than all $s^*_{\beta'}$ (for $\beta'<\beta$) and $s_\beta\rest
[\xi, \vare_\beta)=r^*_\delta\rest [\xi,\vare_\beta)$.  Clearly $r^*_\delta
\rest \vare_\beta\leq s_\beta$ and also $q^\delta_{t^\beta_\circ}\rest\xi
\leq s_\beta \rest\xi$ (remember $(\boxdot)_a$). Now by induction on
$\zeta\in [\xi,\vare_\beta]$ we argue that $q^\delta_{t^\beta_\circ}\rest
\zeta \leq s_\beta\rest\zeta$. Suppose that $\xi\leq \zeta<\vare_\beta$ and
we know $q^\delta_{t^\beta_\circ}\rest\zeta\leq s_\beta\rest\zeta$. By
$(*)_2+(*)_4+(*)_6$ we know that $s_\beta\rest\zeta\forces (\forall 
i<\delta)(r_i(\zeta)\leq p^\delta_{t^\beta_\circ}(\zeta)\leq
q^\delta_{t^\beta_\circ}(\zeta))$ and therefore we may use $(*)_{11}$ to
conclude that    
\[s_\beta\rest\zeta\forces_{\bbP_\zeta} q^\delta_{t^\beta_\circ}(\zeta)\leq
r_\delta(\zeta)\leq r(\zeta)\leq r^*_\delta(\zeta)=s_\beta(\zeta).\]   
Then the condition $s^*_\beta\in\bbP_{\vare_\beta}$ is determined
$(\boxdot)_b$.  

Now suppose that $\beta=\beta'+1\leq\beta^*$ and we have already defined 
$s_{\beta'}, s^*_{\beta'}\in\bbP_{\vare_{\beta'}}$ and $t^{\beta'}_\circ\in
T_\delta$. It follows from the choice of $\delta$ and $(*)_{14}$ that
$r^*_\delta\rest \vare_{\beta'}\forces\delta\in\bigcap\limits_{\alpha
  <\delta} \mathop{\triangle}\limits_{i< \lambda} \name{A}^\xi_{\alpha,i}$
and hence, by the choice of $r$ and $(*)_{15}$ we have (remember $(\odot)$
of \ref{servant}(2)) 
\[s^*_{\beta'}\forces_{\bbP_{\vare_{\beta'}}}\mbox{`` }\big(\exists\vare<
\name{\vare}_{\delta, \vare_{\beta'}}\big)\big(\name{\bar{q}}_{\delta,
  \vare_{\beta'}}(\vare)\leq r^*_\delta(\vare_{\beta'})\big)\mbox{ ''.}\]
Therefore we may use $(*)_9$ to choose $\vare=(t)_{\vare_{\beta'}}$ and a
condition $s_\beta\in\bbP_{\vare_\beta}$ such that 
\begin{itemize}
\item $t^\beta_\circ\stackrel{\rm def}{=}t^{\beta'}_\circ\cup\{(
  \vare_{\beta'},\vare)\}\in T_\delta$, $s^*_{\beta'} \leq s_\beta
  \rest\vare_{\beta'}$ and    
\[s_\beta\rest\vare_{\beta'}\forces_{\bbP_{\vare_{\beta'}}}\mbox{`` }
q^\delta_{t^\beta_\circ}(\vare_{\beta'})\leq r^*_\delta(\vare_{\beta'})=
s_\beta(\vare_{\beta'})\mbox{ '',}\]
\item $r^*_\delta\rest (\vare_{\beta'},\vare_\beta)=s_\beta\rest
  (\vare_{\beta'},\vare_\beta)$. 
\end{itemize}
We finish exactly like in the limit case.

After the inductive construction is completed, look at $t=t^{\beta^*}_\circ$ 
and $s_{\beta^*}$.  
\end{proof}

\noindent (b)\quad Should be clear at the moment.
\end{proof}

\begin{definition}
[See {\cite[Def. 3.1]{RoSh:860}}]
Let $\bbQ$ be a strategically $({<}\lambda)$--complete forcing notion. 
\begin{enumerate}
\item Let $p\in\bbQ$. A game $\rcbgame(p,\bbQ)$ is defined similarly to
  $\bgame(p,\bbQ)$ (see \ref{Bplus}) except that the winning criterion
  $(\circledast)^\bp_{\rm rbB}$ is weakened to 
\begin{enumerate}
\item[$(\circledast)^{\rm rc}_B$] there is a condition $p^*\in\bbQ$
  stronger than $p$ and such that 
\[p^*\forces_{\bbQ}\mbox{`` }\big\{\alpha<\lambda:\big(\exists t\in
I_\alpha\big)\big(q^\alpha_t \in\Gamma_\bbQ\big)\big\}\in D^\bbQ\mbox{
  ''.}\]  
\end{enumerate}
\item A forcing notion $\bbQ$ is {\em reasonably B--bounding over $D$\/}
  if for any $p\in\bbQ$, Generic has a winning strategy in the game
  $\rcbgame(p,\bbQ)$.      
\end{enumerate}
\end{definition}

\begin{observation}
\label{obs3.2}
  If $\bbQ$ is reasonably B--bounding over $D$, then it is reasonably merry
  over $(S,D)$.
\end{observation}

It is not clear though, if forcing notions which are reasonably B--bounding 
over a $\Dl$--parameter $\bp$ are also reasonably merry (see Problem
\ref{prob3}). Also, we do not know if fuzzy properness introduced in
\cite[\S A.3]{RoSh:777} implies that the considered forcing notion is
reasonably merry (see Problem \ref{prob4}), even though the former property
seems to be almost built into the latter one. 

One may ask if being reasonably merry implies being B--bounding. There are
examples that this is not the case. The forcing notion $\tefo$ (see
\ref{3.2} below) was introduced in \cite[Section 6]{RoSh:860} and by \cite[Proposition
6.4]{RoSh:860} we know that it is not reasonably B--bounding over
$D$. However we will see in \ref{tefois merry} that it is reasonably merry
over $(S,D)$.   

\begin{definition}
[See {\cite[Def. 5.1]{RoSh:860}}]
\label{ext}
\begin{enumerate}
\item Let $\alpha<\beta<\lambda$. {\em An $(\alpha,\beta)$--extending
function\/} is a mapping $c:\cP(\alpha)\longrightarrow\cP(\beta)\setminus 
\cP(\alpha)$ such that $c(u)\cap\alpha=u$ for all $u\in\cP(\alpha)$. 
\item Let $C$ be an unbounded subset of $\lambda$. {\em A $C$--extending
sequence\/} is a sequence $\gc=\langle c_\alpha:\alpha\in C\rangle$ such
that each $c_\alpha$ is an $(\alpha,\min(C\setminus(\alpha+1)))$--extending
function.  
\item Let $C\subseteq\lambda$, $|C|=\lambda$, $\beta\in C$, $w\subseteq
\beta$ and let $\gc=\langle c_\alpha:\alpha\in C\rangle$ be a $C$--extending
sequence. We define $\pos^+(w,\gc,\beta)$ as the family of all subsets $u$
of $\beta$ such that  
\begin{enumerate}
\item[(i)]  if $\alpha_0=\min\big(\{\alpha\in C:(\forall\xi\in w)(\xi<
\alpha)\}\big)$, then $u\cap\alpha_0=w$ (so if $\alpha_0=\beta$, then
$u=w$), and  
\item[(ii)] if $\alpha_0,\alpha_1\in C$, $w\subseteq\alpha_0<\alpha_1=
\min(C\setminus(\alpha_0+1))\leq\beta$, then either $c_{\alpha_0}(u\cap
\alpha_0)=u\cap\alpha_1$ or $u\cap\alpha_0=u\cap\alpha_1$, 
\item[(iii)] if $\sup(w)<\alpha_0=\sup(C\cap\alpha_0)\notin C$,
$\alpha_1=\min\big(C\setminus(\alpha_0+1)\big)\leq\beta$, then
$u\cap\alpha_1=u\cap\alpha_0$.
\end{enumerate}
For $\alpha_0\in \beta\cap C$ such that $w\subseteq\alpha_0$, the family
$\pos(w,\gc,\alpha_0,\beta)$ consists of all elements $u$ of $\pos^+(w,\gc,\beta)$
which satisfy also the following condition: 
\begin{enumerate}
\item[(iv)] if $\alpha_1=\min\big(C\setminus(\alpha_0+1)\big)\leq\beta$,
  then $u \cap\alpha_1=c_{\alpha_0}(w)$.
\end{enumerate}
\end{enumerate}
\end{definition}

\begin{definition}
[See {\cite[Def. 6.2]{RoSh:860}}]
\label{3.2}
We define a forcing notion $\tefo$ as follows.\\
{\bf A condition in $\tefo$} is a triple $p=(w^p,C^p,\gc^p)$ such that 
\begin{enumerate}
\item[(i)]    $C^p\in D$, $w^p\subseteq\min(C^p)$,  
\item[(ii)]   $\gc^p=\langle c^p_\alpha:\alpha\in C^p\rangle$ is a
  $C^p$--extending sequence.
\end{enumerate}
{\bf The order $\leq_{\tefo}=\leq$ of $\tefo$} is given by 

$p\leq_{\tefo} q$\qquad if and only if 
\begin{enumerate}
\item[(a)] $C^q\subseteq C^p$ and $w^q\in\pos^+(w^p,\gc^p,\min(C^q))$ and  
\item[(b)] if $\alpha_0,\alpha_1\in C^q$, $\alpha_0<\alpha_1=\min(C^q
\setminus(\alpha_0+1))$ and $u\in\pos^+(w^q,\gc^q,\alpha_0)$, then
$c^q_{\alpha_0}(u)\in \pos(u,\gc^p,\alpha_0,\alpha_1)$. 
\end{enumerate}
For $p\in\tefo$, $\alpha\in C^p$ and $u\in\pos^+(w^p,\gc^p,\alpha)$ we let  
$p\rest_\alpha u\stackrel{\rm def}{=}(u,C^p\setminus\alpha,\gc^p\rest(C^p
\setminus\alpha))$.
\end{definition}

In \cite[Problem 6.1]{RoSh:860} we asked if $\lambda$--support iterations
of forcing notions $\tefo$ are $\lambda$--proper. Now we may answer
this question positively (assuming that $\lambda$ is strongly
inaccessible). First, let us state some auxiliary definitions and facts. 

\begin{proposition}
\label{3.4}
\begin{enumerate}
\item $\tefo$ is a $({<}\lambda)$--complete forcing notion of cardinality 
  $2^\lambda$. 
\item If $p\in\tefo$ and $\alpha\in C^p$, then
\begin{itemize}
\item for each $u\in\pos^+(w^p,\gc^p,\alpha)$, $p\rest_\alpha u\in\tefo$ is
  a condition stronger than $p$, and    
\item the family $\{p\rest_\alpha u: u\in\pos^+(w^p,\gc^p,\alpha)\}$ is
  pre-dense above $p$. 
\end{itemize}
\item Let $p\in\tefo$ and $\alpha<\beta$ be two successive members of
$C^p$. Suppose that for each $u\in\pos^+(w^p,\gc^p,\alpha)$ we are given a
condition $q_u\in\tefo$ such that $p\rest_\beta c^p_\alpha(u)\leq q_u$. 
Then there is a condition $q\in\tefo$ such that letting $\alpha'=\min(
C^q\setminus\beta)$ we have
\begin{enumerate}
\item[(a)] $p\leq q$, $w^q=w^p$, $C^q\cap\beta=C^p\cap\beta$ and
  $c^q_\delta=c^p_\delta$ for $\delta\in C^q\cap\alpha$, and 
\item[(b)] $\bigcup\big\{w^{q_u}:u\in\pos^+(w^p,\gc^p,\alpha)\big\}\subseteq 
  \alpha'$, and 
\item[(c)] $q_u\leq q\rest_{\alpha'} c^q_\alpha(u)$ for every
  $u\in\pos^+(w^p,\gc^p,\alpha)$. 
\end{enumerate}
\item Assume that $p\in\tefo$, $\alpha\in C^p$ and $\name{\tau}$ is a
$\tefo$--name such that $p\forces$``$\name{\tau}\in\bV$''. Then there is a
  condition $q\in\tefo$ stronger than $p$ and such that   
\begin{enumerate}
\item[(a)] $w^q=w^p$, $\alpha\in C^q$ and $C^q\cap\alpha=C^p\cap\alpha$,
  and 
\item[(b)] if $u\in\pos^+(w^q,\gc^q,\alpha)$ and $\gamma=\min(C^q\setminus 
(\alpha+1))$, then the condition $q\rest_\gamma c^q(u)$ forces a value to
$\name{\tau}$. 
\end{enumerate}
\end{enumerate}
\end{proposition}

\begin{proof}  
Fully parallel to \cite[Proposition 5.1]{RoSh:860}.
\end{proof}

\begin{definition}
{\em The natural limit\/} of an $\leq_{\tefo}$--increasing sequence
$\bar{p}=\langle p_\xi:\xi<\gamma\rangle\subseteq\tefo$ (where
$\gamma<\lambda$ is a limit ordinal) is the condition $q=(w^q,C^q,\gc^q)$ 
defined as follows: 
\begin{itemize}
\item $w^q=\bigcup\limits_{\xi<\gamma}w^{p_\xi}$, $C^q=\bigcap\limits_{\xi<
\gamma} C^{p_\xi}$ and 
\item $\gc^q=\langle c^q_\delta:\delta\in C^q\rangle$ is such that for
  $\delta\in C^q$ and $u\subseteq\delta$ we have $c^q_\delta(u)=
  \bigcup\limits_{\xi<\gamma}c^{p_\xi}_\delta(u)$. 
\end{itemize}
\end{definition}

\begin{proposition}
\label{limit}
\begin{enumerate}
\item Suppose $\bar{p}=\langle p_\xi:\xi<\lambda\rangle$ is a
$\leq_{\tefo}$--increasing sequence of conditions from $\tefo$ such that 
\begin{enumerate}
\item[(a)] $w^{p_\xi}=w^{p_0}$ for all $\xi<\lambda$, and 
\item[(b)] if $\gamma<\lambda$ is limit, then $p_\gamma$ is the natural
  limit of $\bar{p}\rest\gamma$, and 
\item[(c)] for each $\xi<\lambda$, if $\delta\in C^{p_\xi}$, $\otp(C^{p_\xi}
\cap\delta)=\xi$, then $C^{p_{\xi+1}}\cap (\delta+1)=C^{p_\xi}\cap
(\delta+1)$ and for every $\alpha\in C^{p_{\xi+1}}\cap\delta$ we have 
$c^{p_{\xi+1}}_\alpha=c^{p_\xi}_\alpha$. 
\end{enumerate}
Then the sequence $\bar{p}$ has an upper bound in $\tefo$.
\item Suppose that $p\in\tefo$ and $\name{h}$ is a $\tefo$--name such that
$p\forces$``$\name{h}:\lambda\longrightarrow\bV$''. Then there is a
condition $q\in\tefo$ stronger than $p$ and such that    
\begin{enumerate}
\item[$(\otimes)$] if $\delta<\delta'$ are two successive points of $C^q$,
$u\in\pos(w^q,\gc^q,\delta)$, then the condition $q\rest_{\delta'}
c^q_\delta(u)$ decides the value of $\name{h}\rest(\delta+1)$. 
\end{enumerate}
\end{enumerate}
\end{proposition}

\begin{proof}
Fully parallel to \cite[Proposition 5.2]{RoSh:860}.
\end{proof}
 
\begin{proposition}
\label{tefois merry}
Assume that $\lambda,S,D$ are as in \ref{consec6}. The forcing notion
$\tefo$ is reasonably merry over $(S,D)$.
\end{proposition}

\begin{proof}
  Let $p\in \tefo$. We will describe a strategy $\st$ for Generic in the
  game $\masgame(p,\tefo)$ - this strategy is essentially the same as the
  one in the proof of \cite[Proposition 5.4]{RoSh:860}, only the argument
  that it is a winning strategy is different.

  In the course of a play the strategy $\st$ instructs Generic to build
  aside an increasing sequence of conditions $\bar{p}^*=\langle
  p^*_\alpha:\alpha< \lambda\rangle\subseteq\tefo$ such that for each
  $\alpha<\lambda$:
\begin{enumerate}
\item[(a)] $p_0^*=p$ and $w^{p^*_\alpha}=w^p$, and  
\item[(b)] if $\alpha<\lambda$ is limit, then $p^*_\alpha$ is the natural
  limit of $\bar{p}^*\rest\alpha$, and  
\item[(c)] if $\delta\in C^{p^*_\alpha}$, $\otp(C^{p^*_\alpha}\cap\delta)=
  \alpha$, then $C^{p^*_{\alpha+1}}\cap (\delta+1)=C^{p^*_\alpha}\cap
  (\delta+1)$ and for every $\xi\in C^{p^*_{\alpha+1}}\cap\delta$ we have
  $c^{p^*_{\alpha+1}}_\xi= c^{p^*_\alpha}_\xi$, and 
\item[(d)] after stage $\alpha$ of the play of $\masgame(p,\tefo)$, the
  condition $p^*_{\alpha+1}$ is determined.
\end{enumerate}
After arriving to the stage $\alpha$, Generic is instructed to pick
$\delta\in C^{p^*_\alpha}$ such that $\otp(C^{p^*_\alpha}\cap\delta)
=\alpha$, put $\gamma=\min(C^{p^*_\alpha}\setminus(\delta+1))$ and play as 
her innings of this stage:
\[I_\alpha=\pos^+(w^{p^*_\alpha},\gc^{p^*_\alpha},\delta)\quad\mbox{ and }
\quad p^\alpha_u=p^*_\alpha\rest_\gamma c^{p^*_\alpha}_\delta(u)\mbox{ for }
u\in I_\alpha.\]  
Then Antigeneric answers with $\langle q^\alpha_u:u\in I_\alpha\rangle
\subseteq\tefo$. Since $p^*_\alpha\rest_\gamma c^{p^*_\alpha}_\delta(u)\leq  
q^\alpha_u$ for each $u\in\pos^+(w^{p^*_\alpha},\gc^{p^*_\alpha},\delta)$,
Generic may use \ref{3.4}(3) (with $\delta,\gamma,p^*_\alpha,q^\alpha_u$ 
here standing for $\alpha,\beta,p,q_u$ there) to pick a condition
$p^*_{\alpha+1}$ such that, letting $\alpha'=\min(C^{p^*_{\alpha+1}}
\setminus\gamma)$, we have  
\begin{enumerate}
\item[(e)] $p^*_\alpha\leq p^*_{\alpha+1}$, $w^{p^*_{\alpha+1}}=w^p$, $C^{
p^*_{\alpha+1}}\cap\gamma=C^{p^*_\alpha}\cap\gamma$ and $c^{p^*_{\alpha+
1}}_\xi=c^{p^*_\alpha}_\xi$ for $\xi\in C^{p^*_{\alpha+1}}\cap\delta$, and 
\item[(f)] $\bigcup\big\{w^{q^\alpha_u}:u\in I_\alpha\big\}\subseteq
  \alpha'$, and  
\item[(g)] $q^\alpha_u\leq p^*_{\alpha+1}\rest_{\alpha'}
  c^{p^*_{\alpha+1}}_\delta(u)$ for every  $u\in I_\alpha$. 
\end{enumerate}
This completes the description of $\st$. Suppose that $\big\langle I_\alpha,
\langle p^\alpha_u,q^\alpha_u:u\in I_\alpha\rangle:\alpha<\lambda
\big\rangle$ is the result of a play of $\masgame(p,\tefo)$ in which Generic
followed $\st$ and constructed aside $\bar{p}^*=\langle p^*_\alpha:\alpha<
\lambda\rangle$. By \ref{limit}, there is a condition $p^*\in\tefo$ stronger
than all $p^*_\alpha$ (for $\alpha<\lambda$). We claim that $p^*$ is an 
$(S,D)$--knighting condition for the servant $\bar{q}=\langle q^\alpha_u:
\alpha\in S \ \&\ u\in I_\alpha\rangle$.  To this end consider the following
strategy $\st^*$ of COM in $\sergame(\bar{q},p^*,\tefo)$. After arriving to
a stage $\alpha\in S$ of a play of $\sergame(\bar{q},p^*,\tefo)$, when
$\langle r_\beta,A_\beta:\beta<\alpha\rangle$ has been already constructed,
COM plays as follows.\\
If $\alpha$ is a successor or $\alpha\notin \bigcap\limits_{\beta<\alpha}
A_\beta$, then she just puts $r_\alpha,A_\alpha$ such that:
\begin{enumerate}
\item[(h)] $r_\beta\leq r_\alpha$ for all $\beta<\alpha$ and if $\delta\in
  C^{p^*_\alpha}$ is such that $\otp(\delta\cap C^{p^*_\alpha})=\alpha$,
  then $w^{r_\alpha}\setminus (\delta+1)\neq\emptyset$, and 
\item[(i)] $A_\alpha=\bigcap\limits_{\beta<\alpha}A_\beta\cap
  \bigcap\limits_{\beta<\alpha} C^{r_\beta}\setminus
  (\sup(w^{r_\alpha})+1)$. 
\end{enumerate}
If $\alpha\in \bigcap\limits_{\beta<\alpha} A_\beta$ is a limit ordinal,
then COM first lets $u=\bigcup\limits_{\beta<\alpha} w^{r_\beta}$. It
follows from (h)+(i) from earlier stages that $u\subseteq\alpha$, $\alpha\in
C^{p^*_\alpha}$ and $\otp(\alpha\cap C^{p^*_\alpha})=\alpha$. Note that
$\alpha\in\bigcap\limits_{\beta<\alpha}C^{r_\beta}$, $u\in
\bigcap\limits_{\beta<\alpha}\pos^+(w^{r_\beta},\gc^{r_\beta},\alpha)$ and 
$u\in I_\alpha$. Let $\alpha'=\min\big(C^{p^*}\setminus (\alpha+1)\big)$ and
$\alpha''\in \bigcap\limits_{\beta<\alpha} C^{r_\beta}\setminus (\alpha+1)$.
It follows from (c)+(g) that for each $\beta<\alpha$ we have    
\[q^\alpha_u\leq p^*_{\alpha+1}\rest_{\alpha'} c^{p^*_{\alpha+1}}_\alpha(u) 
\leq p^*\rest_{\alpha'} c^{p^*}_\alpha(u) \leq r_\beta\rest_{\alpha''}
c^{r_\beta}_\alpha(u).\]   
Hence COM may choose a condition $r_\alpha\geq q^\alpha_u$ stronger than all
$r_\beta$ (for $\beta<\alpha$) and satisfying (h). Then $A_\alpha$ is given
by (i).

It follows directly from the description of $\st^*$ that it is a winning
strategy of COM in $\sergame(\bar{q},p^*,\tefo)$. 
\end{proof}

The master game $\masgame$ used to define the property of being reasonably
merry is essentially a variant of the A--reasonable boundedness game
$\Game^{\rm rcA}$ of \cite[Def. 3.1]{RoSh:860}. The related bounding
property was weakened in \cite[Def. 2.9]{RoSh:890} by introducing {\em double
  {\bf a}--reasonably completeness game\/} $\Game^{{\rm rc}{\bf 2a}}$. We
may use these ideas to introduce a property much weaker than ``reasonably
merry'', though the description of the resulting notions becomes somewhat
more complicated. As an award for additional complication we get a true
preservation theorem, however. In the rest of this section, in addition to
\ref{consec6} we assume also 

\begin{context}
 \label{extcon}
$\bar{\mu}=\langle\mu_\alpha:\alpha<\lambda\rangle$ is a sequence of
cardinals below $\lambda$ such that $(\forall \alpha<\lambda)(\aleph_0
\leq\mu_\alpha=\mu_\alpha^{|\alpha|})$. 
\end{context}

\begin{definition}
  \label{doubleser}
Let $\bbQ$ be a forcing notion.
\begin{enumerate}
\item {\em A double $\bbQ$--servant over $S,\bar{\mu}$\/} is a sequence 
\[\bar{q}=\langle \xi_\delta,q^\delta_\gamma:\delta\in S\ \&\ \gamma<
\mu_\delta\cdot \xi_\delta\rangle\]
such that for $\delta\in S$, 
\begin{itemize}
\item $0<\xi_\delta<\lambda$ and $q^\delta_\gamma\in\bbQ$
  (for $\gamma<\mu_\delta\cdot \xi_\delta$), 
\item $\big(\forall i,i'<\xi_\delta\big)\big(\forall j< \mu_\delta\big)
  \big(i'<i\ \Rightarrow\ q^\delta_{\mu_\delta\cdot i'+j}\leq
  q^\delta_{\mu_\delta\cdot i+j}\big)$.
\end{itemize}
(Here $\mu_\delta$ is treated as an ordinal and $\mu_\delta\cdot\xi_\delta$
is the ordinal product of $\mu_\delta$ and $\xi_\delta$.)
\item Let $\bar{q}$ be a double $\bbQ$--servant over $S,\bar{\mu}$ and let 
  $q\in\bbQ$. We define a game $\dsgame(\bar{q},q,\bbQ)$ as follows. A play of
  $\dsgame(\bar{q},q,\bbQ)$ lasts at most $\lambda$ steps during which the
  players, COM and INC, attempt to construct a sequence $\langle r_\alpha,
  A_\alpha:\alpha<\lambda\rangle$ such that 
  \begin{itemize}
\item $r_\alpha\in\bbQ$, $q\leq r_\alpha$, $A_\alpha\in D$ and $\alpha<\beta
  <\lambda\ \Rightarrow\ r_\alpha\leq r_\beta$.
  \end{itemize}
The terms $r_\alpha,A_\alpha$ are chosen successively by the two players so
that 
\begin{itemize}
\item if $\alpha\notin S$, then INC picks $r_\alpha,A_\alpha$, and 
\item if $\alpha\in S$, then COM chooses $r_\alpha,A_\alpha$.
\end{itemize}
If at some moment of the play one of the players has no legal move, then INC
wins; otherwise, if both players always had legal moves and the sequence
$\langle r_\alpha,A_\alpha:\alpha<\lambda\rangle$ has been constructed, then
COM wins if and only if 
\[(\heartsuit)\quad \big(\forall\delta\in S\big)\big([\delta\in
  \bigcap\limits_{\alpha<\delta} A_\alpha\ \&\ \delta\mbox{ is limit }]\
  \Rightarrow\ (\exists j<\mu_\delta)(\forall i<\xi_\delta)(
  q^\delta_{\mu_\delta\cdot i+j}\leq r_\delta)\big).\] 
\item If COM has a winning strategy in the game $\dsgame(\bar{q},q,\bbQ)$,
  then we will say that {\em $q$ is an knighting condition for the double 
    servant $\bar{q}$}. 
\end{enumerate}
\end{definition}

\begin{definition}
\label{doublemaster}
Let $\bbQ$ be a strategically $({<}\lambda)$--complete forcing notion.
\begin{enumerate}
\item For a condition $p\in\bbQ$ we define a game $\dmgame(p,\bbQ)$ between
two players, Generic and Antigeneric. A play of $\dmgame(p,\bbQ)$ lasts
$\lambda$ steps and during a play a sequence     
\[\big\langle\xi_\alpha,\langle p^\alpha_\gamma, q^\alpha_\gamma:
\gamma<\mu_\alpha\cdot \xi_\alpha\rangle:\alpha<\lambda\big\rangle.\]  
is constructed. (Again, $\mu_\alpha\cdot\xi_\alpha$ is the ordinal product
of $\mu_\alpha$ and $\xi_\alpha$.) Suppose that the players have arrived to
a stage $\alpha<\lambda$ of the game. First, Antigeneric picks a non-zero 
ordinal $\xi_\alpha<\lambda$. Then the two players start a subgame of length
$\mu_\alpha\cdot\xi_\alpha$ alternately choosing the terms of the sequence
$\langle p^\alpha_\gamma,q^\alpha_\gamma: \gamma<\mu_\alpha\cdot
\xi_\alpha\rangle$. At a stage $\gamma=\mu_\alpha\cdot i+j$ (where
$i<\xi_\alpha$, $j<\mu_\alpha$) of the subgame, first Generic picks a
condition $p^\alpha_\gamma\in\bbQ$ stronger than $p$ and stronger than all
conditions $q^\alpha_\delta$ for $\delta<\gamma$ of the form
$\delta=\mu_\alpha\cdot i'+j$ (where $i'<i$), and then Antigeneric answers
with a condition $q^\alpha_\gamma$ stronger than $p^\alpha_\gamma$.  

At the end, Generic wins the play if 
\begin{enumerate}
\item[(a)] there were always legal moves for both players (so a sequence
  \[\big\langle\xi_\alpha,\langle p^\alpha_\gamma,q^\alpha_\gamma:\gamma
  <\mu_\alpha\cdot\xi_\alpha\rangle: \alpha< \lambda\big\rangle\]
has been constructed) and  
\item[(b)] for each $\alpha\in S$, the conditions in
  $\{p^\alpha_j:j<\mu_\alpha\}$ are pairwise incompatible, and 
\item[(c)] letting 
\[\bar{q}=\langle \xi_\delta,q^\delta_\gamma:\delta\in S\ \&\ \gamma<
\mu_\delta\cdot \xi_\delta\rangle\]
(it is a double $\bbQ$--servant over $S$) we may find a knighting condition
$q\geq p$ for the double servant $\bar{q}$.   
\end{enumerate}
\item A forcing notion $\bbQ$ is {\em reasonably double merry over
    $(S,D,\bar{\mu})$\/} if (it is strategically $({<})$--complete and)
  Generic has a winning strategy in the game $\dmgame(p,\bbQ)$ for any
  $p\in\bbQ$. 
\end{enumerate}
\end{definition}

\begin{theorem}
 \label{presdouble}
Assume that $\lambda,S,D,\bar{\mu}$ are as in \ref{consec6}+\ref{extcon}. Let
$\bar{\bbQ}=\langle\bbP_\alpha,\name{\bbQ}_\alpha:\alpha<\gamma\rangle$ be a
$\lambda$--support iteration such that for each $\alpha<\gamma$: 
\[\forces_{\bbP_\alpha}\mbox{`` }\name{\bbQ}_\alpha\mbox{ is reasonably
  double merry over $(S,D,\bar{\mu})$''.}\]
Then $\bbP_\gamma=\lim(\bar{\bbQ})$ is reasonably double merry over
$(S,D,\bar{\mu})$ (so also $\lambda$--proper). 
\end{theorem}

\begin{proof}
Combine the proof of \cite[Thm 2.12]{RoSh:890} (the description of the
strategy here is the same as the one there) with the end of the proof of
\ref{masiter}(a). 
\end{proof}

\section{Open problems}

\begin{problem}
  \label{prob0}
Let $\bp=(\bar{P},S,D)$ be a $\Dl$--parameter. Does ``reasonably B-bounding
over $\bp$'' (see \ref{Bplus}) imply ``reasonably B-bounding over $D$'' (of
\cite[Def. 3.1]{RoSh:860})?
Does ``reasonably B-bounding over $\bp$'' imply ``B--noble over $\bp$'' ?
(Note \ref{obs3.2}.)
\end{problem}
\begin{problem}
\label{prob1}
  Are there any relations among the notions of properness over $D$--semi
  diamonds (of \cite{RoSh:655}), properness over $D$--diamonds (of
  \cite{Ei03}) and B--nobleness (of \ref{supergame})?
\end{problem}

\begin{problem}
\label{prob2}
Does $\lambda$--support iterations of forcing notions of the form $\fEE$
  add $\lambda$--Cohens? Here we may look at iterations as in
  \ref{nobleiteration} or \ref{lordsthm}. 
\end{problem}

\begin{problem}
\label{prob3}
Does ``reasonably B--bounding over $\bp$'' (for a $\Dl$--parameter $\bp$)
imply ``reasonably merry over $(S,D)$'' (for some $S,D$ as in \ref{consec6})?
\end{problem}

\begin{problem}
 \label{prob4}
Does ``fuzzy proper over quasi $D$--diamonds for $W$'' (see
\cite[Def. A.3.6]{RoSh:777}) imply ``reasonably merry''? (Any result in this
direction may require additional assumptions on $W,\bar{\gY}$ in
\cite[A.3.1, A.3.3]{RoSh:777}.)
\end{problem}

\bibliographystyle{amsalpha}

\end{document}